Alexander S. Yessenin-Volpin
Christer Hennix

# Beware Of The Gödel-Wette Paradox!

## Part I

### (Preliminary Version)

## ABSTRACT.

This paper gives a counter-example to the impossibility, by Gödel's second incompleteness theorem, of proving a formula expressing the consistency of arithmetic in a fragment of arithmetic *on the assumption that the latter be consistent*. This counter-example gives rise to a new type of metamathematical paradox, to be called the *Gödel-Wette paradox*, which E. Wette claims to have established since some time ago (see [Wette, 1971], [Wette, 1974]). Nevertheless, our work is independent of Wette's since we have failed to understand the details of his work although we recognize the possibility of the correctness of the latter. Furthermore, the Gödel-Wette paradox is not the only foundational anomaly which the framework of our approach has uncovered but new questions concerning the decision problem, completeness problem, truth definitions and the status of Richard's paradox in arithmetic and set theory (including type theory) have arisen as well. This work will, eventually, be unified into a single monograph.

**Acknowledgement**. This work was supported by a grant to the second author from the Clay Mathematics Institute, Cambrige, MA, USA, which is hereby gratefully acknowleged.

Alexander S. Yessenin-Volpin
Christer Hennix

Beware Of The Gödel-Wette Paradox!

Part I

(Preliminary Version)

Preface.

**A**. In the field of the foundations of mathematics there has evolved, since 1931, the widespread conviction that any consistent and 'sufficiently rich and smooth' formal system, S, is unavoidably incomplete — so that at least one closed formula (or 'sentence'), F, in the language, LS, of S is 'true' but not S-provable and, furthermore, that one such formula is $Con_S$ , i.e. a/the formula expressing, in LS, the *consistency* of the system S.

Statements of this sort belong to the metatheory, MS, of S, which was, since [GÖDEL, 1931], 'arithmeticized', i.e. imbedded, on the basis of an *enumeration* (or '*Nummerzuordnung*' in the terminology of [HILBERT-BERNAYS, 1939]), $En_S$ , into the usual arithmetic, Ar, such that each object E in LS becomes 'replaced' by an integer $_E$ whereby the study of LS and S is transformed into the arithmetic which is 'contained' in S provided the latter is 'sufficiently rich'. The properies and relations between the objects $E_1$ , $E_2$ , ... are thus 'replaced' by relations between the integers $_{E_1}$ , $_{E_2}$ ,... (which in this work shall be referred to as '$En_S$-*numbers*' or '$En_S$-N's' of these objects) and, thus, can be studied in terms of arithmetical functions. Although otherwise more or less arbitrary, $En_S$ is introduced in a way which makes this study convenient. The basic properties of the objects E — such as being a formula in LS or a deduction (proof) in S — are thus 'expressed' by primitive recursive (p.r.) predicates or by the characteristic p.r. function symbols which '$En_S$-*represent* ' these predicates.

As there is no 'canonical' choice of the enumeration $En_S$ for any fixed system S, one can consider several such enumerations $En_S^1$, $En_S^2$ , ... and interconnections between them. Nevertheless, when one speaks about Gödel's incompleteness theorems, (s)he has in mind one specific $En_S$ which can be, in case of need, easily adjusted to the current tasks. The main results, starting with Gödel's incompleteness theorems, being once proved for, say, $En_S^1$, continue to hold, on the basis of the 'same' proofs, for $En_S^2$. I.e, details of $En_S$'s rarely matter as the basic property of any such enumeration must consist in that when $_{E_1}= {}_{E_2}$, then the objects $E_1$, $E_2$ must coincide (i.e., must be the 'same'). The inverse 'mapping' is not always supposed 1-1: Rosser ([ROSSER, 1936]) has shown that $En_S$ can be allowed to contain repetitions, i.e., an object, E, can be allowed to have more than one En-N: $\overset{1}{E}$, $\overset{2}{E}$, ..., in which case the notation $_E$ requires additional specifications: say, it may stand for the least of $\overset{1}{E}$, $\overset{2}{E}$, ... — though it is not always easy, or even possible, to indicate this 'unique' integer for a given object E. (When only a single system, S, is at issue the subscript S in $En_S$ ($En_S$-N) will often be omitted). In addition, in the same work Rosser also proved that Gödel's proofs apply to any sufficiently rich system S if the $En_S$-numbers of all S-provable formulae form a set coinciding with the set of all values of a p.r. function; this condition is at issue when S is considered as 'sufficienly smooth'.

As is well known, Gödel's theorems survive Rosser's generalization. Gödel's first incompleteness theorem is deduced from the construction of a formula, G, which expresses, via $En_S$, that it is not provable in S — and, therefore, if it is provable, then it is false and also provable while if G is not provable in S it is true and not provable. If S is consistent and its theorems are true (i.e., S is sound), then only the second possibility can take place so that S contains a true and S-unprovable formula, G. (Assuming the -consistency of S, this theorem can be extended also for the formula ¬G — the negation of G. That is, also ¬G is not provable in S). That is the content of the first incompleteness theorem.

Gödel's second incompleteness theorem consists in that the particular formula, $Con_S$, which in a straightforward way expresses (again, via $En_S$) the consistency of S, is not S-provable if S is consistent (i.e., if $Con_S$ is true). This theorem is based on a formal deduction of G from $Con_S$ via a S-proof of the implication

$$Con_S \quad G$$

the S-provability of which entails the S-provability of G given a proof, in S, of $Con_S$. Since, by the first incompleteness theorem, G is not provable if $Con_S$ is true, neither is $Con_S$.

**B**. These and other theorems are closely connected with the use of the intuitive arithmetic from which the very concept of 'integer' is borrowed and also on the *way of presentation* of integers by 'numerals' or other formal objects used for representing the integers in LS. Roughly speaking, this way must correspond, in the metatheory, MS, of S, to the way used in the recursive function theory and also in arithmetical calculations (as these are dealt with in that theory).

This way is in no way unique. In many works the 'numerals' are introduced as the objects

$$0, \ 0', 0'', \ \ldots \ , \ 0''\overset{n}{\cdots}', \ \ldots \ ,$$

where **n** in the superscript of $0''\overset{n}{\cdots}'$ indicates the number, **n**, of occurrences of the symbol ''' (read: *stroke*). The use of the stroke-symbol is not always convenient and I have chosen to replace it by the 'successor function symbol', suc, so that the numerals become the '*traditional figures*':

$$0, \ suc(0), \ suc(suc(0)), \ \ldots \ , \ suc^n(0), \ \ldots \ ,$$

where the superscript **n** is the metatheoretical notation for the number of occurrences of the symbol suc. Thus, in terms of MS, the traditional figures become denoted by

$$suc^0(0), \ suc^1(0), \ suc^2(0), \ \ldots \ , \ suc^n(0), \ \ldots \ ,$$

and, to be considered as objects in LS, they must be fully 'unabbreviated', i.e., contain neither dots nor superscripts (and have the 'length' indicated, in MS, by the superscript; the words 'En$_S$-N of such figure' may refer, strictly speaking, only to these very objects in LS but not to their notations in MS).

However, still another way of presentations of integers dispences with the use of function symbols altogether and uses, instead, the '*numeroids*'

$$0, \ Su0, \ SuSu0, \ \ldots \ , \ Su^n0, \ \ldots \ ,$$

which are considered as *indecomposable* objects (none of which is a 'part' of another one). The superscript **n** is the metatheoretical notation for the 'length' of the numeroid which contains exactly **n** occurences of the '*particle*' Su (which is not considered as a function symbol). The superscript **n** — in $suc^n(0)$ as well as in $Su^n0$ — may be replaced by any meaningful notation of integer **n** in LS. Thus the notations $suc^r(0)$, $Su^r0$ arise for any term — or *termoid*, as I prefer to say — r in LS (which, depending on S, may also contain free variables). The numeroids, $Su^r0$, can also be denoted by $\nu[r]$ or, when the brackets are superflous, simply by $\nu r$.

The distinction between the traditional figures and numeroids is quite decisive because the calculating procedures and most arithmetical formulae deal with the former, not the latter. In order to overcome this deficiency of numeroids, I suggest to introduce in the system S the '*numeroidal axioms*' (for **n** = 1, 2, ...)

su-Su    $suc^n(0) = Su^n 0$

Su-su    $Su^n 0 = suc^n(0)$

(of course, only one of these two axiom schemata suffices).

Now, the use of numeroids can be accepted in constructions of arithmetic in S. In what follows, the Arabic figures

   0, 1, 2, ... ,

shall be considered as the (metatheoretical notations for) numeroids (the length of which they 'express').

Can the notation r which denotes the length of the numeroid $Su^r 0$, in LS, contain free variables? In particular, can this r be replaced by a variable, , for which a substitution is possible and which can be quantified as any other variable?

The answer depends on the specific choice of the system S or on its introduction. In what follows, LS and S shall denote, respectively, the language of the system S and this system itself, both extended by the allowance to admit the use of variables in notations of numeroids and traditional figures. The numeroid $\nu r$ will be called a '*variable numeroid*' when r contains a free variable.

Originally, a system denoted by S and containing numeroids will be supposed to contain no variable numeroid so that LS, S can be considered, in MS, only as extensions of LS, S. Then it will be possible to prove — in MS (to be appropiately extended) the consistency of (the required parts of) S relative to S. The requires S -axioms will be enumerable and the Gödel'-Rosser incompleteness theorems will be provable for S as they are provable for S.

**B 1**.    Yet a third way of representing integers in LS arises when this language contains 0 and 1 (or 0 and Su0) with their usual properties assumed in S. Let also + with its properties be contained in LS, S: then each positive integer, r, can be represented by the r-ary sum

   $1 + \overset{r}{...} + 1$

— or, if only the binary + is included in LS, by the presentation of such sum, when r>2, with the aid of parentheses which can also be omitted on the basis of 'association-to-the-left'. Even more simply, one can write a r-ary sum as

  1 ... 1

or, if need be, as

  1 .ͬ. 1 ,

where r indicates the number of occurreces of the symbol 1 in 1...1. In order to also include 0 in this presentation of integers, let them be modified as

  0 + 1 + .ͬ. + 1 ,

  0 1 .ͬ. 1 ,

respectively. (+ can be omitted only in informal parts of this text; the notation 1 .¹. 1 shall denote 1, etc.).

   Are these notations worse than $Su^r 0$? They can be arithmetized, with the aid of the unary constant function symbol, sg, introduced, in S, by the axioms

$def_1^{sg}$:  sg(0) = 0,

$def_2^{sg}$    (sg(suc( )) = 1)

and of the standard summation symbol,  , as the sums which can be denoted by

  $\sum_{=1}^{r}$ sg( ) .

   However, this notation is not always convenient because, in particular, its use presupposes that the variable   does not occur in r freely, a condition which must be systematically checked whenever this notation is used. Therefore, I suggest to replace the use of this notation by the use of the unary function symbol,   , introduced, in S, by the defining axioms

$def_1$ :    (0) = 0,

$def_2$ :    ( (suc( )) = suc( ( ))

which provide, in S, a formal proof of the formula

 =   :    ( =  ( ))

(see     ; this S-proof uses, as its non-logical axioms, only the given defining axioms and an instance of the induction schema

$$C(0) \ \& \ C(\xi) \supset C(suc(\xi)) \supset C(\eta)$$

with the formula $\eta = \varphi(\eta)$, or $\varphi(\eta) = \eta$, as its $C(\eta)$; the choice of the bound variable in the axiom schemata and defining axioms may now be considered immaterial.) Thus, in particular, for any LS-termoid r, in its occurrence(s) in $\varphi(r)$, $\xi$ may contain freely any variable including $\eta$, and the formula $\eta = \varphi(\eta)$, above, will have an easy S-proof. If **r** is a constant termoid, then its value, $\nu\mathbf{r}$, will be considered, in MS, as determined (though not necessarily as calculable). The formula $\eta = \varphi(\eta)$ shall entail in S (with the aid of the substitution of **r** for $\eta$), the equality

$$=_\mathbf{r}: \quad \mathbf{r} = \varphi(\mathbf{r})$$

whose right side equals the **r**-ary sum denoted, in MS, by the '**dot-numeroid**' $0+1+\overset{\mathbf{r}}{...}+1$ with **r** occurrences of $+1$'s in its right side. In MS, one can use the recurrent procedure which provides the S-proofs of the equalities

$$=^{dn}_{\varphi(\mathbf{r})}: \quad \varphi(\mathbf{r}) = 0+1+\overset{\mathbf{r}}{...}+1$$

and (using $=_\mathbf{r}$)

$$=^{dn}_\mathbf{r}: \quad \mathbf{r} = 0+1+\overset{\mathbf{r}}{...}+1.$$

The right side of $=^{dn}_\mathbf{r}$ must equal $\nu\mathbf{r}$, as well as $suc^\mathbf{r}(0)$; the S-proofs of the equalities referred to as

$$0+1+\overset{\mathbf{r}}{...}+1 = \nu\mathbf{r}\quad,\quad 0+1+\overset{\mathbf{r}}{...}+1 = suc^\mathbf{r}(0)$$

can also be obtained with the aid of a recursive procedure up to (the value of) of the constant termoid **r**. In this way, also each of the equalities can be S-proved equivalent to each of the equalities

$$=^\nu_\mathbf{r}: \quad \mathbf{r} = \nu\mathbf{r}$$

$$=^n_\mathbf{r}: \quad \mathbf{r} = suc^\mathbf{r}(0)$$

can be S-proved, for each constant termoid, **r**, in LS regardless to whether or not its value is calculable in S. While, fully unabbreviated, these formulae belong to LS; this conclusion itself, without being a formula in LS, will be considered as belonging to MS. Furtermore, the following statement holds in MS: when **r** is a quite definite (calculable or not) constant termoid in LS, then the En-N, $_\mathbf{r}$, of **r** is definite as well as the En-N, $_{\nu\mathbf{r}}$, of $\nu\mathbf{r}$, as well as the (least) En-N's of S-proofs of the equalities $=^{\nu}_{\mathbf{r}}$, $=^{n}_{\mathbf{r}}$ ; these En-N's all depend on both $_\mathbf{r}$, **r**, and as **r** is replaced by other constant termoids, $\mathbf{r}_1$, $\mathbf{r}_2$, ... $\mathbf{r}_i$, ..., $_{\mathbf{r}_i}$ will increase when (the value of) $\mathbf{r}_i$ increses as will the En-N's of the equalities $=^{\nu}_{\mathbf{r}_i}$, $=^{n}_{\mathbf{r}}$ and of their mentioned S-proofs.

The arithmetization of the 'grammary' of LS will actually require a consideration of termoids and their values in greater generality than so far considered. That is why **r** in $_{\nu\mathbf{r}}$ and $=^{\nu}_{\mathbf{r}_i}$ sometimes will be replaced by termoid r which may contain free variables so that the case of the variable numeroids arises. For the latter case the formulae

$=^{\nu}_{r}$:  . (r = $\nu$r), i.e.  . (r = Su$^r$0)

$=^{n}_{r}$:  . (r = suc$^r$(0))

— where . stands for any closure string which contains each variable free in r — must be considered; but how to do that when r in Su$^r$0, suc$^r$(0) does not occupy an argument place? The extension S  of S (mentioned on p.x, above) provides an answer to this problem.

The mentioned proof of the consistency of $S_\nu$ relative to S has to deal with an appropiatet revision of LS and its grammary; this revision will be described in Part II of this work, in a way which suffices for the substantiation of each particular $S_\nu$-proofs. For any use of S  in the current Part I , this relative consistency proof will consist 'essentially' in the systematic replacings, throughout the required S -proofs, of all occurrences of Su$^r$0, suc$^r$0 by the occurrences of (r),. for each termoid r, in LS, and of using, in case of need, the replacings based on the equalities $=^{dn}_{(t)}$ with constant termoid t. This use of these replacings transforms each S -proof into a S-proof, with the root modified only by the replacingings, so that if the root is 0 = 1 it remains unchanged.

In $S_\nu$ thus substantiated (w.r.t. S) the axioms dealing with substitution of termoids for variables can also be used when the latter have free occurrences via variables in a numeroid. In particular, 'the least number principle' ([KLEENE, 1952] §40, *149º.) can be applied to variable numeroids, i.e. the closure

:  μA(μ).   (A( ) &  ( <   ¬ A( ))

of the formula *149°. may contain occurrences of their variables μ, , via occurrences of numeroids via A(μ), A( ), A( ).

The equalities of the sort $=_{\mathbf{r}}^{\nu}$ will be used below only in cases where **r** is either 0 or a termoid which represents in LS a En-N, that is, a 'LS-image' of En-N, $_E$, of formal object E in LS. Each En-N, $_E$, of object E in LS must be a definite integer — and in order to consider the intuitive equalities

$$_E = \nu\, _E,$$

(considered in MS) their formalization require that $_E$ must be presented by a termoid, $\pi_E$, in LS, so that this equality can be expressed, in LS, as

$$=_{\pi_E}^{\nu} :\ \pi_E = \nu \pi_E,$$

which is an instance of $=_{\mathbf{r}}^{\nu}$, above. When E contains free, in E, occurrences of variables μ, , ..., $\pi_E$ shall contain, at the respective places, occurrences of LS-images $\pi_\mu$, $\pi$, ... of the En-N's $_\mu$, , ... , of these variables.

**Digression**. Subtleties connected with the use of variable numeroids will be considered closer in Part II of this work. These subtleties arise, in particular, when E is a 'hypothetical' object — in which case the very existence of E is only assumed, in MS, and such an assumption will be represented in LS in terms of the LS-image, $\pi_E$, of the En-N , $_E$, of E. The use of $\pi_E$ shall be connected with the use of $\nu\pi_E$ and the equalities $=_{\pi_E}^{\nu}$. Still, $\pi_E$ and $\nu\pi_E$ are constants — but they can be used, with the required generality, only when the language LS contains variable numeroids or equivalent means, such as the ' -numeroids' above which 'represent' the respective 'dot-numeroids', $0+1+\overset{\pi_E}{\ldots}+1$, where $\pi_E$ (which can be replaced by $(\pi_E)$) is allowed to contain free (in it) occurrences of variables.

Such subtleties, be they only implicit, are hidden behind Gödel's ( half-)proof of his first incompleteness theorem of 1931. At that time the non-S-provability of the 'self-referential' formula, G, (which (supposedly) 'expresses', in LS, its own non-S-provability) was established (in MS, using the current notations) by a *reductio ad absurdum* argument. By this argument it was assumed, again, in MS, that 'the' formula G is S-provable, i.e. that there exists a S-proof, $P_G$, of the formula G in S. That assumption was expresssed, in LS, by an equality

$$\mathscr{b}_S(\nu\mathbf{\mu}_G, \nu_G) = 0$$

in LS, where $\ell_S$ stands for a binary p.r. characteristic function symbol $En_S$-representing the proof relation in S and $\nu\mu_G$, $\nu_G$ are used as the LS-images of the En-N's $\mu_G$, $\nu_G$ of $P_G$ and G, respectively. The formula G was introduced as the negation

$$\neg\, \ell_S(\nu\mu_G, \nu_G) = 0$$

of the assumption formula $\ell_S(\nu\mu_G, \nu_G) = 0$, above — and with the aid of the (assumed in 1931) equalities

$$=\frac{\nu}{\mu_G}: \qquad \mu_G = \nu\mu_G\,; \qquad =\frac{\nu}{\pi_G}: \qquad \pi_G = \nu\pi_G \,—$$

the formula G was found to be, in S, incompatible with the assumption formula $\ell_S(\nu\mu_G, \nu_G) = 0$, above. In more detail; $\nu\pi_G$ must coincide with $\nu_G$ and, if G is S-provable, then $\ell_S(\nu\mu_G, \pi_G) = 0$ must hold — and due to the (metatheoretical) truth of $=\frac{\nu}{\mu_G}$, $=\frac{\nu}{\pi_G}$, also the equality $\ell_S(\nu\mu_G, \pi_G) = 0$ (or, in other notations, $\ell_S(\nu\mu_G, \nu_G) = 0$) must hold, and, the function symbol $\ell_S$ being p.r., it must be, in the case of S-provability of G, S-provable, all while the formula G is $\neg\, \ell_S(\nu\mu_G, \nu_G) = 0$ and is assumed to be S-provable. Thus, if G is S-provable then — using the equalities $=\frac{\nu}{\mu_G}$, $=\frac{\nu}{\pi_G}$ — a contradiction can be S-proved which is not the case if S is consistent. That is the use of the equalities $=\frac{\nu}{\mu_G}$, $=\frac{\nu}{\pi_G}$ in the proof, in MS, of the first incompleteness theorem for S.

Which were assumed in [GÖDEL, 1931] — the truth or the S-provability of $=\frac{\nu}{\mu_G}$, $=\frac{\nu}{\pi_G}$? Strictly speaking, Gödel used other notations, enumeration and terminology and this question will make more sense a soon as his ideas can be well-understood using the current enumeration, $En_S$, and the notations and terminology used in this work. That discusion will be saved for Part II of this work. However, without drawing any conclusions from the argument of this digression, it may be sumarized — as follows.

The formula G is considered, in MS, as 'explicitly displayed' — and, up to the use of practically unavoidable abbreviartions, G is even practically displayed, and the LS-image, $\pi_G$, of the En-N of G has the value $\nu\pi_G$ or $\nu_G$, considered in MS as calculable. $\pi_G$ is even a practically fulfillable object, and huge as the length of $\nu_G$ must be, this numeroid and a S-proof of $=\frac{\nu}{\pi_G}$ can both be considered, in MS, as achieved. However, this achievability is not quite obvious when $\nu\pi_G$ and (a S-proof of) the equality $=\genfrac{}{}{0pt}{}{\nu}{\mu_G}$ are considered — because the S-proof $P_G$ is only a hypothetical object.
(End of digression.)

**C**. In MS, it can be established that LS contains the binary p.r. function symbol $\mathscr{b}_S$ which En-*represents* the proof relation in S, i.e. for which the equality

$$\mathscr{b}_S(\nu s, \nu r) = 0$$

is S-provable iff s equals $En_S$-N of S-proof the end formula of which has the $En_S$-N r, and the equality

$$\mathscr{b}_S(\nu s, \nu r) = 1$$

is S-provable iff *not-* $\mathscr{b}_S(\nu s, \nu r) = 0$. (The subscript S in $\mathscr{b}_S$ may be dropped when only one system S is at issue).

With the aid of the symbol $\mathscr{b}_S$, the consistency of S can be expressed by the formula

$Con_S$ :    $\neg \quad \mathscr{b}_S(\ , \nu\pi_\ell) = 0$

where $\pi_\ell$ denotes the $En_S$-N of the '*anti-axiom*'

$\ell$:    $0 = 1$

and $\nu\pi_\ell$ — often abbriged as $\nu_\ell$ — the numeroid $Su^{\pi}\!/0$. Thus, $Con_S$ can also be displayed as

$\neg \quad \mathscr{b}_S(\ , \nu_\ell) = 0$

**C1**. Now, consider the integer    , which shall equal the least, if any, En-N of S-proof of the anti-axiom $\ell$, and shall equal 0 if no such En-N exists (i.e., if S is consistent). Provided that the disjunction

$TND_\ell$ :    $\mathscr{b}_S(\ , \nu_\ell) = 0 \ \lor \ \neg \quad \mathscr{b}_S(\ , \nu_\ell) = 0$

is recognized, in MS, as true, this definition of makes sense. Provided that TND_ℓ is provable in S and S contains definite descriptions with the axioms given for them in [ROSSER, 1953], can be introduced as the definite description

$$(\mathscr{b}_S(\ ,\nu_\ell) = 0 \ \&\ \neg\ (\ <\ \&\ \mathscr{b}_S(\ ,\nu_\ell) = 0).$$
$$\vee.\ = 0\ \&\ \neg\ \mathscr{b}_S(\ ,\nu_\ell) = 0\ )$$

and satisfies the disjunction

$$M(\ ): \quad \mathscr{b}_S(\ ,\nu_\ell) = 0\ \&\ \neg\ (\ <\ \&\ \mathscr{b}_S(\ ,\nu_\ell) = 0\ ).$$
$$\vee.\ = 0\ \&\ \neg\ \mathscr{b}_S(\ ,\nu_\ell) = 0.$$

which must be S-provable with the aid of the just mentioned axioms. Also, if the description theory is not included in S, the termoid can be included in LS as a 'new indecomposable constant termoid' for which the disjunction M( ) can be adjoined to S as a new axiom without loosing the consistency (which follows from Theorem 42, §74 in [KLEENE, 1952]). Here I shall consider M( ) as an axiom in S; the En-N, π , can be introduced, for definiteness, as if is considered as the description which is displayed above and the En-N, π , of the description symbol is introduced, in En, as equal to an odd prime, say as 113 (i.e., the 29$^{th}$ prime, $p_{29}$). More specifically, π being introduced, the En-N of M( ) becomes quite definite and it would be easy to redefine $\mathscr{b}_S$ as the proof-relation of the extension, S+M( ), of S which can be proved consistent relative to S and then to rewrite M( ) with $\mathscr{b}_S$ replaced by $\mathscr{b}_S'$ where $\mathscr{b}_S'$ denotes the redefined $\mathscr{b}_S$ of the extended system: also must be redefined as ´. This would not change much in the considerations, below, because (in accordance with mentioned Kleene's Theorem 42) the formulae

$$\neg\ \mathscr{b}_S(\ ,\nu_\ell) = 0 \quad \text{and} \quad \neg\ \mathscr{b}_S'(\ ,\nu_\ell) = 0$$

are equivalent and — whether S is contradictory or consistent — this equivalence is 'S-*provable*' (i.e. provable in S). So I'll continue with S and $\mathscr{b}_S$.

By virtue of the equality

$$= \nu\ : \qquad = \nu$$

which, eventually, will be found S-provable (whether its right side, $\nu$ , denotes 0 or a positive numeroid), the equality $\ell_S(\ ,\nu_\ell) = 0$ can be replaced in M( ) by $\ell_S(\nu\ ,\nu_\ell) = 0$ and, thus, the disjunction

M($\nu$ ):     $\ell_S(\nu\ ,\nu_\ell) = 0$ & ¬ ( < & $\ell_S(\ ,\nu_\ell) = 0$ ).
             ∨.     $= 0$ & ¬   $\ell_S(\ ,\nu_\ell) = 0$

can be deduced from M( ) and, hence, found S-provable. The members of this disjunction entail, in the propositional calculus, the members of the simpler disjunction

m($\nu$ ):     $\ell_S(\nu\ ,\nu_\ell) = 0 \vee \neg\ \ell_S(\ ,\nu_\ell) = 0$

which also can be displayed as

d($\nu$ ):     $\ell_S(\nu\ ,\nu_\ell) = 0 \vee \mathrm{Con}_S$ .

The value of the termoid $\ell_S(\nu\ ,\nu_\ell)$ can be calculated, in S, because the function symbol $\ell_S$ is p.r. and both $\nu$ and $\nu_\ell$ denote numeroids. If S is consistent, then this value is 1 and the equality

$$\ell_S(\nu\ ,\nu_\ell) = 1$$

is S-provable — and it remains S-provable, as any other closed formula in LS, also in the case when S is contradictory.

Thus, regardless to whether S is consistent or contradictory, this equality is S-provable. It entails, in S,

$$\neg\ \ell_S(\nu\ ,\nu_\ell) = 0,$$

i.e. the negation of the first member of the S-provable disjunction d($\nu$ ); therefore, the second member

$$\mathrm{Con}_S$$

of this disjunction is S-provable.

**Remark**: A similar statement conflicting with the Gödel's second incompleteness theorem was claimed already in 1966 (any many times since) by E.

Wette — though I could not obtain, so far, any clear expounding of his proofs. Nevertheless, it seems to me appropiate to call the ensuing contradiction the *Gödel-Wette Paradox*.

(End of remark).

**D**. The S-proof described in **C1**., above, for the formula $Con_S$, presupposes the possibility to S-prove, for any constant termoid, r, such as     (in LS), the equality r = $\nu$ r (or its equivalent with a different sort of 'numeroid' used in the role of $\nu$ r). Such S-proofs are available in $S_\nu$ but not necessarely in S; they use a sort of induction on these 'numeroids' which is not necessarily available in S because the language LS is not supposed to provide for the 'variable numeroids': The induction is carried out in MS (the language of which does contain such variables) and it consists in performing a recursion the $s^{th}$ step of which delivers a S-proof of the equality s = $\nu$ s where s may, in particular, be taken as    . This approach presupposes that s — in particular,    — occurs in the sequence of integers 0, 1, 2, ... .

If S is consistent, them    must equal 0 and thus occur at the very start of this sequence, fuzzy as it can be considered in connection with various metatheoretical problems. Also this idea belongs to MS. However, the task now consists in formalizing a proof of    = $\nu$    in S without adressing any philosophical problems.

The recursion on s proceeds by the steps from integers to their (immediate) successors in 'the' series of 'all' integers; unless    = 0,    must be En-N of S-proof. *A priori*, it is thinkable that such En-N's exceed all 'usual' integers. In a case of    , $\nu_\ell$) = 0,    must be a En-N of contradiction proof in S — and such an object is not expected to be 'usual' (unless S is contradictory or otherwise a 'pathological' system of arithmetic). Perhaps, $\ell_S($   , $\nu_\ell$) = 0, and, still, $\ell$, i.e., 0 = 1, is not entailed by that equality because the integer    is not an 'Archemedian' integer — though, it can be admitted that    , as any integer, must be a finite ordinal (in the sense that the set of all integers which are less than.    is finite and well-ordered).

A proponent of Gödel's second incompleteness theorem can insist on this possibility. If    = $\nu$    = suc (0), then, if these equalities are S-provable, then so is $Con_S$ and the contradiction, due to K. Gödel, can be proved in S, and, of course, also the equalities    = $\nu$    and the equalities $\ell_S(\nu$   , $\nu_\ell) = \ell_S($   , $\nu_\ell) = 0$ can then be S-proved (since $\nu$    is choosen as the least such numeroid — cf. §41: \*166, \*172 in [KLEENE, 1952] — and that is normal since ℓ is found S-provable. But this argument presupposes that the equalities = ν = suc (0) are, in fact, S-provable. And what if the recursion cannot reach the th step?

Somewhat must happen with MS as the result of the current considerations of enriched by the recursive proofs of S-provability of the equalities s = νs for all integers s. A volume can be written in order to substantiate the soundness of this recursion. Then it will be clear that the use of identities presupposed by the introduction of En$_S$ are not always reliable.

But is not a revision of the (traditional) idea of identity a too strong means? Does it not simply suffice to assume that that the recursion and induction can be used, in conection with    , only in a less general way than the standard n-to-n+1 approach which presupposes, for each non-zero integer, the 'existence' of its immediate predecessor?

However, it has been found that also a demand of restricting the generality of the use of induction would not help to avoid the confrontation with Gödel's theorem with En-N's — standard or not— like    .

**D 1**.   I recall that Gödel's theorems have been extended, starting with [Rosser, 1936], to the case when S contains infinitary rules of 'generalization' and MS contains arithmetization involving transfinite ordinals. Generalizations of Gödel's theorems in such directions have been important (however vague). Here I don't consider them closer. I admitt the possibility that when En$_S$ involves the use of some ordinals the theory of the latter does not always allow for a proof of the S-provability of the implications of the sort

$$\ell_S(\quad, \nu_\ell) = 0 \quad \ell$$

where    stands for an ordinal.

Here I only wish to claim that the '*proof-truth implications*'

$$\ell_S(r, \nu_E) = 0 \quad E$$

are provable for the systems S with the 'standard proof theory' using En$_S$ which admits the induction on the complexity of S-proofs. The logic of such systems S will be closed, with *modus ponens* (MP) as the only rule of inference. In a standard way, the simplest S-proofs shall be '*trivial*',

{E},

where E is an axiom in S, and they shall be assigned the En-N's, $\pi^{\{E\}}$, introduced as $2^{\pi_E}$ where the exponent $\pi_E$ equals En-N of the formula E.

The non-trivial proofs in S shall be obtained, each, from two S-proofs, $P_C$, $P_{C \; D}$, the end-formula of which are C and C D, respectively, with coinciding C's, by the application of the step

MP: $$\frac{C \quad C \; D}{D}$$

the '*conclusion*' of which shall be the end-formula, D, of the S-proof, , thus obtained; $P_C$, $P_{C \; D}$ shall be considered as the '*branches*' growing in $P_D$ from the '*premises*' $P_C$, $P_{C \; D}$ .of this MP. When the En-N's $\pi^{P_C}$ and $\pi^{P_{C \; D}}$ of the branches $P_C$ and $P_{C \; D}$ are given and $\pi_D$ denotes the En-N of the formula D, then $P_D$ shall be given the En-N, $\pi^{P_D}$, equal to the product

$$Mp(P_C, P_{C \; D}): 2^{\pi_D} \cdot 3^{\pi^{P_C}} \cdot 5^{\pi^{P_{C \; D}}}$$

Thus each S-proof shall be introduced, in S, together with the En-N's attached (or 'assigned') to it in MS. MS shall recognize only such S-proofs which are either 'trivial', i.e., {E} where E is an axiom in S, or obtainable, by a step MP, from two previously obtained S-proofs, $P_C$ and $P_{C \; D}$. In MS, the following 'MP-induction' will be recognized as a general principle for accepting statements, **D**(**P**), about S-proofs, P:

If **D**(**P**) holds for each trivial S-proof, {E}, considered as P, and if, for each S-proof, $P_D$ — of which the end-formula, D, is the conclusion of the last step, MP, as described, above — **D**($P_D$) holds as soon as **D**($P_C$) and **D**($P_{C \; D}$) hold, then **D**(**P**) holds for any S-proof P.

In particular, **D**(**P**) may consist in that P has a certain property formulated in MS. This property may consist, in particular, in the finiteness, and also in the possibility to introduce the En-N, $\pi^{P_D}$, for $P_D$; in the uniqueness of this $\pi^{P_D}$, etc. Other such properties may consist in the 'truth' of the end-formula of D, and then the MP-induction leads to a proof, in MS, of the consistency of S. (The possibility of such a consistency proof presupposes the possibility to formulate, by means of MS, the 'truth' of the end-formula).

The arithmetization of MS, achievable with the aid of $En_S$, entitles one to reduce the MP-induction to a specific case of arithmetical induction, to be designated as 'mp-induction', which, in turn, can be reduced, in the arithmetic, to a course-of-values induction on the En-N's of S-proofs. This course-of-values induction does not nessearily presuppose that for each $En_S$-N, $\pi^{P_D}$, the S-proof of the equality $\pi^{P_D} = \nu\pi^{P_D}$ is S-provable.

The use of mp-induction shall correspond to the use of MP-induction in MS.

**D 2**. The mp-induction will be presented in an extension, $S^\nu$, of the system S the language, $LS^\nu$, of which shall contain, among its atomic formulae, the objects

$F_r$

where r may stand for any termoid in $LS^\nu$ (which may, in particular, contain free variables). When r denotes $En_S$-N of quite definite closed formula in LS, then $F_r$ shall be interpreted, in MS, as denoting this formula; otherwise $F_r$ shall be interpreted as $F_{\nu_\ell}$ i.e. as $\ell$. (Cf. the use of the metatheoretical notations $A_r$ in [KLEENE, 1952], §42, p. 206.) Derivations in $S^\nu$ can be presented with the use of such atomic formulae in which case the variables free in r shall be considered as parameters of these derivational presentations. The postulates in $S^\nu$ used in these presentations, below, will be found provable in $S_\nu$ and thus found justified for use in this work. Any closed formula in $LS^\nu$ used as $S^\nu$-provable in this work will be found $S_\nu$-provable.

The formulae to be found $S^\nu$-provable with the aid of mp-induction shall be (as considered in the open logic) the 'proof-truth-implications'

$pti_{(\mu)_0}$:    $\ell_S(\mu,(\mu)_0) = 0$      $F_{(\mu)_0}$

where μ stands for the variable of this induction; this variable will be held constant till the last formula of this induction (although in the conclusion of this induction μ must be bound). Remark that *only* values ≤ ωβ of μ will be of relevance to the use of this induction in this work (which thereby retains a quite finitistic approach). In MS, the equality ωβ = $\nu$ωβ will be recognized, and the instances, with $\nu$r ≤ ωβ,

$\ell_S(r,(\nu r)_0) = 0$      $F_{(\nu r)_0}$

of this induction shall make a finite set of these implications.
The ωβ-instance of the implications $pti_{(\mu)_0}$,

$$\text{pti}_{(\omega\beta)_0}: \quad \mathscr{b}_S(\omega\beta,(\omega\beta)_0) = 0 \quad\quad F_{(\omega\beta)_0},$$

shall be found $S^\nu$-provable. If $\omega\beta$ is positive so shall be $(\omega\beta)_0$ which, in this case, must equal $\nu_\ell$, so that $F_{(\omega\beta)_0}$ will coincide with $F_{\nu_\ell}$, i.e. with $\ell$, and this implication will be, in this case, the required

$$\mathscr{b}_S(\omega\beta,\nu_\ell) = 0 \quad\quad \ell.$$

The antecedent of $\text{pti}_{(\mu)_0}$ is displayed as

$$\mathscr{b}_S(\mu,(\mu)_0) = 0$$

because, regardless to whether $\mu$ equals En-N of trivial or non-trivial S-proof, the root of the latter must have the En-N which equals $(\mu)_0$. In general, $\mathscr{b}_S(r,s) = 0$ shall entail

$$\mathscr{b}_S^1(r) = 0 \text{ and } (r)_0 = s$$

where $\mathscr{b}_S^1$ is the symbol of the characteristic function which En-represents the proof property in S. $\mathscr{b}_S^1(r)$ shall equal 0 only when either $r$ equals En-N of trivial S-proof, in which case $r = 2^s$ and $(r)_0 = s$ entails

$$r = 2^{(r)_0};$$

or, $r$ equals En-N of S-proof which terminates in MP and has the En-N equal to $Mp((r)_1,(r)_2)$ where $(r)_1, (r)_2$ equal the respective En-N's of the branches growing from this MP's premises; these branches must be S-proofs and these En-N's must therefore satisfy the equalities

$$\mathscr{b}_S^1((r)_1) = 0, \; \mathscr{b}_S^1((r)_2) = 0;$$

which will be important for the use as the induction hypotheses.

The basis of the mp-induction shall deal with the case when $\mu$ equals En-N of trivial S-proof; in this case the antecedent, $\mathscr{b}_S(\mu,(\mu)_0) = 0$ of $\text{pti}_{(\mu)_0}$ shall be equivalent, in S, to the conjunction

$$\mu = 2^{(\mu)_0} \; \& \; \text{ax}_S((\mu)_0) = 0$$

— where $ax_S$ is the p.r. characteristic function symbol for the property of being En-N of axiom in S; and the implication $pti_{(\mu)_0}$ shall be, for this case, entailed by the condition

$$\mu = 2^{(\mu)_0},$$

its incompatibility with $\mu = Mp((\mu)_0, (\mu)_2)$ and the implication

T-Ax$_{(\mu)_0}$:      $ax_S((\mu)_0)$      $F_{(\mu)_0}$.

In this case, $F_{(\mu)_0}$ denotes an axiom in S and the S-provability of T-Ax$_{(\mu)_0}$ shall be entailed by the S-provability of the propositional axiom

$j_{F_{(\mu)_0}}$:      $F_{(\mu)_0}$      $(ax_S((\mu)_0) = 0$      $F_{(\mu)_0})$

In $S^\nu$, for each axiom E of S and each closed formula D in LS, the implication

    D   E

shall be accepted as an axiom on the basis of the general applicability of the argument of the last paragraph. The presence of these axioms in $S^\nu$ will trivialize the consideration of the basis of the required mp-induction.

In the case when $\ell_S(\mu, (\mu)_0) = 0$ and $\mu = Mp((\mu)_1, (\mu)_2)$ are given, the equalities $\ell_S^1(Mp((\mu)_1, (\mu)_2)) = 0$, $\ell_S^1((\mu)_1) = 0$ and $\ell_S^1((\mu)_2) = 0$ can be deduced, as well as the inequalities

$$(\mu)_1 < \mu, \; (\mu)_2 < \mu;$$

thus the induction hypotheses are available as well as the antecedents $\ell_S((\mu)_1, ((\mu)_1)_0) = 0$ and $\ell_S((\mu)_2, ((\mu)_2)_0) = 0$ of the truth-proof implications

$pti_{(\mu)_1, 0}$ : $\ell_S((\mu)_1, ((\mu)_1)_0) = 0$      $F_{((\mu)_1)_0}$

$pti_{(\mu)_2, 0}$ : $\ell_S((\mu)_2, ((\mu)_2)_0) = 0$      $F_{((\mu)_2)_0}$

which are provided by these induction hypotheses; the conditions $\ell_S(\mu, (\mu)_0) = 0$ and $\mu = Mp((\mu)_1, (\mu)_2)$ entail, in accordance with the definition of $\ell_S^1$ for this case, that $F_{((\mu)_1)_0}$ and $F_{((\mu)_2)_0}$ are premises, C and C   D, of a MP the conclusion, D, of which has En-N equal to $(\mu)_0$, so that

D coincides with $F_{(\mu)_0}$ and the required implication

$$pti_{(\mu)_0} : \ell_S(\mu, (\mu)_0) = 0 \quad F_{(\mu)_0}$$

is entailed, in $S^\nu$, by the conditions of this case.

Technically, this argument requires a rather long formalization though not nessearily as long as the one which I have found in 1996 and communicated to my friend and long-time collaborator Catherine Christer Hennix who has abbriged that work by more than half. Below, her version will be exposed; I appreciate her participation as my co-author of the current work.

Besides that, the substantiation of the use of $S_\nu$, $S^\nu$ for establishing the S-provability of $Con_S$ requires a use of an old theorem of mine to the effect that if the consistency of a formal system, W, is proved in a consistent system V, then W is consistent relative to V and to any system, U, in which V is proved consistent. (See my Supplement I to my translation [KLEENE. 1957] of [KLEENE, 1952]. Here W, V, U are formal systems of the same nature as the system presently denoted by S, above. In 1957-61 S. Feferman has considered many such systems and sorts of them; these researches have been exposited in his paper [FEFERMAN, 1960], and, probably, that theorem of mine was covered also by results of that paper. Also in that paper he mentions more than once ambiguities and discrepancies in the litterature of that time (see his remarks on pp. 37 and 78). The clarity in this field did not increse since, and now, in order to achieve it, the study of presentations of the axioms of 'S-like' formal systems must be resumed. All S-proofs leading to the S-provability of $Con_S$ must be exposited explicitly in order to accomplish the introduction of the unary function symbol $ax_S$ on which the axioms defining $\ell_S$ and $\ell_S^1$ depend. That is possible, in particular, because the En-N's can be introduced for these symbols regardless to their interpretations or the choice of of the axioms which define them. But that can be done exhaustively only after the S-proofs of the S-proof of the formula $Con_S$ are completely displayed and then scrutinized. In the present Part I of the current research the deduction of $Con_{Ari}$ (where Ari is sufficiently rich and smooth formal system for arithmetic) from the formula T-Ax (for that system) is exposited in a way which does not require the use of the equalities $=_r^\nu$ with constant termoids **r** in LAri. In Part II also the S-proofs of these equalities will be described and T-Ax (for Ari) will be found to have the Ari-provable closures.

Revere, June 9, 2001, Alexander S. Yessenin-Volpin

# INTRODUCTION.

0. The task of the present paper consists in studying a paradoxical derivation of a (standard) version of Gödel's consistency formula for arithmetic. This derivation is relatively long on account of its finitistic character. However, since a paradoxical result ensues it becomes essential that every step of the derivation is exposed in order to spot any possible logical mistake or, else, identify the cause of the ensuing paradox. The length of this derivation has made it mandatory to carry its termination into a Part II of the present paper which is designated as Part I. Before commencing with the paradoxical derivation itself, we outline the approach which has been taken for the presentation of this derivation in the present work, immediatelly below.

0.1. Let Ari denote a fragment of primitive recursive (p.r.) arithmetic — assumed formally to be consistent — and for which Gödel's incompleteness theorems are provable. Let MAri denote the metatheory of Ari. Further, let **En** denote a 1-1 enumeration in terms of which $\mathscr{b}$ 'En-represents' the p.r. characteristic function of Ari's *proof relation*; let     be an unrestricted variable in the language, LAri, of Ari and $\nu_f$ the Gödel number (in terms of **En**) of the falsum symbol $f$ (which denotes the 'anti-axiom' 0 = 1); then, the following  -term, denoted by    , shall belong to the language LAri:

$$(\mathscr{b}(\ ,\nu_f) = 0\ \&\ \neg\quad (\ <\quad \&\ \mathscr{b}(\ ,\nu_f) = 0).$$
$$\vee.\quad = 0\ \&\ \neg\quad \mathscr{b}(\ ,\nu_f) = 0).$$

(In (approximate) words: given **En**,     shall equal the smallest gödel number of proof, in Ari, of $f$ if there exists such a proof or, else,     shall equal zero and there exists no (Gödel number of) proof, in Ari, of $f$.)

Then, assuming the acceptance, for Ari, of the following instance of *tertium non datur*:

TND$_f$ :     $\mathscr{b}(\ ,\nu_f) = 0\ \vee\ \neg\quad \mathscr{b}(\ ,\nu_f) = 0;$

the following formula , labeled (M-   ), is Ari-provable:

(M-   ):    $\mathscr{b}(\ ,\nu_f) = 0\ \&\ \neg\quad (\ <\quad \&\ \mathscr{b}(\ ,\nu_f) = 0).$

$$\vee. \quad = 0 \;\&\; \neg\; \mathscr{b}(\;,v_{\ell}) = 0.$$

(Intuitionistically, only the double negation of this formula is Ari-provable.).

Since the function symbol $\mathscr{b}$ is p.r. it is decidable if $\mathscr{b}(\;,v_{\ell}) = 0$ or $1$, *provided* that is p.r. computable. By the interpretation given to the symbols in M-, it shall follow that, if $= 0$, then $\mathscr{b}(\;,v_{\ell}) = 1$ because 0 shall, by definition, not be a Gödel number of any proof (or any object in LAri); and if $0 < \; < $, then $\mathscr{b}(\;,v_{\ell}) = 1$ since equals the *smallest* Gödel number of Ari-proof (if any) of $\ell$. Thus, $\mathscr{b}(\;,v_{\ell})$ equals 1 unless Ari is inconsistent in which case $\mathscr{b}(\;,v_{\ell}) = 0$ holds and denotes the smallest Gödel number of Ari-proof of $\ell$. That is, unless Ari is inconsistent, $\mathscr{b}(\;,v_{\ell})$ will always equal 1 and will equal 0. So, in MAri, under the formal assumption that Ari is consistent, the first member of the disjunction M- is disprovable and, therefore, the second member of this disjunction,

$$= 0 \;\&\; \neg\; \mathscr{b}(\;,v_{\ell}) = 0,$$

follows. The second member of this conjunction may also be written as

(L- ) $\qquad \mathscr{b}(\;,v_{\ell}) = 0 \quad \ell$

the Ari-proof of which contradicts Gödel's second incompleteness theorem which assert that the formula L- is provable only if Ari is inconsistent.

**Remark**. Since is unrestricted, the formula $\neg\; \mathscr{b}(\;,v_{\ell}) = 0$ is a standard formula expressing, in LAri, the (classical) consistency of Ari. The formula L- shall be abbreviated by Con $_{Ari}$ or simply by Con$_{Ari}$. This formula is provably equivalent, in Ari, to Gödel's original formula as given in [Hilbert-Bernays, 1939]. Clearly, also the equality $= 0$ is provably equivalent, in Ari, to the latter.

0.2. The proof of the Gödel's consistency formula in a consistent formal system for which also Gödel's second incompleteness theorem can be pro-ved is an incident which we shall refer to as the *Gödel-Wette paradox*. Remark that this is a paradox allowed by the *metatheory*, MAri, *rather* than a direct inconsistency proof of the formal system Ari itself — *provided* that Gödel's second incompleteness theorem is accepted for Ari. In this respect, our claims are not as strong as Wette's — cf.[WETTE, 1971],[WETTE, 1974].

0.2.1. The task of obtaining an Ari-proof of L- will be reduced to the task of obtaining an Ari-proof of the ' -instance',

(L-ωβ): $\mathscr{b}(\ ,\nu_\ell) = 0 \quad \ell$ ,

of L- , which is Ari-provable as formula 667. of the proof 1.-686a.of L- , below. The Ari-provability of the formula L-ωβ suffices for the Ari-provability of $\mathscr{b}(\ ,\nu_\ell) = 0 \quad \ell$ , formula 686., below, since the former formula disproves the first member of the disjunction M- . Inspection of the the formulae 633 - 686a. of the proof 1.-686a., below, shows that only *propositional* steps are involved in the obtaining of L- from L-ωβ.

In order to prepare this last step, the task of proving, in Ari, the formula L-ωβ shall be reduced to a proof by a special form of induction, called 'mp-induction', of the 'open reflection formula' (also refered to a the *proof-truth*-implication (p.t.i.) in the Preface) $\mathscr{b}_1(\mu) = 0 \quad F_{(\mu)_0}$ in which $\mathscr{b}_1$ denotes the characteristic function symbol of Ari's *proof property* and $F_{(\mu)_0}$ denotes a closed formula in LAri per an extended use of Kleene's metamathematical notations for closed arithmetical formulae in [KLEENE, 1952] — and the current use of these notations is captured by an open (conservative) extension, $\text{Ari}^\nu$ , of Ari. The details are given in Sections 5, 6 and Appendices C, D1 and D2, below.

Essential for the convincingness of the proof 1.-686a., below, is the possibility to restrict the free variable denoted by μ to values in the proof of the mp-induction formula and, hence, throughout the proof 1.-686a.. This possibility will be confirmed and so will the Gödel-Wette paradox. However, these confirmations will be completely substantiated only with the completion of Part II of this work (in preparation).

0.3. In the present Part I we shall proceed as follows; (i) first we set up a *deduction* 1.-686a. of Con$_{\text{Ari}}$ from a *hypothetical* formula T-Ax: $\text{ax}_{\text{Ari}}((\mu)_0) = 0 \quad F_{(\mu)_0}$, , where $\text{ax}_{\text{Ari}}$ denotes the characteristic function symbol for Ari's set of (En-N's of) axioms. (ii) Since the shortest proof of T-Ax known to us exceeds the length of the present paper we shall, secondly, rather prove one of its consequences, *viz.*, formula 7. (which serves as the 'basis formula' for the current application of mp-induction); $\text{Ant}_1^{\mathscr{b}_1}(\mu) \quad F_{(\mu)_0}$, where $\text{Ant}_1^{\mathscr{b}_1}(\mu)$ defines the condition for the formula $\mathscr{b}_1(\mu) = 0$ to be provable in case $(\mu)_0$ denotes En-N of axiom in Ari. The proof of formula 7. is sketched in Appendix D2. Finally, (iii) formula 630., an instance of the mp-induction schema, is proved in Appendix C. In this way, all hypotheses of the deduction 1.-686a. are eliminated and this deduction is now also a proof in $\text{Ari}^\nu$. By systematically converting the Kleene notations to the corresponding formulae in LAri, the latter proof, *eventually*, becomes a proof in Ari. By this last step, the Gödel-Wette paradox can now be verified in MAri. However, the details of this verification as well as the resolution of the ensuing paradox shall be given in Part II where also other metamathematical paradoxes closely connected with the present work will be exposed related to completeness, decidability, truth-definitions and a resurgence of Richard's paradox with regard to Ari or its extensions (as here considered).

1. This paper is organized as follows: Section 2 introduces the basic notions of the language, LAri, of arithmetic and fixes a 1-1 enumeration **En** — but, to remain concise, skips over the usual formation rules for terms and formulae; however, the notion of deduction is detailed together with the assignments of the corresponding integers to formal objects in accordance with the fixed enumeration **En**. Section 3 presents the logical axiom schemata, mentiones *some* of the non-logical axioms and states the rules of inference. Section 4 defines the proof property and relation for the formal system Ari. In Section 5, which is central to this work, two (basically conservative) extensions, $\text{Ari}^+$ and $\text{Ari}^\nu$, of Ari, are introduced together with a substantiation of the additional axioms. Finally, Section 6 is a presentation of the proof, 1.-686a., by means of which the Gödel-Wette paradox can be claimed to exist. In addition, there are four appendices the last two of which concern the provability of two (well-known) formulae which are cri-tical for this work. However, the complete proofs of these two formulae will exceed the length of the present paper and therefore they will appear, in full, as Part II of this paper. In addition, Part II will exhaustively list all non-logical axioms used in this work.

Notes are collected as endnotes and are found directly after the main derivation. Referencs are collected at the end of the present paper.

2. Most of the metamathematical terminology will be standard with the only exception being the use of the notions of *numeroid* and *termoid* and the use of En-N's in place of Gödel numbers. These notions will be sufficiently explained, below. In general, the En-N of any object, E, in LAri shall be denoted by ⌜E⌝, and this notation wil be specified below for different cases of E. The object E shall always be determined by ⌜E⌝ so that the equality ⌜C⌝ = ⌜D⌝ shall entail the coincidence of C, D. The uniqueness of ⌜E⌝, for each object E in LAri, shall be entailed by by the introduction of each ⌜E⌝ as a positive integer and the uniqueness of the Gaussian decomposition for positive integers:

$$\mathbf{m} = 2^{(\mathbf{m})_0} \cdot 3^{(\mathbf{m})_1} \cdot \ldots \cdot p_{\ell(\mathbf{m})}^{\ell(\mathbf{m})},$$

for each positive **m** (see pp.)

LAri shall denote a certain fragment of a language of p.r. arithmetic and **En** shall denote a 1-1 *enumeration* of LAri which uniquely assignes an *enumeration number* (En-N) to each object in LAri. Each such number shall be intended as a *Ziffer* (in the sense of [Hilbert-Bernays, 1939]) which in this work shall be represented by an entity called 'numeroid', i.e. the equality between such an integer and the corresponding numeroid belongs to MAri (see below). LAri's expressions partition into two classes, *viz.* (i) the *indecomposable symbols* and (ii) all other expressions.

2.1. Excepting function symbols, LAri shall contain the following four *sorts* of *indecomposable symbols, only*:

(**1**). *variables*: $\mu^1, \mu^2, \mu^3, \ldots$; each variable, $\mu^i$, shall be assigned the En-N $29^i$, for i denoting a *positive* integer 1, 2, 3, ...;

(**2**). *numeroids*: 0, Su0, SuSu0, ..., $Su^n 0$, ...; each numeroid, $Su^n 0$, shall be assigned the En-N $23^{suc(n)}$ (= $23^{n+1}$), where suc denotes the function symbol for succesor and **n** takes all integer values;

(**3**). *binary predicate symbols*: = , < ,   ; assigned, respectively, the En-N's 15, 25 and 35;

(**4**). *logical symbols*:   , & , ∨ ,    ,    ,    ; assigned, respectively, the En-N's 3, 5, 7, 11, 13 and 113.

Terms (here usually refered to as 'termoids', see below) and formulae are inductively defined as certain trees in the usual way so that given the above assignment of En-N's to LAri's indecomposable symbols, any termoid or formula may, uniquely, be assigned an En-N accordingly.

The following *special* notations shall be used: if t is an any term of LAri's, then $_t$ shall denote its En-N while if t is a closed term of LAri's its En-N may also be denoted by $\tau_t$; and if E is an any formula in LAri, then $_E$ shall denote its En-N while if E is a closed formula in LAri, $\nu_E$ may also denote its En-N  Unless otherwise specified, these special notations stand for *numeroids*, i.e. formal objects of LAri.

The function symbols in LAri shall be, only, unary, binary and ternary. Their En-N's will be introduced, for the unary symbols, as $71^1$, $71^2$, ..., $71^{ta_1}$, for the binary symbols, $73^1, 73^2$, ..., $71^{ta_2}$, for the ternary symbols — as $79^1, 79^2$, ..., $79^{ta_3}$, where $ta_1$, $ta_2$, $ta_3$ stand for three integers which always will be supposed given in a Table, Ta (=  $Ta_i$, i = 1,2,3), in which these En-N's will be specified for each function symbol in LAri. In particular, the function symbols

   suc, + , · ,

will be assigned the En-N's $71^1 (= 71)$, $73^1 (= 73)$, $73^2$, respectively. Thus, the symbols for the functions used in the system of [KLEENE, 1952] chapter IV, §19, are introduced, and for each other function symbol they will apper in $Ta_i$, i = 1,2,3, in an order which will be, in case of need, specified. In general, the En-N of function symbol,  , in LAri, shall be denoted by    for its special notation.

**Add**(4).    Since Ari will, eventually, be subjected to a *finitistic* interpretation we must insist that negation (¬) is a defined symbol in LAri (via   and the absurd equality   0 = Su 0  ). So, for any formula E in LAri, the implication E    0 = 1 — or, shorter, E    ∫ , will often be abriged as ¬ (E) — or as ¬ E when the parenthese can easily be restored.

**Add**(2).    Also remark that each numeroid, $Su^n 0$ , with exactly n occurrences of the symbol Su , is not only an indecomposable object but, in addi-tion, is *uniquely* obtainable, except for 0, by concatenating the symbol Su to its predecessor $Su^{n-1}0$. Thus, if **m** < **n**, then $Su^{m}0$ is *not* identified with a 'part' or 'initial segment' of $Su^{n}0$. The bold letters **m**, **n**, above, denote *parameters* taking all non-negative integer values. Remark further that, *in* the *metatheory*, MAri, *each En-N is considered as equal to its numeroid*. For further details, refer to ENDNOTE A.

2.1.1. All objects, i.e. termoids (including, in particular, the indecomposable termoids (**1**), (**2**), above), and formulae, as well as deductions and proofs, are thus to be assigned En-N's in the usual way in terms of their tree constructions (say, as in [KLEENE, 1952]). It suffices here to limit the description of this assignment to those cases which explicitly will be at issue and, thus, we forego all irrelevant details. Thus, only the assignment of En-N's to deductions and proofs will be shown in some detail sufficient for the purpose of the present work.

(The notions of 'termoid', 'formuloid', 'deductoid' and so on (see [YESSENIN-VOLPIN, 1970]) are used by us in contexts for which *no* (genuine) *consistency proof is available in the litterature* and, thus, this (genuine) lack makes the use of the traditional terminology ('terms', 'formulae', 'deductions') (genuinly) *potentially risky* in such contexts — in particular in the vicinity of (potentially genuine) paradoxes. However, in this paper we restrict our use of the anti-traditional or ultra-intuitionistic vocabulary to the notions of *termoid* and *numeroid*, only — but we will reintroduce, when need be, the full anti-traditional terminology as a necessary means for resolving the Gödel-Wette paradox and other paradoxes connected with the ramifications of the approach explored in the present work).

A deduction may be carried out either in the *closed* logic, **Lo**, or in the *open* logic, **lo**. The closed logic requires that all formula in a deduction **E** be closed, (the formulae 'of **E**', being those which are members of the tree of the deduction); in the open logic this restriction is lifted. This concerns, in particular, the axioms of the system; they are introduced differently for the versions of Ari with the logics **Lo** and **lo**. These two versions shall be denoted by Ari and Ari°, respectively. Derivations (i.e. deductions and proofs) in the open logic are useful for the purpose of shortening their presentations.

2.1.2. Let **E** be a deduction in Ari or its extension(s).

Then **E** and its deduction tree, $\mathscr{T}_{\mathbf{E}}$, shall both be called *trivial* if **E** consists of a single formula, E, in which case E must be either (i) a *closed* formula, or, if the deduction is given in the open logic, (ii) E may be any formula; and if **E** is also a proof, E must also be an *axiom* in Ari, Ari° or their respective extension(s) (as defined in this work).

On the assumption that $\mu_E$, the En-N assigned to E, is given, the trivial deduction **E** as well as its deduction tree $\mathscr{T}_E$ shall both be assigned En-N

$$2^{\mu_E}$$

whether the logic is closed or open. If E is the root formula (here also denoted by $\nu E$) of a trivial deduction **E**, then $\mu_E$ shall denote En-N of this deduction. It is easy to see that the following equality, also called the *deduction (proof)-root-correlation*, holds;

($\nu 1$)  $(\mu_E)_0 = \mu_E$ — since $\mu_E = 2^{\mu_E}$ and 2 is the $0^{\text{th}}$ prime, $p_0$.

2.1.3. Let **E** be a non-trivial deduction, then its root, $\nu E$, is the unique member of $\mathscr{T}_E$ which generates the branch, growing from $\nu E$, which coincides with $\mathscr{T}_E$; also, the lowermost node

$$\frac{E_1, \ldots, E_s}{\nu E}$$

of $\mathscr{T}_E$ is unique.

Let, for $i = 1, \ldots, s$, the branch $\mathscr{b}\nu_{E_i}$ of $\mathscr{T}_E$ generated by the *beginning*, $E_i$, of this node be assigned the En-N $\mu_{E_i}$. Then the En-N, $\mu_\mathbf{E}$, of the object $\mathscr{T}_E$ as well as the deduction it represents are both assigned the En-N

$$2^{\mu_E} \cdot \prod_{g=1}^{s} p_g^{\mu_{E_i}},$$

where $p_g$ denotes the $g^{\text{th}}$ prime. ($p_0 = 2$, $p_1 = 3$, $p_3 = 5$, $p_4 = 7$, ...).

Clearly, also for this case the equality (or deduction (proof)-root-correlation)

($\nu 1$)  $(\mu_\mathbf{E})_0 = \mu_E$

must hold.

Since any subtree, $\mathscr{T}_{E^*}$, of $\mathscr{T}_E$ is also a deduction tree with some root formula, $\mathscr{r}E^*$, the above described assignment works for any branch of $\mathscr{T}_E$. In more detail: Each node of a deduction tree shall determine the branch, $\mathscr{br}_E$, of the tree generated by the conclusion, E, of that node. Let $D_1, \ldots D_r$ be all premisses of that node (which shall always be considered in a definite order), and let, for $g = 1, \ldots, r$, $\mathscr{br}_{D_g}$ denote the branch of the tree generated by the premise $D_g$, and a number $\mu_{D_g}$, be given as a En-N of $\mathscr{br}_{D_g}$. Then

$$\mu_{\mathscr{br}_E} = 2^E \cdot \prod_{g=1}^{r} p_g^{\mu_{D_g}},$$

shall be assigned as the En-N to $\mathscr{br}_E$.

Thus, each branch of a deduction tree shall be assigned a En-N; the tree itself shall be assigned the number of the branch generated by its root, and this branch shall be identified with the tree, and this number shall also be considered as the number of the deduction represented by the tree.

2.1.4. For the special case, $s = 2$, when the lowermost node is the *modus ponens* with the premises $E_1$ and $E_1 \supset E$, and — $\mu_{E_1}$ and $\mu_{E_1 \supset E}$ being, respectively, the En-N's of the branches of these premisses, the branch growing from the conclusion, E, shall be assigned the En-N

$$Mp(\mu_{E_1}, \mu_{E_1 \supset E})$$

which equals, in accordance with this stipulation, the product

$$\mu_E = 2^E \cdot 3^{\mu_{E_1}} \cdot 5^{\mu_{E_1 \supset E}} ;$$

The binary function symbol Mp thus introduced is p.r. and it will be included in LAri.

In this way each branch of $\mathscr{T}_E$ is assigned a En-N and a course-of-values induction on the height of $\mathscr{T}_E$, eventually, assignes each deduction, **E**, a unique En-N, $\mu_E$.

2.1.5. For the case when a deduction, **E**, is carried out in the *open* extension $Ari^\nu$ of Ari, the rule **Gen** (see p.32, below) gives, for a formula $E(\mu)$ containing the variable μ freely, the conclusion μE(μ), denoted by E, and the En-N , $\mathscr{E}_E$ , of the object $\mathscr{T}_E$ shall equal

2 $^{\mu E(\mu)}$ . 3 $^{E(\mu)}$ ,

where   is used as a special notation to indicate En-N of *deduction in the open logic*, **lo**. The deduction-root correlation

$$(\mathscr{E}_E)_0 = \mu E(\mu) \quad (= \mathscr{E}_E)$$

2.1.6. As soon as all top formulae, $E_1$, ..., $E_m$ , of a deduction $\mathscr{E}_E$ are *axioms* (accepted for Ari or its extensions in this work), then $\mathscr{E}_E$ shall also be called a *proof* of its root formula .

**Remarks**. The use of the notation $(n)_g$ in the deduction (proof)-root-correlation, ($\varkappa_1$), above, is taken from [KLEENE, 1952], §45; in cases when the use of functional notations is preferable it shall be denoted by $\exp_g$:

$\exp_g = (n)_g$ , (g = 0, 1, 2, ...).

These $\exp_g$ shall be dealt with as the unary function symbols, i.e. g shall be fixed. Also the iterated notations

$(n)_{g,h}$  for  $(n)_g)_h$

$(n)_{g,h,i}$  for  $(n)_g)_h)_i$

are borrowed from [KLEENE, 1952] as are the notations for the corresponding unary function symbols whose functional notations shall be specified as

$\exp_{g,h}$; $\exp_{g,h,i}$, etc.

For the benefit of the reader unfamiliar with [KLEENE, 1952], we review, below, these Kleene notations for the values of the exponents of the prime factors of a Gaussian factorization together with some equalities that often occur throughout this work on (and on the basis of which the deduction (proof)-root-correlation is established in Ari). Readers familiar with the Kleene notation for prime exponents may safely skip remark **A**. Also the following remark **B** may be skipped since its main points are restated (with less details) in Section 4, below.

**A**. Thus, to begin with, unless $\mathbf{m} = 0$, $\exp_i(\mathbf{m})$ denotes the exponent of the $i^{th}$ prime factor among the prime factors of $\mathbf{m}$. If $\mathbf{m} = 0$, then also $\exp_i(\mathbf{m}) = 0$. We may write $((\mathbf{m})_i)_j$ as $(\mathbf{m})_{i,j}$, $((\mathbf{m})_i)_j)_k$ as $(((\mathbf{m})_{i,j,k}$ etc.

Now, $Mp(\mathbf{n},\mathbf{q})$ can be introduced in terms of its arguments, $\mathbf{n}$, $\mathbf{q}$, as satisfying the equality

$$Mp(\mathbf{n},\mathbf{q}) = 2^{(\mathbf{q})_{0,2}} \cdot 3^{\mathbf{n}} \cdot 5^{\mathbf{q}}$$

(cf. [KLEENE, 1952], §52).

Thus., for $Mp(\mathbf{n},\mathbf{q}) = \mathbf{M}$, the following equalities must hold;

(Mp.1)  $(\mathbf{M})_0 = (\mathbf{q})_{0,2}$

(Mp.2)  $(\mathbf{M})_1 = \mathbf{n}$,

(Mp.3)  $(\mathbf{M})_2 = \mathbf{q}$.

Thus, for $Mp(\mathbf{n},\mathbf{q}) = \mathbf{M}$, the equality, above, can be replaced by its equivalent

$$Mp(\mathbf{n},\mathbf{q}) = \mathbf{M} = 2^{(\mathbf{M})_0} \cdot 3^{(\mathbf{M})_1} \cdot 5^{(\mathbf{M})_2}.$$

The En-N of an implication $E_1 \supset E_2$ is given as; $2^3 \cdot 3^{E_1} \cdot 5^{E_2}$. This gives

(Im.1)  $(E_1 \supset E_2)_0 = 3$

(Im.2)  $(E_1 \supset E_2)_1 = E_1$

(Im.3)  $(E_1 \supset E_2)_2 = E_2$.

Thus, for $Mp(\mathbf{n},\mathbf{q}) = \mathbf{M} = Mp(\mu_{E_1}, \mu_{E_1 \supset E_2})$, these equalities and the deduction-root-correlation give the follöwing equalities

$(\mathbf{M})_2 = (Mp(\mathbf{n},\mathbf{q}))_2 = (Mp(\mu_{E_1}, \mu_{E_1 \supset E_2}))_2 = \mu_{E_1 \supset E_2}$;

$(\mathbf{M})_{2,0} = E_1 \supset E_2$ — and, further,

(Mp.4)    $(\mathbf{M})_{2,0,0} = (\ulcorner E_1 \supset E_2 \urcorner)_0 = 3$

(Mp.5)    $(\mathbf{M})_{2,0,1} = (\ulcorner E_1 \supset E_2 \urcorner)_1 = \ulcorner E_1 \urcorner$

(Mp.6)    $(\mathbf{M})_{2,0,2} = (\ulcorner E_1 \supset E_2 \urcorner)_2 = \ulcorner E_2 \urcorner$.

From the last three lines of equalities and $(\mathbf{M})_2 = \mathbf{q}$ one obtains;

(mp.1)    $(\mathbf{q})_{0,0} = 3$

(mp.2)    $(\mathbf{q})_{0,1} = \ulcorner E_1 \urcorner$

(mp.3)    $(\mathbf{q})_{0,2} = \ulcorner E_2 \urcorner$

and, therefore, by (Mp.6) and (mp.3), the exponent $(\mathbf{q})_{0,2}$ in **M** equals both $\ulcorner E_2 \urcorner$ and $(\mathbf{M})_{2,0,2}$ : this and (Mp.1) give

(Mp.7)    $(\mathbf{M})_0 = (\mathbf{M})_{2,0,2}$.

If $ro$ denotes the unary function symbol which extracts the En-N of the root formula from the En-N of non-trivial deduction (of a MP-step), then $ro(\mathbf{M})$ equals $(\mathbf{q})_{0,2}$ since $(\mathbf{q})_0 = \ulcorner E_1 \supset E_2 \urcorner$ yields $(\mathbf{q})_{0,2} = \ulcorner E_2 \urcorner$, and, therefore, $ro(\mathbf{M})$ shall equal $(\mathbf{M})_0$ $(= (\mathbf{M})_{2,0,2} = (\mathbf{q})_{0,2})$. Thus, the equality

($rt$)    $ro(\mathbf{M}) = (\mathbf{M})_0$

shall hold for any **M** which equals En-N of a 'MP-deduction' (i.e. deduction in Ari the last step of which is a MP-step) in Ari or its (open) extension Ari$^\nu$.

(End of remark **A**).

**B**. Actually, the function symbol Mp shall be replaced in this work by another function symbol, mp, which is introduced by a slightly more complex definition than for Mp. The purpose for using this new definition is that it provides for the conditions for the applicability of Mp which are 'factored' into the new function symbol mp to be defined below.

We recall the essential features of a MP-step. *Viz.*, for the case when the deduction has En-N **M** = Mp(**n**,**q**), the formula C which occurs as the minor premise of the last MP-step must be the 'same' whether it is considered as the root of the deduction which has **n** as its En-N or as the antecedent of the root of the other deduction having **q** as its En-N. That is, the En-N of this formula shall be the same in both cases. This En-N shall equal, in these two cases, the numbers

$$(\mathbf{M})_{1,0} \text{ and } (\mathbf{M})_{2,0,1}$$

and therefore these numbers shall be equal. The equality

$$(\mathbf{M})_{2,0,1} = (\mathbf{M})_{1,0}$$

must be satisfied as a precondition of the 'normal' applicability of the operation MP represented by Mp.

One way of formalizing this condition at a MP-step is by using a new function symbol, mp, which shall require for its definition, besides Mp, the additional function symbols ⬚ and sg together with $\overline{sg}$ (the inverse of the signum function, sg) and msd ('modulus of symmetric difference' or 'absolut difference') and the already used symbols $\exp_0$, $\exp_{0,0}$ and $\exp_{0,1}$ so that mp shall be given by the following explicit definition

**def**$_{mp}$:   mp(**n**,**q**)   Mp(**n**,**q**) · ⬚(**n**, **q**) · sg(**n**),

where the symbol ⬚ is defined by

**def** ⬚:   ⬚(**n**,**q**)
$\overline{sg}$ (msd ($\exp_{0,0}(\mathbf{q})$,3)) · $\overline{sg}$ (msd ($\exp_{0,1}(\mathbf{q})$,$\exp_0(\mathbf{n})$))

Thus, the factor $\overline{sg}$ (msd ($\exp_{0,1}(\mathbf{q})$,$\exp_0(\mathbf{n})$)) in the definition **def**$_{mp}$ of mp, formalizes a condition of normal applicability of the operation MP (as described above) iff it equals one.

Furthermore, another such condition shall consist in that the the major (right) premise of the MP-step must be an implication — and the factor $\overline{sg}$ (msd ($\exp_{0,0}(\mathbf{q})$,3)) formalizes this iff it equals one in **def**$_{mp}$.

Lastly, the condition sg(**n**) = 1 in **def**$_{mp}$ assures that **n** is a *positive* integer, so that, in particular, 0 shall not be En-N of any proof or deduction.

It is clear that unless mp(**n**,**q**) equals zero, its value will coincide with Mp(**n**,**q**).

(This remark **B** will be restated by a more compact formulation in Section **4**., below).

(End of remark **B**).

2.2. As to the En-N's of p.r. functions, we remark that **En** shall deal with the function *symbols* rather than the functions they denote since the identification of functions does not always have a clear finitistic meaning. Thus, the usual way of assigning En-N's to function symbols consists in that En-N's are first assigned to the initial p.r. function symbols and any further function symbol is assigned an En-N in accordance with the p.r. description ([KLEENE, 1952], §43) which introduces (or 'defines'] this function symbol. Since, in this work, there will be no need to identify any two function symbols introduced by different p.r. descriptions and while their total number shall not be very large (say, << 500), it will be possible just to present these En-N's in a Table, Ta, mentioned on p..There will be no need to present this table now: moreover, for the sake of generality this presentation will be postponed as long as possible. In this way, with regard to the vocabulary and En-N's of function symbols, the language LAri, **En** and the set of axioms of Ari will remain partially underdefined and admit various versions of specifications. Different tasks may require some freedom of choice between such versions — and it will suffice to fix the table Ta and the list of axioms of Ari and the introductsion of the function symbol ax (and thus of the final specification of the symbols $\ell$' and $\ell$' (discussed in the Preface, p.) just till the moment when such a choice will be done. If there really will be a need in different specifications of LAri, then the use of the notations Ari, ax, $\ell$' and $\ell$ will be specified by using different subscripts for each such choice. For the general purposes of this work, **En** shall be supposed to contain the En-N's $71^1$ ,..., $71^{ta_1}$, $73^1$ ,..., $73^{ta_2}$, $79^1$ ,..., $7^{ta_3}$ without any gaps and the En-N's $_{suc}$, $_+$, $_.$ specified as per p. , above.

There will be no need to associate these En-N's with the p.r. descriptions of these functions denoted by their symbols. The use of such notations for these En-N's as

, , , — or, , ,

will suffice for all purposes of this work.

The equalities

$$= \overset{\nu}{ta_1}, \quad = \overset{\nu}{ta_2}, \quad = \overset{\nu}{ta_3}$$

will be supposes to hold and be Ari-provable (as soon as the integers $ta_1$, $ta_2$, $ta_3$ are denoted by the Arabic $< 500$), and the calculability of the function, ex, of exponentiation shall entail the calculability of these En-N's, i.e. the Ari-provability of the equalities $= \overset{\nu}{ta_1}, = \overset{\nu}{ta_2}, = \overset{\nu}{ta_3}$.

3. The *logical axioms* fall into five *groups*, I.-V., for which the schemata of their matrices will now be given. Each logical axiom shall be a closure of such a 'matrix', and the matrix may be any formula in LAri obtainable as an instance of one of these schemata subject, for some of these schemata, to restrictions on the use of the variable indicated in the presentation of the schemata by $\bar{\mu}$. The bold letters **A**, **B**, **C** of these presentations must be replaced, throughout each schema of group I., by any formula in LAri (the same for each occurrence of each of these bold letters), in order to obtain the (arbitrary) instance of the matrix schema. The letter $\ell$ denotes the closed formula $0 = 1$ (called the '*anti-axiom*'). The labels on the left, below, are also, when *unbolded* and *without parenthesis* used as names for the corresponding *instances* of the axiom matrices in the running (right) column of proof comments in the proof 1.- 686a., pp. 54-70, as well as in Appendices B, C and D.

### I. Propositional Axiom Schemata Matrices.

(**Imp**1):   **A** → (**B** → **A**),

(**Imp**2):   (**A** → (**B** → **C**)) → ((**A** → **B**) → (**A** → **C**)),

(**Imp**3):   $\ell$ → **B**,

(**Con**1):   **A** & **B** → **A**,

(**Con**2):   **A** & **B** → **B**,

(**Con**3):   **A** → (**B** → **A** & **B**),

(**Dis**1):   **A** → **A** ∨ **B**,

(**Dis**2):   **B** → **A** ∨ **B**,

(**Dis**3):   (**A** → **B**) & (**A** → **C**) → (**A** ∨ **B** → **C**)

## II. Quantifier Axiom Schemata Matrices.

Below, $\bar{\mu}$, denotes an arbitrarily choosen variable in LAri, **A** denotes any formula in LAri in which $\bar{\mu}$ does not occur freely and $C(\bar{\mu})$ — any formula in LAri in which $\bar{\mu}$ does occur freely. In the schemata (**Fi1**),(**Fi2**) $\bar{\mu}$ must be fictitious while in the the schemata (**Fi3**),(**Fi4**) $\bar{\mu}$ must be ficticious. The schema (**Dis**) has no such restriction.

(**Fi1**): $\quad\bar{\mu}\mathbf{A} \quad \mathbf{A}$,

(**Fi2**): $\quad \mathbf{A} \quad \bar{\mu}\mathbf{A}$,

(**Fi3**): $\quad \mathbf{A} \quad \bar{\mu}\mathbf{A}$,

(**Fi4**): $\quad \bar{\mu}\mathbf{A} \quad \mathbf{A}$,

(**Dis**): $\quad \bar{\mu}(B \quad C) \quad (\bar{\mu}B \quad \bar{\mu}C)$,

(**WBA**): $\quad \bar{\mu}C(\bar{\mu}) \quad C(\bar{\mu})$,

(**WBA**): $\quad C(\bar{\mu}) \quad \bar{\mu}C(\bar{\mu})$,

(**BR**): $\quad \bar{\mu}(C(\bar{\mu}) \quad \mathbf{A}) \quad (\bar{\mu}C(\bar{\mu}) \quad \mathbf{A})$.

## III. Equality Axiom Schemata Matrices.

Below, r, t, s are any termoids, P any predicate symbol and any function symbol in LAri in which r, s occupy the same argument places in the consequent of (LEA - rp). The schema (Ref ) may be choosen as a nonlogical axiom schema.

(**Ref**): $\quad r = r$,

(**Ref**): $\quad r \quad r$,

(**Sym$_=$**): $\quad r = s \quad s = r$,

(**LEA$_1^=$**): $\quad r = s \quad (P(r,t) \quad P(s,t))$,

($LEA_2^=$):   r = s   (P(t,r)   P(t,s)),

($LEA$-rp): r = s   ( ... r ... ) =   ( ... s ... ).

## IV. Equality-Substitution Axiom Schema Matrix.

Below, t denotes any termoid not containing $\bar{\mu}$ freely. The arrow $\bar{\mu}$ t is a condition which expresses that t's value(s), of which there has to be at least one, belongs (belong) to the range of $\bar{\mu}$ (which, thus, must be assumed to be non-empty). The small inverted turn-stile sign, ⊣, is used to indicate that the arrow expresses, in MAri, a condition which makes it 'correct' and under which the formula to the right of ⊣ is to be considered as accepted as an axiom. This is a metatheoretical notation which does *not* belong to any of the formal languages in this work; in the continuation of this work the arrows which 'accompany' some axioms will be collected, for each formal proof which contains these axioms, in order to explain or prevent arising contradictions.

(**E-S**):    $\bar{\mu}$   t ⊣   $\bar{\mu}$ ( $\bar{\mu}$ = t)

## V. Strong Bernays Axiom Schemata Matrices.

Below, $\bar{\mu}$ occurs freely in **C**($\bar{\mu}$) and with the same restrictions on $\bar{\mu}$ and t as in (E-S), above, *i.e.*: t is free for $\bar{\mu}$ in **C**($\bar{\mu}$) (cf. § 18, [KLEENE, 1952]) — i.e. that the occurrences in C(t) of variables via the occurrence of t which don't stem from occurrences (if any) of t in.**C**($\bar{\mu}$) which are free in t continue to be free also in C(t); the notation denotes the formula obtainable from **C**($\bar{\mu}$) as the result of replacing each occurrence of $\bar{\mu}$ by the occurrence of t.

(**SBA** ):    $\bar{\mu}$   t ⊣      $\bar{\mu}$ **C**($\bar{\mu}$)    **C**(t),

(**SBA** ):    $\bar{\mu}$   t ⊣      **C**(t)    $\bar{\mu}$ **C**($\bar{\mu}$).

**Remark**. These schemata are derivable from the equality-substitution schema (with the same arrows) while the latter is derivable from these schemata; hence, they are interchangable so that it suffices to choose either the former or the latter for quantification theory. (End of remark).

This ends the presentation of the logical axiom schemata belonging to the *Strong Predicate Calculus*. For some comments on the limitations of the Weak Predicate Calculus, see [YESSENIN-VOLPIN, 1981] or ENDNOTE B.

**3**.1.  The **rule of inference** shall, in the *closed* logic, **Lo**, be, as usual, only the rule of *modus ponens*;

(MP)    $\dfrac{\mathbf{B},\ \mathbf{B}\ \mathbf{C}}{\mathbf{C}}$ ;

nevertheless, *derived rules* of inference will be used often to shorten proof presentations (see Appendix B for the list of the derived rules of inference). However, in the *open* logic, **lo**, the further rule of *generalization*

(**Gen**)    $\dfrac{\mathbf{C}(\bar{\mu})}{\bar{\mu}\ \mathbf{C}(\bar{\mu})}$

is admitted so that if $\mathbf{C}(\bar{\mu})$, in which $\bar{\mu}$ occurs freely, has an Ari-proof (or Ari-deduction whose hypotheses do not contain $\bar{\mu}$ freely), then also $\bar{\mu}\ \mathbf{C}(\bar{\mu})$ has an Ari-proof (or Ari-deduction with the same hypotheses) obtainable by an application of **Gen** whose premises is the root of that Ari-proof (-deduction). The formula denoted, in these schemata, by **B**, **C**, $\mathbf{C}(\bar{\mu})$ can be choosen in LAri arbitrarely provided that $\bar{\mu}$ occurs freely in $\mathbf{C}(\bar{\mu})$.

If the logic is closed, then the rule (**Gen**) cannot actually be used in proofs constructed in accordance with its rules, and the rule (MP) can be used only with closed **B**, **C** in its instances; the admissability of not closed **B, C** or of (**Gen**) with variables other than $\bar{\mu}$ and free in $\mathbf{C}(\bar{\mu})$ in deductions can be tolerated but such cases will never occur below. Without any harm to this work, **B** and **C** in instances of (MP) can be formally restricted to be closed when the logic is **Lo**.
(The open logic consists of the above groups of axioms *except* that the axioms shall not be supposed closed and, hence, be presented *without* closure strings, i.e. the schemata given above for matrices can be considered as the (logical) axiom schemata.

3.2.  Among the non-logical axioms there shall be only a few *instances* of *tertium non datur* (**TND**) and *double negation elimination* (**DNE**) while we don't consider these axioms in their full generality and, hence, abstain from listing the corresponding schemata. There shall also be a single instance of *definite descriptions* for which the axioms given in [ROSSER, 1953], chapter VIII, axiom schemata 8.-11., shall be accepted. The axioms governing Ari's proof property contain some unusual clauses which, however, lead to some metatheoretical simplifications. Ari's remaining non-logical axioms shall be restricted to those actually used in this work (including Part II) — which are the usual arithmetical axioms for p.r. arithmetic extended by those special axioms for termoids and function symbols needed for the present work. They will be introduced in the contexts in which they are used. In most cases, we will, for convenience, refer the reader to [KLEENE, 1952].

In any case, each non-logical axiom in Ari other than (**TND**) and (**DNE**), the Peano axioms, induction or recursion axioms or axiom-definitions, shall be a formula provable in the arithmetic of [KLEENE, 1952].

Ari's axioms shall form a recursively enumerable set.

4. We now introduce, in MAri, the *axiom-definitions* for the respective characteristic function symbols, $b_1$, $b$, which '**En**-represent' Ari's *proof property*, $\mathcal{B}_1$, and *proof relation*, $\mathcal{B}$, respectively. Other (equivalent) axiom-definitions are also possible although the present ones have been found to be the most convenient in the present context. Furthermore, observe that the presentation of these axiom-definitions is first given in MAri while in Section 4.3., below, the corresponding axioms in Ari will be introduced.

This introduction shall require the definitions of the already mentioned p.r. function symbols $ax_{Ari}$ and mp, respectively, which '**En**-represent' the property of being an Ari-axiom and a conclusion of a MP-step, respectively.

The unary $ax_{Ari}$ shall **En**-represent the property of being an Ari-axiom; only members of a finite (or, at most, recursively enumerable) set of closed formulae in LAri shall be included in this property; thus, it will be possible to include, in Ari, all axioms required for the introduction of the function symbol $ax_{Ari}$. However, the p.r. character of the function it introduces will remain problematic untill the provability of all required instances of (T-Ax) — see p.— will be, eventually, established. Besides that, it will suffice that a p.r. function, , enumerate all those **n**'s for which the equality $ax_{Ari}(\mathbf{n}) = 0$ holds. See more on p.

The binary function symbol mp will be introduced, independently of $ax_{Ari}$, in order to prepare the introduction of the unary function symbol $b_1$ En-representing the proof property of the system Ari.

The Ari-proofs shall be 'trivial' and 'non-trivial' (see Sections 2.1.2.-2.1.3., above); the binary function symbol mp was already 'described' in Section 2.1.4. and Remark B at the end of that section, above, and it will be used in order to 'En-represent' the property $\text{Ant}_2^{b_1}$ of being a non-trivial proof. More specifically, the definitions (see p. 28, above)

**def**$_{mp}$:  $mp(\mathbf{n},\mathbf{q})$    $Mp(\mathbf{n},\mathbf{q}) \cdot$    $(\mathbf{n},\mathbf{q}) \cdot sg(\mathbf{n})$

**def** :    $(\mathbf{n},\mathbf{q})$

$\overline{sg}(msd(exp_{0,0}(\mathbf{q}),3)) \cdot \overline{sg}(msd(exp_{0,1}(\mathbf{q}),exp_0(\mathbf{n})))$,

are given in MAri and use the *parameters* **n**, **q** which range over (all) non-negative integers. (The use of such parameters will be continued in other definitions in MAri).
(For any termoids r, s, msd(r,s) denotes the absolut difference or 'modul of symmetric difference' of r, s — for which the standard notation usually is $|r - s|$; $\overline{sg}$ denotes the function symbol inverse to sg ($\overline{sg}(\mathbf{n})$ can be introduced as $1 \dotminus sg(\mathbf{n})$ if sg is introduced, explicitly, by cases of $\mathbf{n} = 0$, $\neg \mathbf{n} = 0$, and if $\overline{sg}$ is thus explicitly introduced by cases then $sg(\mathbf{n})$ can be introduced as $\overline{sg}(\overline{sg}(\mathbf{n}))$.)

As it was stated by the end of **Remark B** in Section 2.1., above, unless $mp(\mathbf{n},\mathbf{q}) = 0$, the two functions Mp and mp coincide, i.e. have the same values; that is achieved by including the factor $sg(\mathbf{n})$ in **def**$_{mp}$. The right side of **def** is the product of two factors; the first of them makes $mp(\mathbf{n},\mathbf{q})$ equal to zero in the case when **q** equals En-N of deduction whose end formula is not an implication while the second factor of this product makes $mp(\mathbf{n},\mathbf{q}) = 0$ in the case when **n**, **q** are En-N's of two deductions whose roots cannot be used, respectively, as the minor and major premises of a MP because the root of the minor premise does not coincide with the antecedent of the root of the major premise. The equality

$mp(\mathbf{n},\mathbf{q}) = Mp(\mathbf{n},\mathbf{q})$

must hold iff $(\mathbf{n},\mathbf{q}) = 1 = sg(\mathbf{n})$.

4.1. Continuing in MAri (as indicated by the use of parameters instead of variables), the axiom-definitions for $b_1$ shall be:

**df**$_1^{b_1}$ :    $b_1(\mathbf{n}) = 0$   if   $\mathbf{n} = 2^{(\mathbf{n})_0}$   &   $ax_{Ari}((\mathbf{n})_0) = 0$ ,

$\mathbf{df}_2^{b_1}:$  $\quad b_1(\mathbf{n}) = 0 \quad \text{if} \quad \mathbf{n} = \text{mp}\,((\mathbf{n})_1, (\mathbf{n})_2) \quad \& \quad \neg\,\mathbf{n} = 0 \quad \&$

$$b_1((\mathbf{n})_1) = 0 \quad \& \quad b_1((\mathbf{n})_2) = 0,$$

$\mathbf{df}_3^{b_1}:$  $\quad b_1(\mathbf{n}) = 1 \quad \text{if} \quad \neg\,(\text{Ant}_1^{b_1}(\mathbf{n}) \vee \text{Ant}_2^{b_1}(\mathbf{n}))$

where $\text{Ant}_i^{b_1}(\mathbf{n})$ denotes the 'antecedents' in $\mathbf{df}_i^{b_1}$, $i = 1, 2$.

It will also be convenient to define $\text{Ant}_3^{b_1}(\mathbf{n})$ as a compact notation for the formula

$$\neg\,(\text{Ant}_1^{b_1}(\mathbf{n}) \vee \text{Ant}_2^{b_1}(\mathbf{n})).$$

4.2. $b$ can now be introduced by the following definition:

$\mathbf{def}_1^{b}:$  $\quad b(\mathbf{n},\mathbf{q}) = 0 \quad b_1(\mathbf{n}) = 0 \quad \& \quad (\mathbf{n})_0 = \mathbf{q},$

$\mathbf{def}_2^{b}:$  $\quad b_1(\mathbf{n}) = 0 \quad \& \quad (\mathbf{n})_0 = \mathbf{q} \quad b(\mathbf{n},\mathbf{q}) = 0,$

$\mathbf{def}_3^{b}:$  $\quad \neg\,b(\mathbf{n},\mathbf{q}) = 0 \quad b(\mathbf{n},\mathbf{q}) = 1.$

4.3. In order to introduce the corresponding axiom-definitions, as above, in Ari, it suffices to replace the distinct parameters by distinct variables and thereafter close each formula by a suitable string of universal quantifiers (for all axiom-definitions must be axioms in Ari and therefore be obtained as *closed formulae*). Thus, Ari will include the following axiom-definitions for $b_1$.;

$\text{df}_1^{b_1}:$  $\quad \mu(\text{Ant}_1^{b_1}(\mu) \quad b_1(\mu) = 0);$

$\text{df}_2^{b_1}:$  $\quad \mu(\text{Ant}_2^{b_1}(\mu) \quad b_1(\mu) = 0);$

$\text{df}_3^{b_1}:$  $\quad \mu(\text{Ant}_3^{b_1}(\mu) \quad b_1(\mu) = 1)$

where $\text{Ant}_i^{\mathcal{b}_1}(\mu)$, i = 1, 2, 3, is defined as $\text{Ant}_i^{\mathcal{b}_1}(\mathbf{n})$, i = 1, 2, 3, Section 4.1., above, except that $\mathbf{n}$ is replaced throughout by $\mu$ — and the following axiom-definitions for $\mathcal{b}$ ;

$\text{def}_1^{\mathcal{b}}$:    $\mu$    ( $\mathcal{b}(\mu, ) = 0$    $\mathcal{b}_1(\mu) = 0$ & $(\mu)_0 =$    )

$\text{def}_2^{\mathcal{b}}$:    $\mu$    ( $\mathcal{b}_1(\mu) = 0$ & $(\mu)_0 =$    $\mathcal{b}(\mu, ) = 0$ ),

$\text{def}_3^{\mathcal{b}}$:    $\mu$    ( $\neg\, \mathcal{b}(\mu, ) = 0$    $\mathcal{b}(\mu, ) = 1$ ).

In all of these axiom-definitions $\mu$, may denote any distinct variables.

**5**. In this section we define first a closed extension, $\text{Ari}^+$, of Ari and then an *open* extension, $\text{Ari}^\nu$, of $\text{Ari}^+$, by an extended use of Kleene's notation for enumerated formulae introduced in [Kleene, 1952], §42. Here $\text{Ari}^+$ denotes a *closed* extension of Ari, obtained from Ari by adding, for each Ari-axiom E, the implication

(t-ax):    $\text{ax}_{\text{Ari}}(\nu_E) = 0$    E

to the axioms of Ari. (Here the (special) notation $\nu_E$ denotes the numeroid of the En-N $_E$ of E). The implication t-ax is Ari-provable and the proof of this statement is given in Appendix D1.

It is also possible to consider $\text{Ari}^+$ as obtained from Ari by adding, for each axiom, E, in Ari and any closed formula D in LAri, the formula D    E    as the new axiom of Ari (of which t-ax is a special case). In either case, $\text{Ari}^+$ shall have the same language, LAri, as Ari. The purpose of any of these extensions is to aid in the substantiation of the schema T-Ax which will be concluded in Part II of this work.

In Section 5.1., below, the set of formulae belonging to $\text{LAri}^\nu$ are defined and in Section 5.2. $\text{Ari}^\nu$'s additional axioms are introduced including a mention of the critical schema T-Ax. Section 5.3 is concerned with the substantiations of the new axioms belonging to $\text{Ari}^\nu$ (except T-Ax) as well as considerations of the important mp-induction formula. Section 5.4. concerns the possibility of considering M-    not as an axiom but as an Ari-theorem. The section ends with a summary overview of the proof, in Ari, of the formula $\text{Con}_{\text{Ari}}$.

**5**.1. The unary characteristic p.r. function symbol  cfor  shall En-represent, in LAri, the property of being a *closed* formula in LAri. (That is, cfor($n$) = 0 if $n$ equals En-N of such formula and cfor($n$) = 1 otherwise). The unary p.r. function symbol  fl  shall be introduced by the definition

> fl($n$)   $n$     if $n$ equals En-N of closed formula in LAri (*i.e.*, if cfor($n$) = 0),
>
> fl($n$)   $v_{\ell}$ , otherwise.

In MAri, for any En-N,  r , of closed formula, E , in LAri,

> $F_r$  shall denote  E

and for any other interger  r

> $F_r$  shall denote (the 'anti-axiom') $\ell$.

This definition is borrowed from [Kleene, 1952], §42, except that the use of this notation is now extended also to cases when the integer r is not En-N of closed formula in LAri. Below, the use of r as a subscript in the notation $F_r$ must be restricted to termoids in LAri (which may contain free variables).

The notations $F_r$ (with various termoids  $r_1$ , $r_2$ , s , t , u ..., in LAri, admitted as their subscripts r ) shall be used below in derivational presentations, i.e. in the presentations, in MAri, of derivations (deductions and proofs) in Ari$^+$. These presentations considered as such will not necessarily be derivations in Ari$^+$ — in order to obtain the latter, each $F_r$ must be replaced by the formula, in LAri, whose En-N equals  fl(r). When the termoid  r  is a constant then it shall be considered, in MAri, as (denoting) a quite definite integer and also fl(r) will denote a quite definite integer, and the formula whose En-N equals fl(r) will be quite definite. (The definite desciptions in LAri shall be considered as having quite definite integers as their values specified as 0 when otherwise these descriptions have no definite value; remark that the 'definiteness' of the values does not presuppose their 'calculability' — the dependence on TND with non-restricted quantifiers can make the values problematic). When the termoid  r  contains free variables, these variables, in all their free occurrences, will be considered as  the variables held constant ([Kleene, 1952], §23) throughout the derivational presentation in which $F_r$  occurs; that is, they will be considered as the parameters of these presentations. For each fixed value(s) of the parameter(s), $F_r$  shall denote a quite definite formula, *viz*., the closed formula in LAri whose En-N equals fl(r). (In particular, when fl(r) = $\nu_\ell$, the closed formula shall be $\ell$.)

Below, the language LAri$^\nu$ will be used for presentaions of derivations in LAri. These presentations may contain (distinct) parameters — and they shall, with a few (explicit) exceptions, be replaced in LAri$^\nu$ by (the distinct) free variables to be held constant throughout each derivational presentation given in LAri$^\nu$.

**5**.2.  LAri$^\nu$ shall contain LAri and, further, Ari$^\nu$ shall be an extension of each Ari, Ari$^+$. This extension shall consist in the adjoining, on the level with the atomic formulae, of the notation $F_r$ (where r may be any termoid in LAri). That is, the formation rules for LAri$^\nu$ can be described just as these rules for LAri, the only change being that the concept of atomic formula is extended so that also all notations $F_r$ shall be considered, in LAri$^\nu$, as atomic formulae; the variables free (not free) in r shall be considered as occurring, in $F_r$, freely or not freely. Thus, Ari$^\nu$ is, as far as its closed formulae are concerned, basically a notational variant of Ari$^+$ (but with the important difference that Ari$^\nu$ provides for variable numeroids).

LAri$^\nu$ shall be designed for presentations of Ari-proofs of the consistency formula, Con$_{Ari}$ — or, more specifically Con$_{Ari}$, where     denotes the variable bound in this formula,

Con$_{Ari}$ ;    $\neg$    $\ell(\ ,\nu_\ell) = 0$;

besides that, LAri$^\nu$ may be used for presentations of other deductions in Ari, including the deduction of Con$_{Ari}$ from the μ-closure of the hypothesis

(T-Ax):         $ax_{Ari}((\mu)_0) = 0$     $F_{(\mu)_0}$

which can be used as a 'temporary' axiom in the extension Ari$^\nu$ of Ari. Only in Part II of this paper will the definite status of T-Ax be considered so that in the present Part I it will play no role. Here we prove a consequence of T-Ax, namely formula 7. in the proof 1.-686a., Section 6, below. This proof is given in Appendix D2., below.

However, considerable interest is attached to the schema T-Ax since the implication

$$\mu(\text{T-Ax}) \supset \text{Con}_{\text{Ari}}$$

will be found Ari-provable by virtue of the Deduction Theorem.

5.2.1. Each axiom in Ari$^+$ shall be an axiom in Ari$^\nu$ and each instance in LAri$^\nu$ of axiom scheme of Ari$^+$ shall be — if closed — considered as an axiom in Ari$^\nu$. Besides that, Ari$^\nu$ shall contain, as axioms, the closures of the following three '$\nu$-schemata'

(LEA$^\nu_1$): $\quad r_1 = r_2 \supset (F_{r_1} \supset F_{r_2})$,

(LEA$^\nu_2$): $\quad \text{fl}((\mu)_0) = (\mu)_0 \supset {^\nu}F_{\text{fl}((\mu)_0)} = {^\nu}F_{(\mu)_0}$,

(LEA$^\nu_{\text{MP}}$): $\quad (r)_0 = 3 \supset (F_r \supset (F_{(r)_1} \supset F_{(r)_2}))$.

The LAri$^\nu$-instances of these 'matrices' are obtainable by replacing throughout these schemata (including the subscripts) the letters

$r_1, r_2, r, \mu$

by any termoids in LAri.

5.3. The letters $r_1, r_2, r$ can be considered, in these schemata, as standing for any termoids in LAri; accordingly, all function symbols in the termoids $r, (r)_1, (r)_2, r_1, r_2$ shall be p.r.: in the cases at issue, also no $\iota$-description occurs in these termoids. Any replacing, throughout these schemata, of all free occurrences of variables by numeroids shall give a *numeroidal instance* of the $\nu$-scheme which does not contain $\iota$-descriptions.

In order to prove that, in MAri, let, **r**$_1$ and **r**$_2$ be the termoids obtainable from the termoids $r_1, r_2$ of the schema (LEA$^\nu_1$) which contain no $\iota$-description, by replacing, throughout $r_1, r_2$, their variables by numeroids. Thus, LEA$^\nu_1$ becomes

$$\mathbf{r}_1 = \mathbf{r}_2 \supset (F_{\mathbf{r}_1} \supset F_{\mathbf{r}_2})$$

which is Ari$^+$-provable when the antecedent $\mathbf{r}_1 = \mathbf{r}_2$ is disproved in Ari$^+$ (so that the propositional axiom Imp3 applies), and when the antecedent is proved in Ari, then $\mathbf{r}_1, \mathbf{r}_2$ and $F_{\mathbf{r}_1}, F_{\mathbf{r}_2}$ denote the same formula, say, E, in LAri so that the consequent $F_{\mathbf{r}_1} \supset F_{\mathbf{r}_2}$ becomes the Ari$^+$-provable implication

$$E \supset E,$$

from now on labeled as Imp0, whose easy proof (see Section 5.3.3., below) can be continued, in Ari$^+$, by using the axiom Imp1:

$$(E \supset E) \supset (\mathbf{r}_1 = \mathbf{r}_2 \supset (E \supset E)).$$

These Imp0 and Imp1 can be used as premises of the MP whose conclusion is the required instance of $(LEA_1^\nu)$.

5.3.1. The axiom $(LEA_1^\nu)$ occurs in this work only as the (closures of) the formulae 543., 556., 571. and 653. of the LAri$^\nu$-presentation of the deduction of Con$_{Ari}$ from T-Ax. The termoids $\mathbf{r}_1$ and $\mathbf{r}_2$ are, in these formulae, respectively:

543.  $((\mu)_2)_0$ and $2^3 \cdot 3^{((\mu)_1)_0} \cdot 5^{(\mu)_0}$ ;

556.  $((\mu)_1)_0$ and $(2^3 \cdot 3^{((\mu)_1)_0} \cdot 5^{(\mu)_0})_1$ ;

571.  $(2^3 \cdot 3^{((\mu)_1)_0} \cdot 5^{(\mu)_0})_2$ and $(\mu)_0$ ;

653.  $(\iota)_0$ and $\nu_\ell$ .

In each of the first three cases it is easy to observe that $\mathbf{r}_1$ and $\mathbf{r}_2$ contain no $\iota$-description and besides multiplication and exponentiation, only the following p.r. function symbols:

exp$_0$ , exp$_1$ , exp$_2$ .

Therefore, the substantiation given in Section 5.3, above, applies to these instances of $(LEA_1^\nu)$. In these cases, the numeroidal instances of the $\nu$-scheme contain no occurrence of $\iota$ (or any other definite description).

In the case of formula 653., i.e.

$$( \ )_0 = \nu_\ell \quad (F_{( \ )_0} \quad F_{\nu_\ell}),$$

no numeroidal instance arises because the termoids $r_1$, $r_2$ are constants. There is no need in replacing by $\nu$ and the formula 653. can be Ari$^+$-proved, by cases, as follows;

Case 1. In MAri: Assume  **a)**  = 0.

Then $( \ )_0 = 0 \quad \nu_\ell$ and the formulae $F_{( \ )_0}$, $F_{\nu_\ell}$ are, respectively, $F_0$ and $F_{\nu_\ell}$. As $cfor(0) = 1$, $fl(0) = \nu_\ell$, both of these formulae coincide with $\ell$ (as $F_{\nu_\ell}$ does by virtue of its definition). Therefore, the formula 653. has the consequent

$$\ell \quad \ell$$

which is an Imp0 implication from which formula 653. can be deduced, in Ari, using the axiom

$$(\ell \quad \ell) \quad (( \ )_0 = \nu_\ell \quad (\ell \quad \ell)),$$

i.e. Imp1, and MP.

Case 2. In MAri: Assume  **b)** ¬  = 0.

Then the equality  = 0 is false and, since the disjunction M( ) is true, its first disjunct holds; it entails $\ell( \ , \nu_\ell) = 0$ which entails (by virtue of the truth of $def_1^\ell$: (p.34, above) and of Con2 etc.) the truth of the equality $( \ )_0 = \nu_\ell$. In this case, $cfor(( \ )_0) = 0$ is entailed by $\ell( \ , \nu_\ell) = 0$ and by the assumption **b)** of the current case and $cfor(\nu_\ell) = 0$; therefore, the assumption entails $fl(( \ )_0) = ( \ )_0$ and $fl(( \ )_0) = \nu_\ell$; it follows (see pp.43-44, below) that the formula $F_{( \ )_0}$ coincides with $F_{\nu_\ell}$, i.e., with $\ell$, and so does also the last consequent of the formula 653..

It follows that the consequent, $F_{( \ )_0} \quad F_{\nu_\ell}$, of this formula coincides with $\ell \quad \ell$ also in the case of the assumption **b)**, and that the formula 653. has the same Ari-proof as in the case with assumption **a)**.

Thus, in MAri, also the $LEA_1^\nu$, i.e. the formula labeled 653., admits a substantiation by using the disjunction of the cases with assumptions **a)**, **b)**, above, which is an instance of TND (actually equivalent to $TND_\ell$).

**5**.3.2. The substantiation of $(LEA_2^\nu)$, used as line d19. of Appendix D2, can be given as follows. For any constant termoid, **r**, the **r**-instance of $(LEA_2^\nu)$ is

$$fl((r)_0) = (r)_0 \quad \pi_{F_{fl((r)_0)}} = \pi_{F_{(r)_0}};$$

it shall be found deducible, in the logic of $Ari^\nu$, from the disjunction

t.n.d.: $\quad fl((r)_0) = (r)_0 \lor \neg fl((r)_0) = (r)_0$

by cases, **ca**.1, **ca**.2, of this disjunction — as follows:

**ca**.1: $\quad fl((r)_0) = (r)_0$.

In this case the formulae

$$F_{fl((r)_0)} \text{ and } F_{(r)_0}$$

(considered as definite objects because the termoid **r**, in LAri, is constant) coincide, and therefore so do their En-N's, $\pi_{F_{fl((r)_0)}}$ and $\pi_{F_{(r)_0}}$.

1. $\quad \pi_{F_{fl((r)_0)}} = \pi_{F_{(r)_0}}$ \hfill (Ref),

2. $\quad 1. \quad 3.$ \hfill (Imp1),

3. $\quad fl((r)_0) = (r)_0 \quad \pi_{F_{fl((r)_0)}} = \pi_{F_{(r)_0}}$ \hfill (1,2;MP).

**ca**2: $\quad \neg fl((r)_0) = (r)_0$.

1. $\quad fl((r)_0) = (r)_0 \quad /$ \hfill (Def$_\neg$; **ca**2),

2. $\quad / \quad \pi_{F_{fl((r)_0)}} = \pi_{F_{(r)_0}}$ \hfill (Imp3),

3.-7.  $fl((r)_0) = (r)_0$    $\pi_{F_{fl((r)_0)}} = \pi_{F_{(r)_0}}$                                           (1,2;ch.in.).

**5**.3.3. The numeroidal instances of the scheme ($LEA^\nu_{MP}$) are the implications

$LEA^\nu_{MP}(r)$:   $(r)_0 = 3 \supset (F_r \supset (F_{(r)_1} \supset F_{(r)_2}))$

in which **r** may be any constant termoid in LAri. In the only matrix, formula 552., used in this work, the termoid r is

$$(2^3 \cdot 3^{((\mu)_1)_0} \cdot 5^{(\mu)_0})_0$$

and each numeroidal instance, **r**, of it is the termoid

$$(2^3 \cdot 3^{((\mu)_1)_0} \cdot 5^{(\mu)_0})_0$$

where μ denotes a numeroid; **r** is calculable.

In MAri:

The disjunction $cfor(r) = 0 \lor cfor(r) = 1$. Assume that $cfor(r) = 0$. Then $fl(r) = r$ and $F_r$ denotes the closed formula whose En-N equals **r**; the antecedent $(r)_0 = 3$, of the μ-instance of the implication $LEA^\nu_{MP}(r)$ entails that this formula, E, is an implication, $E_1 \supset E_2$, the sides, $E_1$, $E_2$, of which are closed formulae with En-N's $(r)_1$, $(r)_2$, respectively; the conditions $cfor((r)_1) = 0$, $cfor((r)_2) = 0$ are satisfied and entail, respectively, $fl((r)_1) = (r)_1$ and $fl((r)_2) = (r)_2$ so that the sides, $E_1$, $E_2$, of the implication $F_r$ — i.e. E — coincide with the closed formulae with En-N's $fl((r)_1)$, $fl((r)_2)$, i.e. with the formulae denoted by the notations $F_{(r)_1}$, $F_{(r)_2}$, respectively. It follows that $F_r$ — i.e. E — coincides with $F_{(r)_1} \supset F_{(r)_2}$ and $LEA^\nu_{MP}(r)$ has the consequent

          $F_r \supset F_r$                                                              (Imp0),

which is $Ari^+$-provable (or used as an additional (redundant) logical axiom) and which, with the aid of the axiom

$$(F_\mathbf{r} \quad F_\mathbf{r}) \quad ((\mathbf{r})_0 = 3 \quad (F_\mathbf{r} \quad F_\mathbf{r})) \qquad (\text{Imp1})$$

and MP, gives the Ari-proof of the implication $\text{LEA}^\nu_{\text{MP}}(\mathbf{r})$. (In this case, the antecedent $(\mathbf{r})_0 = 3$ is, in each (numeroidal) µ-instance, Ari-provable, but this observation is not used in the substantiation of the axiom scheme $\text{LEA}^\nu_{\text{MP}}(\mathbf{r})$.)

Now, assume that $\text{cfor}(\mathbf{r}) = 1$. Then, by definition of fl and $F_\mathbf{r}$ (in LAri), $\text{fl}(\mathbf{r}) = \nu_\ell$, the closed formula denoted by $F_\mathbf{r}$ is $F_{\nu_\ell}$, i.e. $\ell$ and $\text{LEA}^\nu_{\text{MP}}(\mathbf{r})$ denotes the LAri formula

$$(\mathbf{r})_0 = 3 \quad (\ell \quad (F_{(\mathbf{r})_1} \quad F_{(\mathbf{r})_2}))$$

the consequent of which is Ari-provable as an instance of (Imp3) so that the formula denoted by $\text{LEA}^\nu_{\text{MP}}(\mathbf{r})$ can be Ari-proved as follows:

1. $\ell \quad (F_{(\mathbf{r})_1} \quad F_{(\mathbf{r})_2})$ \hfill (Imp3),

2. 1. 3. \hfill (Imp1),

3. $(\mathbf{r})_0 = 3 \quad (\ell \quad (F_{(\mathbf{r})_1} \quad F_{(\mathbf{r})_2}))$ \hfill (1., 2.; MP).

Notice that the established coincidence of $F_{\nu_\ell}$ with $\ell$ entails that the implication

661. $\qquad F_{\nu_\ell} \quad \ell$

is Ari-provable as an instance of the (redundant) schema (Imp0) (proved in Appendix A).

**5**.3.4.

The $\text{Ari}^\nu$-presentation 1.-686a.(pp. , below) is given in the system $\text{Ari}^\nu$ with the open logic: so, it is assumed that its top-formula,(T-ax), deduced from the (redundant) 'axiom'

(T-Ax):             ($ax_{Ari}$( ) = 0 ⊃ F )

(to be considered as the hypothesis which subsequently will be proved with the aid of the 'SBA-Mp'

  $(\mu)_0$ ⊢ SBA : (T-Ax) ⊃ T-ax                           (SBA)

    T-ax: $ax_{Ari}((\mu)_0)$ = 0 ⊃ $F_{(\mu)_0}$

Further, it will be supposed that its formulae 7. and 630. are used as the Gen-premises, for obtaining their -closures, 7a. and 630a., which will be used as the 'mp-ind-premises', that is, the formulae 7. and 630. will be considered as followed by the formulae

7a.    $\mu(Ant_1^{\mathscr{b}_1}(\mu)$  ⊃ $F_{(\mu)_0})$                                  (7.;Gen),

630a.  $\mu(Ant_2^{\mathscr{b}_1}(\mu)$ ⊃ $(F_{((\mu)_1)_0})$ ⊃ $(F_{((\mu)_2)_0})$ ⊃ $F_{(\mu)_0})))$    (630.;Gen)

which can be used as antecedents of the mp-induction formula with $F_{(\mu)_0}$ used as $C(\mu)$ of **mp-ind** (cf. Appendix C). As also the proof of **mp-ind** is presented in the open logic of Ari, its formula $C(\mu)$ must be replaced by $F_{(\mu)_0}$ in order to obtain from the $Ari^o$-proof the $Ari^\nu$-proof of the formula

631.   7a. ⊃ (630a. ⊃ 633.)                                             (mp-ind).

The top-formula of the $Ari^\nu$-deduction 1.-686a. and in the Ari-proof of the scheme **mp-ind** shall be, besides (T-Ax), either non-logical axioms in the (open) system Ari, or, else, instances of logical axioms, so that, when the Kleene notations occur in them these formulae will denote (for any numeroids substituted for their indices) the logical axioms by the same schemata, or, finally, the instances of the $\nu$-schemata considered in 5.3.

5.4.    Although the formula M- is treated as an axiom in 1.-686a. — formula 634.— it can actually be proved as an Ari-theorem. Recall that this formula is

(M- ):  $\mathscr{b}$( ,$\nu_\ell$) = 0 & ¬ ( < & $\mathscr{b}$( ,$\nu_\ell$) = 0).

        ∨ . = 0 & ¬ $\mathscr{b}$( ,$\nu_\ell$) = 0,

where the constant termoid   is formally introduced by the definite description **def** ; ($\mathscr{b}(\ ,\nu_\ell) = 0$ & ¬ ( < & $\mathscr{b}(\ ,\nu_\ell) = 0$).

∨. = 0 & ¬ $\mathscr{b}(\ ,\nu_\ell) = 0$),

in LAri. (M- ) shall follow from the Ari-provable equivalence

= ∼ M-

which can be proved with the aid of Rosser's axiom 11 (in [ROSSER, 1953]), the SBA ( ) and an application (for the use of this axiom) of the

TND$_\ell$ $\mathscr{b}(\ ,\nu_\ell) = 0\ \vee\ \neg\ \mathscr{b}(\ ,\nu_\ell) = 0$.

which is postulated in Ari as a non-logical axiom

5.4.1. In more detail, the termoid denotes the integer (if any) whose existence can be inferred, in the intuitionistic Ari, from the disjunction

TND$_\ell$ : $\mathscr{b}(\ ,\nu_\ell) = 0\ \vee\ \neg\ \mathscr{b}(\ ,\nu_\ell) = 0$.

More specifically, this TND$_\ell$ entails, intuitionistically, the formula

(A1.) ( $\mathscr{b}(\ ,\nu_\ell) = 0$ $\mathscr{b}(\ ,\nu_\ell) = 0$)

using §27., *34; §35, *90; §27., *59a; §32., *70 in [KLEENE, 1952]. Using TND$_\ell$ again, formula A1. entails, together with the least number principle, *149°, the formula

(A2.) (( $\mathscr{b}(\ ,\nu_\ell) = 0$ $\mathscr{b}(\ ,\nu_\ell) = 0$)

& ¬ ( < & ( $\mathscr{b}(\ ,\nu_\ell) = 0$ $\mathscr{b}(\ ,\nu_\ell) = 0$)))

which expresses the existence of a least integer, , which satisfies the scope of the external in A2. The uniqueness of this least can be expressed by the formula

(A3.) ! ($\mathscr{b}(\ ,\nu_\ell) = 0$ & ¬ ( < & $\mathscr{b}(\ ,\nu_\ell) = 0$). ∨

. = 0 & ¬ $\mathscr{b}(\ ,\nu_\ell) = 0$)

which is deducible from TND$_f$ intuitionistically, only the double negation of A3. can thus be found to be Ari-provable.

The last formula and Rosser's axiom schema 11 for definite descriptions in [ROSSER, 1953] together with the definition of    give, in Ari, first

(A4.)    ( = ~ . ( $\mathscr{b}$ ( ,v$_f$) = 0 & ¬   ( <   & $\mathscr{b}$ ( ,v$_f$) = 0 ).
         ∨. = 0 & ¬   $\mathscr{b}$ ( ,v$_f$) = 0 ))

and then, with the aid of SBA   (       ) and MP, the equivalence

         =   ~ : ( $\mathscr{b}$ ( ,v$_f$) = 0 & ¬   ( <   & $\mathscr{b}$ ( ,v$_f$) = 0 ).
   ∨. = 0 & ¬   $\mathscr{b}$ ( ,v$_f$) = 0 ),

the left side of which is a logical axiom, Ref. (Section 3., p. 30). Hence, the right side of this equivalence, *viz.*,

(M- )    $\mathscr{b}$ ( ,v$_f$) = 0 & ¬   ( <   & $\mathscr{b}$ ( ,v$_f$) = 0 ).
         ∨. = 0 & ¬   $\mathscr{b}$ ( ,v$_f$) = 0

is deducible from TND$_f$ in Ari; it follows that

    TND$_f$     (M-   )

is provable in Ari, and so is (see [KLEENE, 1952]: *27, *58a and *60h.,i.)

    ¬ ¬ TND$_f$     ¬ ¬ (M-   );

since ¬ ¬ TND$_f$ is provable intuitionistically (cf. §27., *51a in [KLEENE, 1952]), the formula ¬ ¬ M-   is also provable in the intuitionistic Ari.

5.5.  We now begin studying the important theorem M-   in some detail. To begin with, observe that

I.        $\mathscr{b}$ ( ,v$_f$) = 0 &     = 0

leads to a contradiction entailed by $\mathscr{b}(0, v_\ell) = 0$ (which entails $(0)_0 = v_\ell$ and, hence, $0 = 1$). Furthermore;

II. $\qquad \mathscr{b}(\ ,v_\ell) = 0 \ \& \ \neg \quad \mathscr{b}(\ ,v_\ell) = 0$

also leads to a contradiction because $\mathscr{b}(\ ,v_\ell) = 0$ and SBA $(\quad)$ give $\mathscr{b}(\ ,v_\ell) = 0$ contradicting $\neg \ \mathscr{b}(\ ,v_\ell) = 0$.

The provability of the negations of I. and II. entails the provability of the implications

$$\mathscr{b}(\ ,v_\ell) = 0 \quad \neg \quad = 0,$$

$$\mathscr{b}(\ ,v_\ell) = 0 \quad \neg\neg \quad \mathscr{b}(\ ,v_\ell) = 0,$$

$$\mathscr{b}(\ ,v_\ell) = 0 \ . \ \neg \quad = 0 \ \vee \ \neg\neg \quad \mathscr{b}(\ ,v_\ell) = 0$$

and (using [KLEENE, 1952] §27, *34, *62a) also of the implication

III. $\qquad \mathscr{b}(\ ,v_\ell) = 0 \quad \neg(\quad = 0 \ \& \ \neg \quad \mathscr{b}(\ ,v_\ell) = 0)$

the consequent of which yields, together with M-$\quad$, both

$$\mathscr{b}(\ ,v_\ell) = 0 \ \& \ \neg \quad (\quad < \quad \& \ \mathscr{b}(\ ,v_\ell) = 0)$$

and

$$\neg \quad (\quad < \quad \& \ \mathscr{b}(\ ,v_\ell) = 0);$$

therefore, the implication

IV. $\qquad \mathscr{b}(\ ,v_\ell) = 0 \quad \neg \quad (\quad < \quad \& \ \mathscr{b}(\ ,v_\ell) = 0)$

is provable in Ari.

By §35., *86 in [KLEENE, 1952], $\neg \quad$ can be replaced in IV. by $\quad \neg$ and by §27., *58b and §32., *69, of the same reference, the consequent in IV. can be replaced by

$$(\quad < \quad \quad \neg \ \mathscr{b}(\ ,v_\ell) = 0);$$

thus, the implication

V.  $\quad b(\ ,v_\ell) = 0 \quad (\ <\quad \neg b(\ ,v_\ell) = 0)$

is provable (in the classical Ari). Now, observe that the definition of $b$ (Section 4.3, p.32) contains the item, $(\text{def}_3^b)$, which entails $(\neg b(\ ,v_\ell) = 0 \quad b(\ ,v_\ell) = 1)$; therefore, the provability of V. entails the provability of the implication

Va. $\quad b(\ ,v_\ell) = 0 \quad (\ <\quad b(\ ,v_\ell) = 1).$

On the other hand, we have that either conjunct in

$$= 0 \ \& \ \neg \ b(\ ,v_\ell) = 0$$

entails the negation of $b(\ ,v_\ell) = 0$ as well as the negation of the conjunction

$$b(\ ,v_\ell) = 0 \ \& \ \neg\ (\ <\ \& \ b(\ ,v_\ell) = 0),$$

while §27., *61a in [KLEENE, 1952] together with M- entail the implication

VI. $\quad \neg (b(\ ,v_\ell) = 0 \ \& \ \neg\ (\ <\ \& \ b(\ ,v_\ell) = 0))$

$$= 0 \ \& \ \neg\ b(\ ,v_\ell) = 0$$

and therefore both of the implications

$$= 0 \quad\quad = 0 \ \& \ \neg\ b(\ ,v_\ell) = 0$$

and

$$\neg\ b(\ ,v_\ell) = 0 \quad\quad = 0 \ \& \ \neg\ b(\ ,v_\ell) = 0$$

are provable in Ari. The consequent of these implications entails each of $= 0$ and $\neg\ b(\ ,v_\ell) = 0$. Therefore, each of the implications

VII. $\quad = 0 \quad \neg\ b(\ ,v_\ell) = 0$

and

VIII.  ¬ $\mathscr{b}($ , $\nu_\ell) = 0$    = 0

is provable in the classical Ari. The antecedent of the last implication is a standard formula expressing, in LAri, the consistency of Ari; accordingly, it shall be denoted by $\mathrm{Con}_{\mathrm{Ari}}$. Hence, the last implications can be, respectively, rewritten as

VII.a.    = 0    $\mathrm{Con}_{\mathrm{Ari}}$

VIII.a.  $\mathrm{Con}_{\mathrm{Ari}}$    = 0

by means of which it becomes obvious that    = 0 and $\mathrm{Con}_{\mathrm{Ari}}$ are, in Ari, equivalent formulae. Remark that the equality    = 0 is quantifier-free.

5.6. Lastly, we shall make some observations about the relation between $\mathrm{Con}_{\mathrm{Ari}}$ and the *numeroidal* instance of the reflection principle. Remark first that    = 0 entails ¬   ( <   & $\mathscr{b}($ , $\nu_\ell) = 0$ ) (and that the implication IV. can be inverted only if ¬   = 0, i.e. only if $\mathrm{Con}_{\mathrm{Ari}}$ is false). Observing next that

$\mathscr{b}($ , $\nu_\ell) = 0$ & ¬   ( <   & $\mathscr{b}($ , $\nu_\ell) = 0$ )    $\mathscr{b}($ , $\nu_\ell) = 0$

is a logical axiom, and that its contraposition gives the implication

IX.    ¬ $\mathscr{b}($ , $\nu_\ell) = 0$
    ¬ ( $\mathscr{b}($ , $\nu_\ell) = 0$ & ¬   ( <   & $\mathscr{b}($ , $\nu_\ell) = 0$ ))

which with chain inference with VI., above, gives

X.    ¬ $\mathscr{b}($ , $\nu_\ell) = 0$    = 0 & ¬    $\mathscr{b}($ , $\nu_\ell) = 0$ —

and the definition $\mathrm{Con}_{\mathrm{Ari}}$ entails that the Ari-provability of this implication entails that the implication

X.a.    ¬ $\mathscr{b}($ , $\nu_\ell) = 0$    $\mathrm{Con}_{\mathrm{Ari}}$

is Ari-provable.

Remark that the provability of the last implication can actually be inferred immediately because M- entails the disjunction

$$\mathcal{E}(\ ,\nu_{\ell}) = 0 \quad \neg\ \mathcal{E}(\ ,\nu_{\ell}) = 0.$$

using §27., *61a in [KLEENE, 1952], together with the definition **def** gives Xa.

It is easy to check that the only use of non-intuitionistic axiom in the proof of X.a. is that of the $TND_{\ell}$ in the proof of M- ; therefore, $\neg\neg$X.a. is intuitionistically provable (recall that by §27., *49b and *60g-i in [KLEENE, 1952], a formula of the shape $\neg A \quad \neg B$ is equivalent to its double negation). X.a. can also be written as

X.b.  $(\mathcal{E}(\ ,\nu_{\ell}) = 0 \quad \ell) \quad \text{Con}_{Ari}$ .

Therefore, it suffices to prove the antecedent,

-p.t.i.:  $\mathcal{E}(\ ,\nu_{\ell}) = 0 \quad \ell$ ,

of X.b. in order to obtain, in Ari, also the resulting proof of the formula $\text{Con}_{Ari}$ and — if the second Gödel theorem has a proof applicable to this case — also of a contradiction formula. (In the Introduction, Section 0.1., the formula -p.t.i. was labeled L- in order to stress its connection with L- , i.e. $\text{Con}_{Ari}$ ).

On the other hand, the obtaining of a *numeroidal instance* of formula -p.t.i requires that the equality

$$= \upsilon$$

is actually provable *in* Ari. Of course, if $= 0$ that may not be so difficult — but since the case where is *positive* must also be considered the proof, in Ari, of this equality is not a trivial task. This task will be considered in Part II of this work — meanwhile, immediatly below, we give an argument, *in* MAri, which does not assume the accomplishment of this task. However, notice that numeroidal instances of certain formulae *are* needed in Appendices D1 and D2 — nevertheless, this need may be circumvented and this point will be taken up in Part II (at which point a complete revision of the present proof presentation will be considered).

.6. Below, follows the *preliminary* proof-presentation 1.-686a. of the formula Con$_{Ari}$. This preliminary proof-presentation is intended to be read in conjunction with the comments given in Section 5, above. Thus, *given* that (i) the additional lines of Section 5.3.4., above, are edited into the preliminary proof-presentation, (ii) that proofs of the induction hypotheses are implanted at formulae 7. and 630., respectively (with due attentio to editing matters regarding the choice of variable), (iii) implanting Ari-proofs of all instances of (t-ax); (iv) all Kleene formulae are translated into corresponding formulae in LAri and, finally, (v) all formulae must be closed, *then* the following preliminary proof-presentation is also a proof-presentation in Ari.

Thus, the reader may *first* consider the presentation as a *deduction* carried out in the extension Ari$^\nu$; and *provided* all hypotheses have been Ari$^\nu$-proved, *secondly*, as a Ari$^\nu$-proof — and only *thereafter* as a *proof* in Ari$^+$ or the ground theory Ari obtained by systematically replacing all Kleene formulae F$_\nu$ in LAri$^\nu$ by the corresponding formulae in LAri.

Thus, we are using a *single* text, 1-686a., to play the role of *several distinct proofs* of which the proof in Ari is obtainable by an intuitive translation from LAri$^\nu$ into LAri (a translation which, for the non-finitist, of course, is always possible to write down 'in principle'). The final editing of this preliminary proof-presentation will be given in Part II. In the meantime, the reader is invited to write out the proof with all details and, in addition, consider any possible simplification(s). (Some substantial simplifications have already been found by the present authors which will be incorporated in this work with the publication of its Part II).

6.1. Some of the proof comments in 1.-686a., below, are reduced to the expression [...] — which means that the comment exceeded the typographical space available on the page. All such proof comments refer to elementary arithmetical axioms so far not taken into attention and their full texts are relegated to Appendix **B** while the remaining labels not accounted for are Ga$_0$, Ga$_1$ and Ga$_2$ which refer to Gauss' laws for prime exponents.

In order to define the later an additional function symbol, so far not introduced (and which also occurs, below; *viz.*, at line 421)., $\ell$, denotes the function symbol for which $\ell(\mathbf{n})$, $\mathbf{n} > 1$, shall denote the largest (if any) integer f such that $\mathbf{n}$ is divisible by $p_f$) — where the latter notation denotes the $f^{th}$ prime. For $\mathbf{n} = 0, 1$, $\ell(\mathbf{n})$ shall equal zero.

For each positive $\mathbf{n}$, the Gaussian Theorem claims the uniqueness of the decomposition of $\mathbf{n}$ in the prime factors:

$$\mathbf{n} = \prod_{g=0}^{\ell(\mathbf{n})} p_g(\mathbf{n})g.$$

On the other hand, for each unary function  and integers h, k such that h  k , the equality

(Ga$_h$)    $(\bigwedge_{g=0}^{k} p_g^{(g)})_h$ =  (h)

holds. For h = 0, 1, 2, …, instances of these equalities are labeled Ga$_0$, Ga$_1$, Ga$_2$, ….

Further, a number of derived rules of inference labeled ch.in., ch.in.$_2$ , ch-fla, Int-Ant, Contrap. and MTP$_2$ occurs and they are Ari-proved in Appendix **A**.
    Formula 7., below, is proved in Appendix D2 and the mp-induction formula, i.e. formula 630., below, is proved in Appendix C. In Part II of this work, the full proof will be completed and the ensuing paradox(es) resolved.
.

1.　　　$ax_{Ari}((\mu)_0) = 0$　　$F_{(\mu)_0}$　　　　　　　　　　　　　　　　　　　　(T-Ax),

2.　　　$\mu = 2^{(\mu)_0}$ & $ax_{Ari}((\mu)_0) = 0$　　$ax_{Ari}((\mu)_0) = 0$　　　　　　(Con 2),

　　2a.　$Ant_1^{\ell_1}(\mu)$　　$ax_{Ari}((\mu)_0) = 0$　　　　　　　　　　　　　　　　[$\mathbf{def}_{Ant_1^{\ell_1}}$; 2],

3.-7.　$Ant_1^{\ell_1}(\mu)$　　$F_{(\mu)_0}$　　　　　　　　　　　　　　　　　　　　(2a., 1; ch.in.),

8.　　　$\mu = mp\,((\mu)_1, (\mu)_2)$ & $\neg\,\mu = 0$ & $\ell_1((\mu)_1) = 0$ & $\ell_1((\mu)_2) = 0$

　　　　$\ell_1((\mu)_2) = 0$　　　　　　　　　　　　　　　　　　　　　　　　　　(Con 2),

　　8a.　$Ant_2^{\ell_1}(\mu)$　　$\ell_1((\mu)_2) = 0$　　　　　　　　　　　　　　　　[$\mathbf{def}_{Ant_2^{\ell_1}}$; 8],

9.　　　$\mu = mp\,((\mu)_1, (\mu)_2)$ & $\neg\,\mu = 0$ & $\ell_1((\mu)_1) = 0$ & $\ell_1((\mu)_2) = 0$

　　　　$\mu = mp\,((\mu)_1, (\mu)_2)$ & $\neg\,\mu = 0$ & $\ell_1((\mu)_1) = 0$　　　(Con 1),

　　9a.　$Ant_2^{\ell_1}(\mu)$　　$\mu = mp\,((\mu)_1, (\mu)_2)$ & $\neg\,\mu = 0$ & $\ell_1((\mu)_1) = 0$　[$\mathbf{def}_{Ant_2^{\ell_1}}$; 9],

10.　　$\mu = mp\,((\mu)_1, (\mu)_2)$ & $\neg\,\mu = 0$ & $\ell_1((\mu)_1) = 0$　　$\ell_1((\mu)_1) = 0$　(Con 2),

11.-15.　$Ant_2^{\ell_1}(\mu)$　　$\ell_1((\mu)_1) = 0$　　　　　　　　　　　　　　　(9a., 10; ch.in.),

16.　　$\mu = mp\,((\mu)_1, (\mu)_2)$ & $\neg\,\mu = 0$ & $\ell_1((\mu)_1) = 0$

　　　　$\mu = mp\,((\mu)_1, (\mu)_2)$ & $\neg\,\mu = 0$　　　　　　　　　　　　　(Con 1),

17.-21.　$Ant_2^{\ell_1}(\mu)$　　$\mu = mp\,((\mu)_1, (\mu)_2)$ & $\neg\,\mu = 0$　　　　(9a., 16.; ch.in.),

22.　　$\mu = mp\,((\mu)_1, (\mu)_2)$ & $\neg\,\mu = 0$　　$\neg\,\mu = 0$　　　　　　　(Con 2),

23.-27.　$Ant_2^{\ell_1}(\mu)$　　$\neg\,\mu = 0$　　　　　　　　　　　　　　　　　　(21., 22.; ch.in.),

28.　　$\mu = mp\,((\mu)_1, (\mu)_2)$ & $\neg\,\mu = 0$　　$\mu = mp\,((\mu)_1, (\mu)_2)$　(Con 1),

29.-33.　$Ant_2^{\ell_1}(\mu)$　　$\mu = mp\,((\mu)_1, (\mu)_2)$　　　　　　　　　　　(21.,28.; ch.in.),

| | | |
|---|---|---|
| 34. | mp $((\mu)_1, (\mu)_2)$ = Mp $((\mu)_1, (\mu)_2) \cdot$ $((\mu)_1, (\mu)_2) \cdot$ sg $((\mu)_1)$ | [**def**$_{mp}$], |
| 35. | 34. 36. | (LEA$_{\bar{2}}$), |
| 36. | $\mu$ = mp $((\mu)_1, (\mu)_2)$ $\mu$ = Mp $((\mu)_1, (\mu)_2)\cdot$ $((\mu)_1, (\mu)_2) \cdot$ sg $((\mu)_1)$ | (34.,35.; MP), |
| 37.-41. | Ant$_2^{\ell_1}(\mu)$ $\mu$ = Mp $((\mu)_1, (\mu)_2) \cdot$ $((\mu)_1, (\mu)_2) \cdot$ sg $((\mu)_1)$ | (33., 36.; ch.in.), |
| 42. | sg $((\mu)_1) = 0$ sg$((\mu)_1) = 1$ | […], |
| 43. | sg $((\mu)_1) = 0$ | |
| | Mp $((\mu)_1, (\mu)_2) \cdot$ $((\mu)_1, (\mu)_2) \cdot$ sg $((\mu)_1)$ = Mp $((\mu)_1, (\mu)_2) \cdot$ $((\mu)_1, (\mu)_2) \cdot 0$ | (LEA-rp) |
| 44. | Mp $((\mu)_1, (\mu)_2) \cdot$ $((\mu)_1, (\mu)_2) \cdot 0 = 0$ | […], |
| 45. | 44. 46. | (LEA$_{\bar{2}}$), |
| 46. | Mp $((\mu)_1, (\mu)_2) \cdot$ $((\mu)_1, (\mu)_2) \cdot$ sg $((\mu)_1)$ = Mp $((\mu)_1, (\mu)_2) \cdot$ $((\mu)_1, (\mu)_2) \cdot 0$ | |
| | Mp $((\mu)_1, (\mu)_2) \cdot$ $((\mu)_1, (\mu)_2) \cdot$ sg $((\mu)_1) = 0$ | (44., 45.; MP), |
| 47.-51. | sg$((\mu)_1) = 0$ Mp $((\mu)_1, (\mu)_2) \cdot$ $((\mu)_1, (\mu)_2) \cdot$ sg $((\mu)_1) = 0$ | (43., 46.; ch.in.), |
| 52. | Mp $((\mu)_1, (\mu)_2) \cdot$ $((\mu)_1, (\mu)_2) \cdot$ sg $((\mu)_1) = 0$ | |
| | ($\mu$ = Mp $((\mu)_1, (\mu)_2) \cdot$ $((\mu)_1, (\mu)_2) \cdot$ sg $((\mu)_1)$ $\mu = 0$) | (LEA$_{\bar{2}}$), |
| 53.–57. | sg$((\mu)_1) = 0$ ($\mu$ = Mp $((\mu)_1, (\mu)_2) \cdot$ $((\mu)_1, (\mu)_2) \cdot$ sg $((\mu)_1)$ $\mu = 0$) | (51., 52.; ch.in.), |
| 58.-65. | ($\mu$ = Mp $((\mu)_1, (\mu)_2) \cdot$ $((\mu)_1, (\mu)_2) \cdot$ sg $((\mu)_1)$ (sg $((\mu)_1) = 0$ $\mu = 0$) | (57.; Int.-Ant.), |
| 66.-70. | Ant$_2^{\ell_1}(\mu)$ (sg $((\mu)_1) = 0$ $\mu = 0$) | (41., 65.; ch.in.), |
| 71.-85. | (sg $((\mu)_1) = 0$ $\mu = 0$) ($\neg \mu = 0$ $\neg$ sg $((\mu)_1) = 0$) | (contrap.), |
| 86.-90. | Ant$_2^{\ell_1}(\mu)$ ($\neg \mu = 0$ $\neg$ sg $((\mu)_1) = 0$) | (70., 85.; ch.in.), |
| 91. | 90. (27. 93.) | (Imp 2), |
| 92. | 27. 93. | (90., 91.; MP), |

| | | | |
|---|---|---|---|
| 93. | $\text{Ant}_2^{\ell_1}(\mu)$ | $\neg\,\text{sg}\,((\mu)_1) = 0$ | (27., 92.; MP), |
| 94.-113. | $\neg\,\text{sg}\,((\mu)_1) = 0$ | $\text{sg}\,((\mu)_1) = 1$ | (42.; Mtp$_2$), |
| 114.-118. | $\text{Ant}_2^{\ell_1}(\mu)$ | $\text{sg}\,((\mu)_1) = 1$ | (93., 113.; ch.in.), |

119.      $\text{sg}\,((\mu)_1) = 1$

       $\text{Mp}\,((\mu)_1, (\mu)_2) \cdot\ ((\mu)_1, (\mu)_2) \cdot \text{sg}\,((\mu)_1) =$

       $\text{Mp}\,((\mu)_1, (\mu)_2) \cdot\ ((\mu)_1, \cdot 1$                              (LEA$_{\overline{2}}$),

120.      $\text{Mp}\,((\mu)_1, (\mu)_2) \cdot\ ((\mu)_1, (\mu)_2)\ 1 = \text{Mp}\,((\mu)_1, (\mu)_2) \cdot\ ((\mu)_1, (\mu)_2)$     […],

121.            120.     122.                                       (LEA$_{\overline{2}}$),

122.      $\text{Mp}\,((\mu)_1, (\mu)_2) \cdot\ ((\mu)_1, (\mu)_2) \cdot \text{sg}\,((\mu)_1) =$

       $\text{Mp}\,((\mu)_1, (\mu)_2) \cdot\ ((\mu)_1, (\mu)_2)\ 1$

       $\text{Mp}\,((\mu)_1, (\mu)_2) \cdot\ ((\mu)_1, (\mu)_2) \cdot \text{sg}\,((\mu)_1) =$

       $\text{Mp}\,((\mu)_1, (\mu)_2) \cdot\ ((\mu)_1, (\mu)_2)$                    (120., 121.; MP),

123.-127.    $\text{sg}\,((\mu)_1) = 1$

       $\text{Mp}\,((\mu)_1, (\mu)_2) \cdot\ ((\mu)_1, (\mu)_2) \cdot \text{sg}\,((\mu)_1) =$

       $\text{Mp}\,((\mu)_1, (\mu)_2) \cdot\ ((\mu)_1, (\mu)_2)$                    (119., 122.; ch.in.),

128.-132.    $\text{Ant}_2^{\ell_1}(\mu)$

       $\text{Mp}\,((\mu)_1, (\mu)_2) \cdot\ ((\mu)_1, (\mu)_2) \cdot \text{sg}\,((\mu)_1) =$

       $\text{Mp}\,((\mu)_1, (\mu)_2) \cdot\ ((\mu)_1, (\mu)_2)$                    (118., 127.; ch.in.),

133.     $\text{Mp}\,((\mu)_1, (\mu)_2) \cdot\ ((\mu)_1, (\mu)_2) \cdot \text{sg}\,((\mu)_1) =$

       $\text{Mp}\,((\mu)_1, (\mu)_2) \cdot\ ((\mu)_1, (\mu)_2)$

       $(\mu = \text{Mp}\,((\mu)_1, (\mu)_2) \cdot\ ((\mu)_1, (\mu)_2) \cdot \text{sg}\,((\mu)_1)$

$$\mu = \text{Mp}((\mu)_1, (\mu)_2) \cdot ((\mu)_1, (\mu)_2))) \qquad (\text{LEA}_{\overline{2}}^{=}),$$

134.-138. $\quad \text{Ant}_2^{\ell_1}(\mu)$

$$(\mu = \text{Mp}((\mu)_1, (\mu)_2) \cdot ((\mu)_1, (\mu)_2) \cdot \text{sg}((\mu)_1)$$

$$\mu = \text{Mp}((\mu)_1, (\mu)_2) \cdot ((\mu)_1, (\mu)_2))) \qquad (132., 133.; \text{ch.in.}),$$

| | | | |
|---|---|---|---|
| 139. | 138. (41. 141.) | | (Imp 2), |
| 140. | 41. 141. | | (138., 139.; MP), |
| 141. | $\text{Ant}_2^{\ell_1}(\mu) \quad \mu = \text{Mp}((\mu)_1, (\mu)_2) \cdot ((\mu)_1, (\mu)_2)$ | | (41., 140.; MP), |
| 142. | $((\mu)_1, (\mu)_2) = \overline{\text{sg}}(\text{msd}((\mu)_2)_{0,0}, 3)) \cdot \overline{\text{sg}}(\text{msd}((\mu)_2)_{0,1}, ((\mu)_1)_0))$ | | [**def**], |
| 143. | $\overline{\text{sg}}(\text{msd}((\mu)_2)_{0,1}, ((\mu)_1)_0)) = 0 \quad \overline{\text{sg}}(\text{msd}((\mu)_2)_{0,1}, ((\mu)_1)_0)) = 1$ | | […], |
| 144. | $\overline{\text{sg}}(\text{msd}((\mu)_2)_{0,1}, ((\mu)_1)_0)) = 0$ | | |

$$\overline{\text{sg}}(\text{msd}((\mu)_2)_{0,0}, 3)) \cdot \overline{\text{sg}}(\text{msd}((\mu)_2)_{0,1}, ((\mu)_1)_0)) =$$

$$\overline{\text{sg}}(\text{msd}((\mu)_2)_{0,0}, 3)) \cdot 0 \qquad (\text{LEA-rp}),$$

145. $\quad \overline{\text{sg}}(\text{msd}((\mu)_2)_{0,0}, 3)) \cdot 0 = 0 \qquad \qquad [\ldots],$

146. $\qquad\qquad 145. \quad 147. \qquad\qquad\qquad (\text{LEA}_{\overline{2}}^{=}),$

147. $\quad \overline{\text{sg}}(\text{msd}((\mu)_2)_{0,0}, 3)) \cdot \overline{\text{sg}}(\text{msd}((\mu)_2)_{0,1}, ((\mu)_1)_0)) =$

$$\overline{\text{sg}}(\text{msd}((\mu)_2)_{0,0}, 3)) \cdot 0$$

$$\overline{\text{sg}}(\text{msd}((\mu)_2)_{0,0}, 3)) \cdot \overline{\text{sg}}(\text{msd}((\mu)_2)_{0,1}, ((\mu)_1)_0)) = 0 \qquad (145., 147.; \text{MP}),$$

148.-152. $\quad \overline{\text{sg}}(\text{msd}((\mu)_2)_{0,1}, ((\mu)_1)_0)) = 0$

$$\overline{\text{sg}}(\text{msd}((\mu)_2)_{0,0}, 3)) \cdot \overline{\text{sg}}(\text{msd}((\mu)_2)_{0,1}, ((\mu)_1)_0)) = 0 \qquad (144., 147.; \text{ch.in.}),$$

153. $\quad \overline{\text{sg}}(\text{msd}((\mu)_2)_{0,0}, 3)) \cdot \overline{\text{sg}}(\text{msd}((\mu)_2)_{0,1}, ((\mu)_1)_0)) = 0$

$$(142. \quad ((\mu)_1, (\mu)_2) = 0) \qquad (\text{LEA}_{\overline{2}}^{=}),$$

| | | |
|---|---|---|
| 154.-158. | $\overline{sg} \, (msd \, ((\mu)_2)_{0,1}, ((\mu)_1)_0)) = 0 \quad (142. \quad ((\mu)_1, (\mu)_2) = 0)$ | (152., 153.; ch.in.), |
| 159.-168. | $\overline{sg} \, (msd \, ((\mu)_2)_{0,1}, ((\mu)_1)_0)) = 0 \quad ((\mu)_1, (\mu)_2) = 0$ | (158., 142.; ch.in.$_2$), |
| 169. | $((\mu)_1, (\mu)_2) = 0 \quad Mp \, ((\mu)_1, (\mu)_2) \cdot ((\mu)_1, (\mu)_2) = Mp \, ((\mu)_1, (\mu)_2) \cdot 0$ | (LEA-rp), |
| 170. | $Mp \, ((\mu)_1, (\mu)_2) \cdot 0 = 0$ | […], |
| 171. | 170. 172. | (LEA$_2^=$), |
| 172. | $Mp \, ((\mu)_1, (\mu)_2) \cdot ((\mu)_1, (\mu)_2) = Mp \, ((\mu)_1, (\mu)_2) \cdot 0$ | |
| | $Mp \, ((\mu)_1, (\mu)_2) \cdot ((\mu)_1, (\mu)_2) = 0$ | (170., 171.; MP), |
| 173.-177. | $((\mu)_1, (\mu)_2) = 0 \quad Mp \, ((\mu)_1, (\mu)_2) \cdot ((\mu)_1, (\mu)_2) = 0$ | (169., 172.: ch.in.), |
| 178. | $Mp \, ((\mu)_1, (\mu)_2) \cdot ((\mu)_1, (\mu)_2) = 0$ | |
| | $(\mu = Mp \, ((\mu)_1, (\mu)_2) \cdot ((\mu)_1, (\mu)_2) \quad \mu = 0)$ | (LEA$_2^=$), |
| 179.-183. | $((\mu)_1, (\mu)_2) = 0$ | |
| | $(\mu = Mp \, ((\mu)_1, (\mu)_2) \cdot ((\mu)_1, (\mu)_2) \quad \mu = 0)$ | (177., 178; ch.in.), |
| 184.-191. | $\mu = Mp \, ((\mu)_1, (\mu)_2) \cdot ((\mu)_1, (\mu)_2) \quad (((\mu)_1, (\mu)_2) = 0 \quad \mu = 0)$ | (183.; Int-Ant.), |
| 192.-196. | $Ant_2^{\ell_1}(\mu) \quad (((\mu)_1, (\mu)_2) = 0 \quad \mu = 0)$ | (141., 191.; ch.in.), |
| 197.-211. | $(((\mu)_1, (\mu)_2) = 0 \quad \mu = 0) \quad (\neg \mu = 0 \quad \neg(((\mu)_1, (\mu)_2) = 0)$ | (contrap.), |
| 212.-216. | $Ant_2^{\ell_1}(\mu) \quad (\neg \mu = 0 \quad \neg(((\mu)_1, (\mu)_2) = 0)$ | (196., 211.; ch.in.), |
| 217. | 216. (27. 219.) | (Imp 2), |
| 218. | 27. 219. | (216., 217.; MP), |
| 219. | $Ant_2^{\ell_1}(\mu) \quad \neg \, ((\mu)_1, (\mu)_2) = 0$ | (27., 218.; MP), |
| 220.-231. | $\neg \, ((\mu)_1, (\mu)_2) = 0 \quad \neg \overline{sg} \, (msd \, ((\mu)_2)_{0,1}, ((\mu)_1)_0)) = 0$ | (168.; Contrap), |

| | | | |
|---|---|---|---|
| 232.-236. | $\mathrm{Ant}_2^{\ell_1}(\mu)$ | $\neg\ \overline{sg}\ (\mathrm{msd}\ ((\mu)_2)_{0,1}\ ,\ ((\mu)_1)_0)) = 0$ | (219., 231.; ch.in.), |
| 237.-256. | $\neg\ \overline{sg}\ (\mathrm{msd}\ ((\mu)_2)_{0,1}\ ,\ ((\mu)_1)_0)) = 0$ | $\overline{sg}\ (\mathrm{msd}\ ((\mu)_2)_{0,1}\ ,\ ((\mu)_1)_0)) = 1$ | (143.; Mtp$_2$), |
| 257.-261. | $\mathrm{Ant}_2^{\ell_1}(\mu)$ | $\overline{sg}\ (\mathrm{msd}\ ((\mu)_2)_{0,1}\ ,\ ((\mu)_1)_0)) = 1$ | (236., 256.; ch.in.), |
| 262. | $\overline{sg}\ (\mathrm{msd}\ ((\mu)_2)_{0,1}\ ,\ ((\mu)_1)_0)) = 1$ | $((\mu)_2)_{0,1} = ((\mu)_1)_0$ | […] |
| 263.-267. | $\mathrm{Ant}_2^{\ell_1}(\mu)$ | $((\mu)_2)_{0,1} = ((\mu)_1)_0$ | (261., 262.; ch.in.), |

268.   $\overline{sg}\ (\mathrm{msd}\ ((\mu)_2)_{0,1}\ ,\ ((\mu)_1)_0)) = 1$

$\overline{sg}\ (\mathrm{msd}\ ((\mu)_2)_{0,0}\ ,\ 3)) \cdot \overline{sg}\ (\mathrm{msd}\ ((\mu)_2)_{0,1}\ ,\ ((\mu)_1)_0)) =$

$\overline{sg}\ (\mathrm{msd}\ ((\mu)_2)_{0,0}\ ,\ 3)) \cdot 1$ \hfill (LEA-rp),

269.-273.   $\mathrm{Ant}_2^{\ell_1}(\mu)$

$\overline{sg}\ (\mathrm{msd}\ ((\mu)_2)_{0,0}\ ,\ 3)) \cdot \overline{sg}\ (\mathrm{msd}\ ((\mu)_2)_{0,1}\ ,\ ((\mu)_1\ )_0)) =$

$\overline{sg}\ (\mathrm{msd}\ ((\mu)_2)_{0,0}\ ,\ 3)) \cdot 1$ \hfill (261., 268.; ch.in.),

274.   $\overline{sg}\ (\mathrm{msd}\ ((\mu)_2)_{0,0}\ ,\ 3)) \cdot 1 = \overline{sg}\ (\mathrm{msd}\ ((\mu)_2)_{0,0}\ ,\ 3))$ \hfill […],

275.   274.   276. \hfill (LEA$_2^{\bar{=}}$),

276.   $\overline{sg}\ (\mathrm{msd}\ ((\mu)_2)_{0,0}\ ,\ 3)) \cdot \overline{sg}\ (\mathrm{msd}\ ((\mu)_2)_{0,1}\ ,\ ((\mu)_1\ )_0)) =$

$\overline{sg}\ (\mathrm{msd}\ ((\mu)_2)_{0,0}\ ,\ 3)) \cdot 1$

$\overline{sg}\ (\mathrm{msd}\ ((\mu)_2)_{0,0}\ ,\ 3)) \cdot \overline{sg}\ (\mathrm{msd}\ ((\mu)_2)_{0,1}\ ,\ ((\mu)_1\ )_0)) =$

$\overline{sg}\ (\mathrm{msd}\ ((\mu)_2)_{0,0}\ ,\ 3))$ \hfill (274., 275.; MP),

277.-281.   $\mathrm{Ant}_2^{\ell_1}(\mu)$

$\overline{sg}\ (\mathrm{msd}\ ((\mu)_2)_{0,0}\ ,\ 3)) \cdot \overline{sg}\ (\mathrm{msd}\ ((\mu)_2)_{0,1}\ ,\ ((\mu)_1\ )_0)) =$

$\overline{sg}\ (\mathrm{msd}\ ((\mu)_2)_{0,0}\ ,\ 3))$ \hfill (273., 276.; ch.in.),

282.   $\overline{sg}\ (\mathrm{msd}\ ((\mu)_2)_{0,0}\ ,\ 3)) \cdot \overline{sg}\ (\mathrm{msd}\ ((\mu)_2)_{0,1}\ ,\ ((\mu)_1\ )_0)) =$

$$\overline{sg}\,(msd\,((\mu)_2)_{0,0},\,3))$$

| | | |
|---|---|---|
| | (142. $((\mu)_1,(\mu)_2) = \overline{sg}\,(msd\,((\mu)_2)_{0,0},\,3))$ | (LEA$_2^=$), |
| 283.-287. | Ant$_2^{b_1}(\mu)$ (142. $((\mu)_1,(\mu)_2) = \overline{sg}\,(msd\,((\mu)_2)_{0,0},\,3))$ | (281., 282.; ch.in.), |
| 288.-292. | Ant$_2^{b_1}(\mu)$ $((\mu)_1,(\mu)_2) = \overline{sg}\,(msd\,((\mu)_2)_{0,0},\,3))$ | (287., 142.; ch.in.$_2$), |
| 293. | $\overline{sg}\,(msd\,((\mu)_2)_{0,0},\,3)) = 0$  $\overline{sg}\,(msd\,((\mu)_2)_{0,0},\,3)) = 1$ | [...], |
| 294.-313. | $\neg\,\overline{sg}\,(msd\,((\mu)_2)_{0,0},\,3)) = 0$  $\overline{sg}\,(msd\,((\mu)_2)_{0,0},\,3)) = 1$ | (293.; Mtp$_2$), |
| 314. | $\overline{sg}\,(msd\,((\mu)_2)_{0,0},\,3)) = 0$  ( $((\mu)_1,(\mu)_2) = \overline{sg}\,(msd\,((\mu)_2)_{0,0},\,3))$ | |
| | $((\mu)_1,(\mu)_2) = 0)$ | (LEA$_2^=$), |
| 315.-322. | $((\mu)_1,(\mu)_2) = \overline{sg}\,(msd\,((\mu)_2)_{0,0},\,3))$ | |
| | ($\overline{sg}\,(msd\,((\mu)_2)_{0,0},\,3) = 0$  $((\mu)_1,(\mu)_2) = 0)$ | (314.; Int-Ant.), |
| 323.-327. | Ant$_2^{b_1}(\mu)$ | |
| | ($\overline{sg}\,(msd\,((\mu)_2)_{0,0},\,3)) = 0$  $((\mu)_1,(\mu)_2) = 0)$ | (292., 322.; ch.in.), |
| 328.-342. | ($\overline{sg}\,(msd\,((\mu)_2)_{0,0},\,3)) = 0$  $((\mu)_1,(\mu)_2) = 0)$ | |
| | ($\neg\ ((\mu)_1,(\mu)_2) = 0)$  $\neg\,\overline{sg}\,(msd\,((\mu)_2)_{0,0},\,3)) = 0)$ | (contrap), |
| 343.-347. | Ant$_2^{b_1}(\mu)$  ($\neg\ ((\mu)_1,(\mu)_2) = 0)$  $\neg\,\overline{sg}\,(msd\,((\mu)_2)_{0,0},\,3)) = 0)$ | (327., 342.; ch.in.), |
| 348. | 347. (219. 350.) | (Imp 2), |
| 349. | 219.  350. | (347., 348.; MP), |
| 350. | Ant$_2^{b_1}(\mu)$  $\neg\,\overline{sg}\,(msd\,((\mu)_2)_{0,0},\,3)) = 0$ | (219., 349.; MP), |
| 351.-355. | Ant$_2^{b_1}(\mu)$  $\overline{sg}\,(msd\,((\mu)_2)_{0,0},\,3)) = 1$ | (350., 313.; ch.in.), |
| 356. | $\overline{sg}\,(msd\,((\mu)_2)_{0,0},\,3)) = 1$  $((\mu)_2)_{0,0} = 3$ | [...], |
| 357.-361. | Ant$_2^{b_1}(\mu)$  $((\mu)_2)_{0,0} = 3$ | (355., 356.; ch.in.), |

362. $\overline{sg}$ (msd $(((\mu)_2)_{0,0}, 3)) = 1$ ( $((\mu)_1, (\mu)_2) = \overline{sg}$ (msd $(((\mu)_2)_{0,0}, 3))$

$((\mu)_1, (\mu)_2) = 1)$ (LEA$_{\overline{2}}$),

363.-367. Ant$_2^{\ell_1}(\mu)$

( $((\mu)_1, (\mu)_2) = \overline{sg}$ (msd $(((\mu)_2)_{0,0}, 3))$

$((\mu)_1, (\mu)_2) = 1)$ (355., 362.; ch.in.),

368. 367. (292. 370.) (Imp 2),

369. 292. 370. (367., 368.; MP),

370. Ant$_2^{\ell_1}(\mu)$ $((\mu)_1, (\mu)_2) = 1$ (292., 369.; MP),

371. $((\mu)_1, (\mu)_2) = 1$

Mp $((\mu)_1, (\mu)_2) \cdot$ $((\mu)_1, (\mu)_2) =$ Mp $((\mu)_1, (\mu)_2) \cdot 1$ (LEA-rp),

372.-376. Ant$_2^{\ell_1}(\mu)$ Mp $((\mu)_1, (\mu)_2) \cdot$ $((\mu)_1, (\mu)_2) =$ Mp $((\mu)_1, (\mu)_2) \cdot 1$ (370., 371; ch.in.),

377. Mp $((\mu)_1, (\mu)_2) \cdot 1 =$ Mp $((\mu)_1, (\mu)_2)$ […],

378. 377. 379. (LEA$_{\overline{2}}$),

379. Mp $((\mu)_1, (\mu)_2) \cdot$ $((\mu)_1, (\mu)_2) =$ Mp $((\mu)_1, (\mu)_2) \cdot 1$

Mp $((\mu)_1, (\mu)_2) \cdot$ $((\mu)_1, (\mu)_2) =$ Mp $((\mu)_1, (\mu)_2)$ (377., 379.; MP),

380.-384. Ant$_2^{\ell_1}(\mu)$

Mp $((\mu)_1, (\mu)_2) \cdot$ $((\mu)_1, (\mu)_2) =$ Mp $((\mu)_1, (\mu)_2)$ (376., 379.; ch.in.),

385. Mp $((\mu)_1, (\mu)_2) \cdot$ $((\mu)_1, (\mu)_2) =$ Mp $((\mu)_1, (\mu)_2)$

$(\mu =$ Mp $((\mu)_1, (\mu)_2) \cdot$ $((\mu)_1, (\mu)_2)$ $\mu =$ Mp $((\mu)_1, (\mu)_2))$ (LEA$_{\overline{2}}$),

| | | | |
|---|---|---|---|
| 386.-390. | $Ant_2^{\ell_1}(\mu)$ | | |
| | ($\mu = Mp\,((\mu)_1, (\mu)_2) \cdot \quad ((\mu)_1, (\mu)_2) \qquad \mu = Mp\,((\mu)_1, (\mu)_2))$ | | (384., 385.; ch.in.), |
| 391. | 390.  (141.  393.) | | (Imp 2), |
| 392. | 141.  393. | | (390., 391.; MP), |
| 393. | $Ant_2^{\ell_1}(\mu) \quad \mu = Mp\,((\mu)_1, (\mu)_2)$ | | (141., 392.; MP), |
| 394. | $Mp\,((\mu)_1, (\mu)_2) = 2^{((\mu)_2)_{0,2}} \cdot 3^{(\mu)_1} \cdot 5^{(\mu)_2}$ | | ($\mathbf{def_{Mp}}$), |
| 395. | 394.  396. | | ($LEA_{\bar{2}}^{=}$), |
| 396. | $\mu = Mp\,((\mu)_1, (\mu)_2) \qquad \mu = 2^{((\mu)_2)_{0,2}} \cdot 3^{(\mu)_1} \cdot 5^{(\mu)_2}$ | | (394., 395.; MP), |
| 397.-401. | $Ant_2^{\ell_1}(\mu) \qquad \mu = 2^{((\mu)_2)_{0,2}} \cdot 3^{(\mu)_1} \cdot 5^{(\mu)_2}$ | | (393., 396.; ch.in.), |
| 402. | $\mu = 2^{((\mu)_2)_{0,2}} \cdot 3^{(\mu)_1} \cdot 5^{(\mu)_2} \qquad (\mu)_0 = (2^{((\mu)_2)_{0,2}} \cdot 3^{(\mu)_1} \cdot 5^{(\mu)_2})_0$ | | (LEA-rp), |
| 403.-406. | $Ant_2^{\ell_1}(\mu) \qquad (\mu)_0 = (2^{((\mu)_2)_{0,2}} \cdot 3^{(\mu)_1} \cdot 5^{(\mu)_2})_0$ | | (401., 402.; ch.in.), |
| 407. | $(2^{((\mu)_2)_{0,2}} \cdot 3^{(\mu)_1} \cdot 5^{(\mu)_2})_0 = ((\mu)_2)_{0,2}$ | | ($Ga_0$), |
| 408. | 407.  409. | | ($LEA_{\bar{2}}^{=}$), |
| 409. | $(\mu)_0 = (2^{((\mu)_2)_{0,2}} \cdot 3^{(\mu)_1} \cdot 5^{(\mu)_2})_0 \qquad (\mu)_0 = ((\mu)_2)_{0,2}$ | | (407., 408.; MP), |
| 410.-413. | $Ant_2^{\ell_1}(\mu) \qquad (\mu)_0 = ((\mu)_2)_{0,2}$ | | (406., 409.; ch.in.), |
| 414. | $(\mu)_0 = ((\mu)_2)_{0,2} \qquad ((\mu)_2)_{0,2} = (\mu)_0$ | | ($Sym_=$), |
| 415.-419. | $Ant_2^{\ell_1}(\mu) \qquad ((\mu)_2)_{0,2} = (\mu)_0$ | | (413., 414.; ch.in.), |
| 420. | $\ell_1((\mu)_2) = 0 \qquad cfor\,((\mu)_2)_0) = 0$ | | [...], |
| 421. | $cfor\,((\mu)_2)_0 = 0 \qquad \ell((\mu)_2)_0) = 2$ | | [...], |
| 422.-426. | $\ell_1((\mu)_2) = 0$ | $\ell((\mu)_2)_0) = 2$ | (420., 421.; ch.in.), |
| 427. | $\ell((\mu)_2)_0) = 2$ | $((\mu)_2)_0 = 2^{((\mu)_2)_0)_0} \cdot 3^{(((\mu)_2)_0)_1} \cdot 5^{(((\mu)_2)_0)_2}$ | [...], |
| 428.-432. | $\ell_1((\mu)_2) = 0$ | $((\mu)_2)_0 = 2^{((\mu)_2)_0)_0} \cdot 3^{(((\mu)_2)_0)_1} \cdot 5^{(((\mu)_2)_0)_2}$ | (426., 427.; ch.in.), |
| 433.-437. | $\mathrm{Ant}_2^{\ell_1}(\mu)$ | $((\mu)_2)_0) = 2^{((\mu)_2)_0)_0} \cdot 3^{(((\mu)_2)_0)_1} \cdot 5^{(((\mu)_2)_0)_2}$ | (8a., 432.; ch.in.), |
| 438. | $(((\mu)_2)_0)_0 = ((\mu)_2)_{0,0}$ | | ($\mathbf{def}^2{}_{\exp_{0,0}}$), |
| 439. | $((\mu)_2)_{0,0} = 3$ | (438. $(((\mu)_2)_0)_0 = 3$) | (LEA$_{\bar 2}$), |
| 440.-444. | $((\mu)_2)_{0,0} = 3$ | $(((\mu)_2)_0)_0 = 3$ | (439., 438.; ch.in.$_2$), |
| 445. | $(((\mu)_2)_0)_0 = 3$ | $2^{(((\mu)_2)_0)_0} = 2^3$ | (LEA-rp), |
| 446.-450. | $((\mu)_2)_{0,0} = 3$ | $2^{(((\mu)_2)_0)_0} = 2^3$ | (444., 445.; ch.in.), |
| 451.-455. | $\mathrm{Ant}_2^{\ell_1}(\mu)$ | $2^{(((\mu)_2)_0)_0} = 2^3$ | (361., 450.; ch.in.), |
| 456. | $(((\mu)_2)_0)_1 = ((\mu)_2)_{0,1}$ | | ($\mathbf{def}^2{}_{\exp_{0,1}}$), |
| 457. | $((\mu)_2)_{0,1} = ((\mu)_1)_0$ | (456. $(((\mu)_2)_0)_1 = ((\mu)_1)_0$) | (LEA$_{\bar 2}$), |
| 458.-462. | $((\mu)_2)_{0,1} = ((\mu)_1)_0$ | $(((\mu)_2)_0)_1 = ((\mu)_1)_0$ | (457., 456.; ch.in.$_2$), |
| 463. | $(((\mu)_2)_0)_1 = ((\mu)_1)_0$ | $3^{(((\mu)_2)_0)_1} = 3^{((\mu)_1)_0}$ | (LEA-rp,) |
| 464.-468. | $((\mu)_2)_{0,1} = ((\mu)_1)_0$ | $3^{(((\mu)_2)_0)_1} = 3^{((\mu)_1)_0}$ | (462., 463.; ch.in.), |
| 469.-473. | $\mathrm{Ant}_2^{\ell_1}(\mu)$ | $3^{(((\mu)_2)_0)_1} = 3^{((\mu)_1)_0}$ | (267., 468.; ch.in.), |
| 474. | $2^{(((\mu)_2)_0)_0} = 2^3$ | $2^{(((\mu)_2)_0)_0} \cdot 3^{(((\mu)_2)_0)_1} = 2^3 \cdot 3^{(((\mu)_2)_0)_1}$ | (LEA-rp), |
| 475.-479. | $\mathrm{Ant}_2^{\ell_1}(\mu)$ | $2^{(((\mu)_2)_0)_0} \cdot 3^{(((\mu)_2)_0)_1} = 2^3 \cdot 3^{(((\mu)_2)_0)_1}$ | (455., 474.; ch.in.), |
| 480. | $3^{(((\mu)_2)_0)_1} = 3^{((\mu)_1)_0}$ | $2^3 \cdot 3^{(((\mu)_2)_0)_1} = 2^3 \cdot 3^{((\mu)_1)_0}$ | (LEA-rp), |
| 481.-485. | $\text{Ant}_2^{b_1}(\mu)$ | $2^3 \cdot 3(((\mu)_2)_0)_1 = 2^3 \cdot 3((\mu)_1)_0$ | (473., 480.; ch.in.), |
| 486. | $2^3 \cdot 3(((\mu)_2)_0)_1 = 2^3 \cdot 3((\mu)_1)_0$ | | |
| | $(2(((\mu)_2)_0)_0 \cdot 3(((\mu)_2)_0)_1 = 2^3 \cdot 3(((\mu)_2)_0)_1$ | | |
| | $2(((\mu)_2)_0)_0 \cdot 3(((\mu)_2)_0)_1 = 2^3 \cdot 3((\mu)_1)_0)$ | | $(\text{LEA}_{\overline{2}})$, |
| 487.-491. | $\text{Ant}_2^{b_1}(\mu)$ | | |
| | $(2(((\mu)_2)_0)_0 \cdot 3(((\mu)_2)_0)_1 = 2^3 \cdot 3(((\mu)_2)_0)_1$ | | |
| | $2(((\mu)_2)_0)_0 \cdot 3(((\mu)_2)_0)_1 = 2^3 \cdot 3((\mu)_1)_0)$ | | (485., 486.; ch.in.), |
| 492. | 491.    (479.    494.) | | (Imp 2), |
| 493. | 479.    494. | | (491., 492.; MP), |
| 494. | $\text{Ant}_2^{b_1}(\mu)$ | $2(((\mu)_2)_0)_0 \cdot 3(((\mu)_2)_0)_1 = 2^3 \cdot 3((\mu)_1)_0$ | (479., 493.; MP), |
| 495. | $(((\mu)_2)_0)_2 = ((\mu)_2)_{0,2}$ | | $(\mathbf{def}^2_{\exp_{0,2}})$, |
| 496. | $((\mu)_2)_{0,2} = (\mu)_0$    (495.    $(((\mu)_2)_0)_2 = (\mu)_0)$ | | $(\text{LEA}_{\overline{2}})$, |
| 497.-501. | $((\mu)_2)_{0,2} = (\mu)_0$    $(((\mu)_2)_0)_2 = (\mu)_0$ | | (496., 495.; ch.in.$_2$), |
| 502. | $(((\mu)_2)_0)_2 = (\mu)_0$    $5^{(((\mu)_2)_0)_2} = 5^{(\mu)_0}$ | | (LEA-rp), |
| 503.-507. | $((\mu)_2)_{0,2} = (\mu)_0$    $5^{(((\mu)_2)_0)_2} = 5^{(\mu)_0}$ | | (501., 502.; ch.in.), |
| 508.-512. | $\text{Ant}_2^{b_1}(\mu)$    $5^{(((\mu)_2)_0)_2} = 5^{(\mu)_0}$ | | (419., 507.; ch.in.), |
| 513. | $5^{(((\mu)_2)_0)_2} = 5^{(\mu)_0}$    $2^3 \cdot 3((\mu)_1)_0 \cdot 5^{(((\mu)_2)_0)_2} = 2^3 \cdot 3((\mu)_1)_0 \cdot 5^{(\mu)_0}$ | | (LEA-rp), |
| 514.-518. | $\text{Ant}_2^{b_1}(\mu)$    $2^3 \cdot 3((\mu)_1)_0 \cdot 5^{(((\mu)_2)_0)_2} = 2^3 \cdot 3((\mu)_1)_0 \cdot 5^{(\mu)_0}$ | | (512., 513.; ch.in.), |
| 519. | $2(((\mu)_2)_0)_0 \cdot 3(((\mu)_2)_0)_1 = 2^3 \cdot 3((\mu)_1)_0$ | | |

$$2^{(((\mu)_2)_0)_0} \cdot 3^{(((\mu)_2)_0)_1} \cdot 5^{(((\mu)_2)_0)_2} = 2^3 \cdot 3^{((\mu)_1)_0} \cdot 5^{(((\mu)_2)_0)_2} \qquad \text{(LEA-rp)},$$

520.-524. $\quad \text{Ant}_2^{\ell_1}(\mu)$

$$2^{(((\mu)_2)_0)_0} \cdot 3^{(((\mu)_2)_0)_1} \cdot 5^{(((\mu)_2)_0)_2} = 2^3 \cdot 3^{((\mu)_1)_0} \cdot 5^{(((\mu)_2)_0)_2} \qquad \text{(494., 519.; ch.in.)},$$

525. $\quad 2^3 \cdot 3^{((\mu)_1)_0} \cdot 5^{(((\mu)_2)_0)_2} = 2^3 \cdot 3^{((\mu)_1)_0} \cdot 5^{(\mu)_0}$

$\quad\quad\quad (2^{(((\mu)_2)_0)_0} \cdot 3^{(((\mu)_2)_0)_1} \cdot 5^{(((\mu)_2)_0)_2} = 2^3 \cdot 3^{((\mu)_1)_0} \cdot 5^{(((\mu)_2)_0)_2}$

$\quad\quad\quad (2^{(((\mu)_2)_0)_0} \cdot 3^{(((\mu)_2)_0)_1} \cdot 5^{(((\mu)_2)_0)_2} = 2^3 \cdot 3^{((\mu)_1)_0} \cdot 5^{(\mu)_0}) \qquad (\text{LEA}_{\overline{2}})$,

526.-530. $\quad \text{Ant}_2^{\ell_1}(\mu)$

$\quad\quad\quad (2^{(((\mu)_2)_0)_0} \cdot 3^{(((\mu)_2)_0)_1} \cdot 5^{(((\mu)_2)_0)_2} = 2^3 \cdot 3^{((\mu)_1)_0} \cdot 5^{(((\mu)_2)_0)_2}$

$\quad\quad\quad (2^{(((\mu)_2)_0)_0} \cdot 3^{(((\mu)_2)_0)_1} \cdot 5^{(((\mu)_2)_0)_2} = 2^3 \cdot 3^{((\mu)_1)_0} \cdot 5^{(\mu)_0}) \qquad (\text{518., 525.; ch.in.})$,

531. $\quad\quad\quad$ 530. $\quad$ (524. $\quad$ 533.) $\qquad\qquad\qquad\qquad$ (Imp 2),

532. $\quad\quad\quad$ 524. $\quad$ 533. $\qquad\qquad\qquad\qquad\qquad\qquad$ (530., 531.; MP),

533. $\quad \text{Ant}_2^{\ell_1}(\mu)$

$\quad\quad\quad 2^{(((\mu)_2)_0)_0} \cdot 3^{(((\mu)_2)_0)_1} \cdot 5^{(((\mu)_2)_0)_2} = 2^3 \cdot 3^{((\mu)_1)_0} \cdot 5^{(\mu)_0} \qquad (\text{524., 532.; MP})$,

534. $\quad 2^{(((\mu)_2)_0)_0} \cdot 3^{(((\mu)_2)_0)_1} \cdot 5^{(((\mu)_2)_0)_2} = 2^3 \cdot 3^{((\mu)_1)_0} \cdot 5^{(\mu)_0}$

$\quad\quad\quad (((\mu)_2)_0 = 2^{(((\mu)_2)_0)_0} \cdot 3^{(((\mu)_2)_0)_1} \cdot 5^{(((\mu)_2)_0)_2}$

$\quad\quad\quad ((\mu)_2)_0 = 2^3 \cdot 3^{((\mu)_1)_0} \cdot 5^{(\mu)_0}) \qquad (\text{LEA}_{\overline{2}})$,

535.-539. $\quad \text{Ant}_2^{\ell_1}(\mu)$

$\quad\quad\quad (((\mu)_2)_0 = 2^{(((\mu)_2)_0)_0} \cdot 3^{(((\mu)_2)_0)_1} \cdot 5^{(((\mu)_2)_0)_2}$

|   |   |   |
|---|---|---|
|  | $((\mu)_2)_0 = 2^3 \cdot 3^{((\mu)_1)_0} \cdot 5^{(\mu)_0}$ | (533., 534.; ch.in.), |
| 540. | 539. (437. 542.) | (Imp 2), |
| 541. | 437. 542. | (539., 540; MP), |
| 542. | $\text{Ant}_2^{\ell_1}(\mu)$ $((\mu)_2)_0 = 2^3 \cdot 3^{((\mu)_1)_0} \cdot 5^{(\mu)_0}$ | (437., 541.; MP), |
| 543. | $((\mu)_2)_0 = 2^3 \cdot 3^{((\mu)_1)_0} \cdot 5^{(\mu)_0}$ $(F_{((\mu)_2)_0}\ F_{2^3 \cdot 3^{((\mu)_1)_0} \cdot 5^{(\mu)_0}})$ | $(LEA_1^\nu)$, |
| 544.–548. | $\text{Ant}_2^{\ell_1}(\mu)$ $(F_{((\mu)_2)_0}\ F_{2^3 \cdot 3^{((\mu)_1)_0} \cdot 5^{(\mu)_0}})$ | (542., 543.; ch.in.), |
| 549. | $(2^3 \cdot 3^{((\mu)_1)_0} \cdot 5^{(\mu)_0})_0 = 3$ | $(Ga_0)$, |
| 550. | $(2^3 \cdot 3^{((\mu)_1)_0} \cdot 5^{(\mu)_0})_1 = ((\mu)_1)_0$ | $(Ga_1)$, |
| 551. | $(2^3 \cdot 3^{((\mu)_1)_0} \cdot 5^{(\mu)_0})_2 = (\mu)_0$ | $(Ga_2)$, |
| 552. | 549. $(F_{2^3 \cdot 3^{((\mu)_1)_0} \cdot 5^{(\mu)_0}}$ $F_{(2^3 \cdot 3^{((\mu)_1)_0} \cdot 5^{(\mu)_0})_1}\ F_{(2^3 \cdot 3^{((\mu)_1)_0} \cdot 5^{(\mu)_0})_2})$ | $(LEA_{MP}^\nu)$, |
| 553. | $F_{2^3 \cdot 3^{((\mu)_1)_0} \cdot 5^{(\mu)_0}}$ $(F_{(2^3 \cdot 3^{((\mu)_1)_0} \cdot 5^{(\mu)_0})_1}\ F_{(2^3 \cdot 3^{((\mu)_1)_0} \cdot 5^{(\mu)_0})_2})$ | (549., 552.; MP), |
| 554. | 550. 555. | $(Sym_=)$, |
| 555. | $((\mu)_1)_0 = (2^3 \cdot 3^{((\mu)_1)_0} \cdot 5^{(\mu)_0})_1$ | (550, 554.; MP), |
| 556. | 555. 557. | $(LEA_1^\nu)$, |
| 557. | $F_{((\mu)_1)_0}$ $F_{(2^3 \cdot 3^{((\mu)_1)_0} \cdot 5^{(\mu)_0})_1}$ | (555., 556.; MP), |
| 558.–565. | $F_{(2^3 \cdot 3^{((\mu)_1)_0} \cdot 5^{(\mu)_0})_1}$ $(F_{2^3 \cdot 3^{((\mu)_1)_0} \cdot 5^{(\mu)_0}}\ F_{(2^3 \cdot 3^{((\mu)_1)_0} \cdot 5^{(\mu)_0})_2})$ | (553.; Int-Ant), |
| 566.–570. | $F_{((\mu)_1)_0}$ $(F_{2^3 \cdot 3^{((\mu)_1)_0} \cdot 5^{(\mu)_0}}\ F_{(2^3 \cdot 3^{((\mu)_1)_0} \cdot 5^{(\mu)_0})_2})$ | (557., 565.; ch.in.), |
| 571. | 551. 572. | $(LEA_1^\nu)$, |
| 572. | $F_{(2^3 \cdot 3^{((\mu)_1)_0} \cdot 5^{(\mu)_0})_2}$ $F_{(\mu)_0}$ | (551., 571.; MP), |

| | | | | | | |
|---|---|---|---|---|---|---|
| 573.-579. | 572. | 580. | | | | (ch.in. - fla$_2$), |
| 580. | $(F_2 3 \cdot {}_3((\mu)_1)_0 \cdot {}_5(\mu)_0$ | $F_{(2^3 \cdot {}_3((\mu)_1)_0 \cdot {}_5(\mu)_0)_2})$ | | | | |
| | $(F_2 3 \cdot {}_3((\mu)_1)_0 \cdot {}_5(\mu)_0$ | $F_{(\mu)_0})$ | | | | (572., 579.; MP), |
| 581.-585. | $F_{((\mu)_1)_0}$ | $(F_2 3 \cdot {}_3((\mu)_1)_0 \cdot {}_5(\mu)_0$ | $F_{(\mu)_0})$ | | | (570., 580.; ch.in.), |
| 585.-593. | $F_2 3 \cdot {}_3((\mu)_1)_0 \cdot {}_5(\mu)_0$ | $(F_{((\mu)_1)_0}$ | $F_{(\mu)_0})$ | | | (585.; Int-Ant.), |
| 594.-600. | 593. | 601. | | | | […], |
| 601. | $(F_{((\mu)_2)_0}$ | $F_2 3 \cdot {}_3((\mu)_1)_0 \cdot {}_5(\mu)_0)$ | $(F_{((\mu)_2)_0}$ | $(F_{((\mu)_1)_0}$ | $F_{(\mu)_0}))$ | (593., 600.; MP), |
| 602.-606. | $\text{Ant}_2^{\ell_1}(\mu)$ | $(F_{((\mu)_2)_0}$ | $(F_{((\mu)_1)_0}$ | $F_{(\mu)_0}))$ | | (548., 601.; ch.in.), |
| 607.-625. | $(F_{((\mu)_2)_0}$ | $(F_{((\mu)_1)_0}$ | $F_{((\mu)_0)}))$ | $(F_{((\mu)_1)_0}$ | $(F_{((\mu)_2)_0}$ $F_{(\mu)_0}))$ | ( int-ant.), |
| 626.-630. | $\text{Ant}_2^{\ell_1}(\mu)$ | $(F_{((\mu)_1)_0}$ | $(F_{((\mu)_2)_0}$ | $F_{(\mu)_0}))$ | | (606.-625.; ch.in.), |
| 631. | 7. | (630. | 633.) | | | (mp-ind), |
| 632. | 630. | 633. | | | | (7., 631.; MP), |
| 633. | $\ell_1(\mu) = 0$ | $F_{(\mu)_0}$ | | | | (630., 632.; MP) |

| | | |
|---|---|---|
| 633. | $\mathscr{b}_1(\ ) = 0 \quad F_{(\ )_0}$ | (633.; μ  ), |
| 634. | $\mathscr{b}(\ ,\nu_\ell) = 0\ \&\ \neg\ (\ <\ \&\ \mathscr{b}(\ ,\nu_\ell) = 0)\ .\lor.$ | |
| | $\quad = 0\ \&\ \neg\ (\mathscr{b}(\ ,\nu_\ell) = 0),$ | (M- ), |
| 635. | $\mathscr{b}(\ ,\nu_\ell) = 0 \quad \mathscr{b}_1(\ ) = 0\ \&\ (\ )_0 = \nu_\ell$ | (**def**$_1^{\mathscr{b}}$), |
| 636. | $\mathscr{b}_1(\ ) = 0\ \&\ (\ )_0 = \nu_\ell \quad \mathscr{b}_1(\ ) = 0$ | (Con 1), |
| 637.-641. | $\mathscr{b}(\ ,\nu_\ell) = 0 \quad \mathscr{b}_1(\ ) = 0$ | (635., 636.; ch.in.), |
| 642.-646 | $\mathscr{b}(\ ,\nu_\ell) = 0 \quad F_{(\ )_0}$ | (641., 633. ; ch.in.), |
| 647. | $\mathscr{b}_1(\ ) = 0\ \&\ (\ )_0 = \nu_\ell \quad (\ )_0 = \nu_\ell$ | (Con 2), |
| 648.-552. | $\mathscr{b}(\ ,\nu_\ell) = 0 \quad (\ )_0 = \nu_\ell$ | (635., 647.; ch.in.), |
| 653. | $(\ )_0 = \nu_\ell \quad (F_{(\ )_0} \quad F_{\nu_\ell})$ | (LEA$_1^\nu$), |
| 654.-658. | $\mathscr{b}(\ ,\nu_\ell) = 0 \quad (F_{(\ )_0} \quad F_{\nu_\ell})$ | (652., 653.; ch.in.), |
| 659. | 658. (646  661.) | (Imp 2), |
| 660. | 646  660. | (658., 659.; MP), |
| 661. | $\mathscr{b}(\ ,\nu_\ell) = 0 \quad F_{\nu_\ell}$ | (646, 660; MP), |
| 662. | $F_{\nu_\ell}\ \ell$ | [...], |
| 663.-667. | $\mathscr{b}(\ ,\nu_\ell) = 0\ \ \ell$ | (661., 662; ch.in.), |
| 668. | $\mathscr{b}(\ ,\nu_\ell) = 0\ \&\ \neg\ (\ <\ \&\ \mathscr{b}(\ ,\nu_\ell) = 0)$ | |
| | $\mathscr{b}(\ ,\nu_\ell) = 0$ | (Con 1), |
| 669.-673. | ℓ ( ,ν_ℓ) = 0 & ¬ ( < & ℓ ( ,ν_ℓ) = 0) ℓ | (668., 667.; ch.in.), |
| 674. | ℓ ¬ ℓ ( ,ν_ℓ) = 0 | (Imp 3), |
| 675.-679. | ℓ ( ,ν_ℓ) = 0 & ¬ ( < & ℓ ( ,ν_ℓ) = 0) | |
| | ¬ ℓ ( ,ν_ℓ) = 0 | (673., 674.; ch.in.), |
| 680. | = 0 & ¬ ℓ ( ,ν_ℓ) = 0 ¬ ℓ ( ,ν_ℓ) = 0 | (Con 2), |
| 681. | 679 (680. 683.) | (Con 3), |
| 682. | 679. 682. | (679., 681.; MP), |
| 683. | 679. & 680. | (680., 682.; MP), |
| 684. | 683 (634. 686.) | (Dis 3), |
| 685. | 634. 686. | (683., 684; MP), |
| 686. | ¬ ℓ ( ,ν_ℓ) = 0 | (634., 685.; MP), |
| 686a. | Con_Ari | [**def** Con_Ari ]. |

(End of proof).

# ENDNOTE A.

a1. In this work the formal objects considered as numeroids are also denoted by the usual arabic figures 0, 1, 2, ..., while integers, considered as intuitive objets, shall, *in case of need*, be denoted by <u>underlining</u> the same figures, i.e., by the symbols $\underline{0}, \underline{1}, \underline{2}, \ldots$ . In MAri, it shall be assumed that these two sequences of symbols are isomorphic so that for any integer $\underline{n}$ there corresponds a numeroid, n ( i.e. $Su^n 0$), which is said to stand for the *value* of the former. Further, for any constant termoid t in LAri, its *value*, $\nu t$, shall be the numeroid $Su^t 0$. This does not mean that the equality $t = \nu t$ is necessarily provable in Ari — in fact, it even looks impossible to prove this equality for arbitrary constant termoid t in LAri (unless Ari is inconsistent).

a2. It will be convenient to have several distinct notations for the same En-N depending on the context in which it occurs. Thus, for example, if E denotes an equality, $t = \nu t$, it shall be assigned the En-N $2^{15} \cdot 3^{\tau_t} \cdot 523^{suc(t)}$, in accordance with the general stipulations for assigning En-N's to (closed) formulae, where $\tau_t$ is a special notation for En-N of (closed) termoid t. However, neither this notation nor the object which it denotes belong to the language LAri — in which the corresponding object would be exceedingly complex. Thus, already the first factor, $2^{15}$, would, in LAri, be quite inconvenient to write down more than once as it is actually the object

ex(SuSu0, SuSuSuSuSuSuSuSuSuSuSuSuSuSu0)

where ex stands for the function symbol of exponentiation in LAri the value of which in the present case would correspond to an expression consisting of $2^{15}$ occurrences of the symbol Su.
Therefore, the notation $\nu_E$, like $\tau_t$, above, shall be considered as only a *special* notation which may be used in place of a general notation for explicitly exhibiting the formal object denoted by $\nu_E$ — which, obviously, is more than unpractical to write down even once. Thus, certain contexts in this work shall require that the *explicit notation* in LAri for $\nu_E$ can be indicated — in which case the notation $\pi_E$ shall indicate that very object explicitly (in LAri) and be refered to as an 'En-NG' or 'En-N-*generator*'. In this work, the latter notation shall be confined to Appendix **D** (and further used in Part II of this paper).

(This matter is ignored even in [Hilbert-Bernays, 1939] where the passage from the metatheoretically considered termoid like $2^{15} \cdot 3^{\tau_t} \cdot 523^{suc(t)}$, above, to the (notation for the) numeroid corresponding to the En-NG for the equality $t = \nu t$ is considered unproblematical (cf. the penultimate remark in §5.2 of [Hilbert-Bernays, 1939])

However, even apart from feasability, the consideration of the provability of the equality $t = \nu t$, for arbitrary termoid t, *in* Ari is non-trivial. For example, it would be a major obstruction to deal with a demand that proscribed the use of any termoid t in LAri as a En-NG unless also the equality $t = \nu t$ has been proved in Ari. Let Ari* denote the formal system obtained from Ari by imposing this demand and let $\mathscr{E}^*$ denote the characteristic function of Ari*'s proof relation. The question is whether or not each En-N, $\nu_E$, in MAri, can be considered also as a En-N in MAri*? If so, then already in MAri* $\nu_E$ can be replaced by $Su^{\nu_E} 0$ so that for any constant termoids r, s, in LAri, $\mathscr{E}(r,s) = 0$ must hold in MAri iff it holds in MAri* and in MAri* it shall hold iff so does $\mathscr{E}(Su^r 0, Su^s 0) = 0$. Does this mean that the equality $\mathscr{E}(r,s) = \mathscr{E}(Su^r 0, Su^s 0)$ holds for all constant termoids r, s, in LAri including the case where r = and s = $\nu$? Can it hold without being Ari-provable?
Of course, if Ari is inconsistent then all equalities are provable — including = $\nu$ . The latter equality can also be written as $= \mathrm{Su}\ 0$ and if Ari is consistent then $= 0$ holds so that the right side of this equality is the numeroid with exactly 0 occurrences of Su, i.e. 0. Therefore, if Ari is consistent then this equality coincides with $= 0$ which turns out to be provably equivalent with the consistency formula $\mathrm{Con}_{\mathrm{Ari}}$. Thus, this equality is problematic as to its provability iff Ari is consistent.

The replacing of $\ell$ by $\ell^*$ might remove this problem — although that is an argument which illustrates how problematic it can be to carry through this replacement. The problem will be reconsidered in Part II of this work).

### ENDNOTE B.

In the open logic, the postulates of the WPC consists of the logical axioms of groups I. and II. together with the rules (MP) and (**Gen**). In the closed logic, the latter rule is excluded but otherwise the closed axioms of the open WPC are retained.

A remarkable property of WPC (with equality) is that it is rich enough to carry out, in the open logic, all derivations of the Predicate Calculus as exposited in Theorem 17 (and its Corollary) in §35. of [KLEENE, 1952]), *except* the formulae *79-*80. Also the propositional derivations can be carried out in the open logic, as well as derivations dealing with the replacement of parts of formulae by their equivalents (§33., Theorem 14, in [KLEENE, 1952]; all closures of formulae provable by such logical derivations are also provable in the WPC of the closed logic. (However, the formula *80. loses its holding unless the range of its variable, x, is supposed to contain a member; or else, the (closure of) *80 can be accepted as the axiom of the WPC which expresses that the range has a member. The renaming formula — i.e. the (closure of) *73, *74 in [KLEENE, 1952], §33, are not provable in the WPC.)

# APPENDIX **A**.

**1**.

Proof of (Imp0) where A denotes any closed formula in LAri:

1.  A ⊃ ( A ⊃ A )                                    ( I m p 1 ) ,

2.  A ⊃ ( ( A ⊃ A ) ⊃ A )                            ( I m p 1 ) ,

3.  2. ⊃ ( 1. ⊃ 5. )                                 ( I m p 2 ) ,

4.  1. ⊃ 5.                                          ( 2 . , 3 . ;  M P ) ,

5.  A ⊃ A                                            ( 1 . ,  4 . : M P )

**2**.

Below are collected the derived rules of inference together with the proof-number of the formula(e) at which it occurs.  <u>Chain inference 1.</u>

(ch.in.): $\dfrac{A \supset B, \quad B \supset C}{A \supset C}$

a) A ⊃ B;  b) B ⊃ C .

1.  b) ⊃ 2.              (Imp1),

2.  A ⊃ (B ⊃ C)          (b), 1.; MP),

3.  2. ⊃ 4.              (Imp2),

4.  a) ⊃ 5.              (2., 3.; MP),

5.  A ⊃ C                (a), 4.; MP).

Instances of ch.in.'s conclusions occur at formulae 7.;15.; 21.; 27.; 33.; 41.; 51.; 57.; 70.; 90; 118.; 127.; 132.; 138.; 152.; 158.; 177.; 183.; 196.; 216.; 236.; 261.; 267.; 273.; 281.; 287.; 327.; 347.; 355.; 361; 367.; 376.; 384.; 390.; 401.; 406. 413.; 419.; 426.; 432.; 437.; 450.; 455.; 468.; 473.; 479; 485.; 491.; 507.; 512.; 518.; 524.; 530.; 539.; 548.; 570.; 585.; 606.; 630.; 641.; 646.; 652.; 658.; 667.; 673. and 679.

<u>Chain inference 2.</u>

(ch.in.$_2$)  $\underline{A \supset (B \supset C) \quad \quad B}$
          $A \supset C$

a): $A \supset (B \supset C)$;   b): B.

1.  $B \supset (A \supset B)$           (Imp2),

2.  $A \supset B$                (b), 1.; MP),

3.  a)   4.               (Imp2),

4.  2.   5.               (a), 3.;MP),

5.  $A \supset C$            (2., 4.; MP).

Instances of ch.in$_2$'s conclusions occur at formulae 168.; 292.; 444; 462. and 501. and, additionally, in the derivation of the scheme contrap., below.

Chain inference 3.

(ch.in. - fla$_2$):   $(B \supset C) \supset ((A \supset B) \supset (A \supset C))$;

1.  $(B \supset C) \supset ((A \supset (B \supset C))$                          (Imp1),

2.  $A \supset (B \supset C) \supset ((A \supset B) \supset (B \supset C))$       (Imp2),

3.-7. $(B \supset C) \supset ((A \supset B) \supset (A \supset C))$              (1., 2.; ch.in.).

Instance of ch.in.-fla$_2$'s conclusion occurs at formula 579.

    Interchange of antecedent 1.

(Int-Ant):   $\underline{A \supset (B \supset C)}$
            $B \supset (A \supset C$

    a) $A \supset (B \supset C)$

1.  A → (B → C) → ((A → B) → (B → C))            (Imp2),

2.  (A → B) → (B → C)                            (Ia), 1.; MP),

3.  B → (A → B)                                   (Imp1),

4.- 8.  B → (A → C)                               (3., 2.; ch.in.).

Instances of (Int-Ant) occur at formula 65; 191; 322; 565.; 593.

Interchange of antecedent 2.

(int-ant):    (A → (B → C)) → (B → (A → C)))

1.  B → (A → B)                                                                      (Imp1),

2.  A → (B → C) → ((A → B) → (B → C))                                                (Imp2),

3.  ((A → B) → (A → C)) → B) → ((A → B) → (A → C))                                   (Imp1),

4.          3.      5.                                                               (Imp1),

5.  2.      A → (B → C)      3.                                                      (3., 4.; MP),

6.          5.      (2.     8.)                                                      (Imp2),

7.          2.      8.                                                               (5., 6.; MP),

8.  (A → (B → C)) → (B → ((A → B) → (A → C)))                                        (7., 2.; MP),

9.          (B → ((A → B) → (A → C))) → ((1. → (B → (A → C)))                        (Imp1),

10.         9.      11.                                                              (Imp1),

11.         ((A → (B → C)) → 9.                                                      (9., 10,; MP),

12.         11.     (8.    14.)                                                      (Imp2),

13.         8.      14.                                                              (11., 12,; MP),

14.         (A → (B → C)) → (1. → (B → (A → C)))                                     (8., 13; MP),

15.         1.      16.                                                              (Imp1),

| 16. | (A → (B → C)) → 1. | (1., 15; MP), |
| 17. | 14. → (16. → 19.) | (Imp2), |
| 18. | 16. → 19. | (14., 17; MP), |
| 19. | (A → (B → C)) → (B → (A → C))) | (16., 18; MP). |

Contraposition 1.

(Contrap): $\dfrac{A \to B}{\neg B \to \neg A}$

a) A → B

| 1.-7. | (B → ℓ) → ((A → B) → (A → ℓ)) | (ch.in.-fla$_2$), |
| 8.-12. | (B → ℓ) → (A → ℓ) | (7., a); ch.in$_2$), |
| 12a. | ¬ B → ¬ A | (**def**$_\neg$ ; 12.). |

The rule Contrap is applied only with the conclusion 231.

Contraposition 2.

(contrap): A → B → (¬ B → ¬ A)

    a) A → B

| 1.-7. | (B → ℓ) → ((A → B) → (A → ℓ)) | (ch.in.-fla$_2$), |
| 8.-12. | (B → ℓ) → (A → ℓ) | (7., a); ch.in$_2$), |

12a. ¬ B    (A    ⨏ )                                      [Def¬; 12.]

12b. ¬ B    ¬ A                                            [Def¬; 12a.].

and

1.-7.    (B    ⨏)    ((A    B)    (A    ⨏))                (ch.in.-fla₂),

8.-15.   (A    B)    ((B    ⨏)    (A    ⨏))                (7.; Int-Ant),

15a. (A    B)    (¬ B    (A    ⨏))                         [Def¬; 15.],

15a. (A    B)    (¬ B    ¬ A)                              [Def¬; 15a.].

Instances of contrap's conclusions occur at formulae 85.; 211.; 231. and 342..

∨-elimination (modus tollens ponens).

(Mtp₂):        $\dfrac{A \vee B, \quad \neg A}{B}$

The rule of Modus tollendo ponens has two forms. The formulae 113., 256. and 313. are the conclusions of the rule

(Mtp₂):        $\dfrac{A \vee B}{\neg A \quad B}$

which has the following derivation:

1.    ⨏    B                                               (Imp3),

2.        1.    3.                                         (Imp1),

3.    A    (⨏    B)                                        (1., 2.; MP),

4..       3.    5.                                         (Imp2),

5.    (A    ⨏    (A    B)                                  (3., 4.; MP),

5a.  ¬A  (A  B)                    [Def¬; 5.],

6.-13.  (A  (¬A  B)                (5a.; Int-Ant),

14.  B  (¬A  B)                    (Imp1),

15..      13.  (14.  17.)          (Con3),

16..      14.  17.                 (13., 15.; MP),

17.       13. & 14.                (14., 16.; MP),

18.       17.  19.                 (Dis3),

19.  A ∨ B  (¬A  B)                (17., 18.; MP),

20.  ¬A  B                         (a., 19.; MP).

(The other form of Mtp is

(Mtp$_1$):        $\underline{A \lor B, \quad \neg A}$
                           B

and it is, so far, not used here. Its derivation has 17 lines:

a)  A ∨ B;  b) ¬A

1.-5a.  ¬A  (A  B)                 (as in Mtp$_2$),

6.  A  B                           ( b), 5a.; MP),

7.-11.  B  B                       (Imp0),

12..      6.  (11.  14.)           (Con3),

13..      11.  14.                 (6., 12.; MP),

14.       6. & 11.                 (11., 13.; MP),

15.                14.    17.                          (Dis3),

16.  A ∨ B       B)                          (14., 15.; MP),

17.  B                                       ( a), 16.; MP).)

# APPENDIX B.

We supply, below, comments for each axiom justification marked [...], above. We use the notation $[\ ]_i$ to refer to the $i^{\text{th}}$ formula to which the axiom justification is attached.

$[\ ]_{42}$;   $sg(r) = 0 \quad sg(r) = 1$ is presupposed as an axiom,

$[\ ]_{44}$;   $r \cdot 0 = 0$ is presupposed as an axiom,

$[\ ]_{120}$;   $r \cdot 1 = r$ is presupposed as an axiom,

$[\ ]_{143}$;   $\overline{sg}(r) = 0 \quad \overline{sg}(r) = 1$ is presupposed as an axiom,

$[\ ]_{145}$;   as $[\ ]_{44}$,

$[\ ]_{170}$;   as $[\ ]_{44}$,

$[\ ]_{262}$;   $\overline{sg}(msd(r, s)) = 1 \quad r = s$ is presupposed as an axiom,

$[\ ]_{274}$;   as $[\ ]_{120}$,

$[\ ]_{293}$;   as $[\ ]_{143}$,

$[\ ]_{356}$;   as $[\ ]_{262}$,

$[\ ]_{377}$;   as $[\ ]_{120}$,

$[\ ]_{420}$; — from the definitions of $ax_{Ari}$, $\mathcal{b}_1$ and cfor,

$[\ ]_{421}$; — from the definitions of cfor and $\mathcal{l}$,

$[\ ]_{427}$; — from the definition of $\mathcal{l}$.

# APPENDIX C.

**C2.** We remind the reader that we are working in the *open* Ari although we have provided no alternative notations for this fragment of Ari. With this reminder out of the way we proceed with a presentation of the open Ari proof of the mp-induction schema.

**Remark**. Formula 57. is a schema for Ari's *course-of-values induction* and provable, in Ari, from the ordinary induction schema (cf. [Kleene, 1951], §40) which, for now, we leave postulated. Notice that *assumption* **d**., below, coincides with the first antecedent of the course-of-values induction formula (formula 162b. in [Kleene, 1952], §40). The schema of **mp**-induction is the formula

**mp**-ind ;    $\mu \ (Ant_1^{b_1}(\mu) \quad C(\mu))$

   $(\ \mu \ (Ant_2^{b_1}(\mu) \quad (C((\mu)_1) \quad (C((\mu)_2) \quad C(\mu)))$

   $(\ b_1(\mu) = 0 \quad C(\mu)))$,

The instance of **mp**-induction at formula 631., above, uses the formula $F_{(\mu)_0}$ for $C(\mu)$.

Proof:

   *Assume*:

**a.**   $\mu \ (Ant_1^{b_1}(\mu) \quad C(\mu))$ ,

**b.**   $\mu \ (Ant_2^{b_1}(\mu) \quad (C((\mu)_1) \quad (C((\mu)_2) \quad C(\mu))))$,

**c.**   $b_1(\mu) = 0$

   *and*

**d.**   $(\ < \mu \quad (b_1(\ ) = 0 \quad C(\ )))$

Then (in open Ari);

1.    $\mu \ (Ant_3^{b_1}(\mu) \quad b_1(\mu) = 1)$                                   (**df** $b_1^{\ 1}_3$),

2.        1.    3.                                                                                        (WBA ),

| | | | |
|---|---|---|---|
| 3. | $\text{Ant}_3^{b_1}(\mu) \quad b_1(\mu) = 1$ | | (1., 2.; MP), |
| 4. | $(b_1(\mu) = 0 \quad (b_1(\mu) = 1 \quad 0 = 1)$ | | (LEA$_1^=$), |
| 5. | $b_1(\mu) = 1 \quad 0 = 1$ | | (**c.**, 4.; MP), |
| 5a. | $b_1(\mu) = 1 \quad \ell$ | | [**def**$_\ell$; 5.], |
| 6.-10. | $\text{Ant}_3^{b_1}(\mu) \quad \ell$ | | (3., 5a.; ch.in.). |
| 10a. | $\neg\,(\text{Ant}_1^{b_1}(\mu) \quad \text{Ant}_2^{b_1}(\mu)) \quad \ell$ | | [**def**$_{\text{Ant}_3^{b_1}}$; 10.], |
| 10b. | $\neg\,\neg\,(\text{Ant}_1^{b_1}(\mu) \quad \text{Ant}_2^{b_1}(\mu))$ | | [**def**$_\neg$], |
| 11. | 10b. 12. | | (DNE), |
| 12. | $\text{Ant}_1^{b_1}(\mu) \quad \text{Ant}_2^{b_1}(\mu)$ | | (10b., 11.; MP), |
| 13. | **a**. 14. | | (WBA ), |
| 14. | $\text{Ant}_1^{b_1}(\mu) \quad C(\mu)$ | | (**a.**, 13.; MP), |
| 15. | $\mu = \text{mp}((\mu)_1, (\mu)_2)\ \&\ \neg\,\mu = 0\ \&\ b_1((\mu)_1) = 0\ \&\ b_1((\mu)_2) = 0$ | | |
| | $b_1((\mu)_2) = 0$ | | (Con 2), |
| 15a. | $\text{Ant}_2^{b_1}(\mu) \quad b_1((\mu)_2) = 0$ | | [**def**$_{\text{Ant}_2^{b_1}}$; 15.], |
| 16. | $\mu = \text{mp}((\mu)_1, (\mu)_2)\ \&\ \neg\,\mu = 0\ \&\ b_1((\mu)_1) = 0\ \&\ b_1((\mu)_2) = 0$ | | |
| | $\mu = \text{mp}((\mu)_1, (\mu)_2)\ \&\ \neg\,\mu = 0\ \&\ b_1((\mu)_1) = 0$ | | (Con 1), |
| 16a. | $\text{Ant}_2^{b_1}(\mu) \quad \mu = \text{mp}((\mu)_1, (\mu)_2)\ \&\ \neg\,\mu = 0\ \&\ b_1((\mu)_1) = 0$ | | [**def**$_{\text{Ant}_2^{b_1}}$; 16.], |
| 17. | $\mu = \text{mp}((\mu)_1, (\mu)_2)\ \&\ \neg\,\mu = 0\ \&\ b_1((\mu)_1) = 0 \quad b_1((\mu)_1) = 0$ | | (Con 2), |
| 18.-22. | $\text{Ant}_2^{\ell_1}(\mu)$ | $\ell_1((\mu)_1) = 0$ | (16a., 17.; ch.in.), |
| 23. | $\mu = mp((\mu)_1, (\mu)_2)$ & $\neg \mu = 0$ & $\ell_1((\mu)_1) = 0$ | | |
| | $\mu = mp((\mu)_1, (\mu)_2)$ & $\neg \mu = 0$ | | (Con 1), |
| 24.-28. | $\text{Ant}_2^{\ell_1}(\mu)$ | $\mu = mp((\mu)_1, (\mu)_2)$ & $\neg \mu = 0$ | (16a., 23.; ch. in.), |
| 29. | $\mu = mp((\mu)_1, (\mu)_2)$ & $\neg \mu = 0$ | $\mu = mp((\mu)_1, (\mu)_2)$ | (Con 1), |
| 30.-34. | $\text{Ant}_2^{\ell_1}(\mu)$ | $\mu = mp((\mu)_1, (\mu)_2)$ | (28., 29.; ch.in.), |
| 35. | $\mu = mp((\mu)_1, (\mu)_2)$ & $\neg \mu = 0$ | $\neg \mu = 0$ | (Con 2), |
| 36.-40. | $\text{Ant}_2^{\ell_1}(\mu)$ | $\neg \mu = 0$ | (28., 35; ch.in.), |
| 41. | $\mu (\neg \mu = 0 \quad (\mu)_1 < \mu)$ | | [...], |
| 42. | $\mu (\neg \mu = 0 \quad (\mu)_2 < \mu)$ | | [...], |
| 43. | 41. 44. | | (WBA ), |
| 44. | $\neg \mu = 0 \quad (\mu)_1 < \mu$ | | (41., 43.; MP), |
| 45.-49. | $\text{Ant}_2^{\ell_1}(\mu)$ | $(\mu)_1 < \mu$ | (40., 44; ch.in), |
| 50. | 42. 51. | | (WBA ), |
| 51. | $\neg \mu = 0 \quad (\mu)_2 < \mu$ | | (42., 50.; MP), |
| 52.-56. | $\text{Ant}_2^{\ell_1}(\mu)$ | $(\mu)_2 < \mu$ | (40., 51.; ch.in), |
| 57. | $\mu ( \quad ( \quad <\mu \quad (\ell_1( ) = 0 \quad C( ))) \quad (\ell_1(\mu) = 0 \quad C(\mu)))$ | | |
| | $\mu (\ell_1(\mu) = 0 \quad C(\mu))$ | | (course-of-values induction), |
| 58. | $\mu (\ell_1(\mu) = 0 \quad C(\mu)) \quad (\ell_1(\mu) = 0 \quad C(\mu))$ | | (WBA ), |
| 59.-63. | μ ( ( < μ (𝔟₁( ) = 0 C( ))) (𝔟₁(μ) = 0 C(μ))) | | |
| | (𝔟₁(μ) = 0 C(μ)) | (57., 58.; ch.in.), | |
| (μ)₁ — 64. | **d**. 65. | (SBA ), | |
| 65. | (μ )₁ < μ (𝔟₁((μ )₁) = 0 C((μ)₁) | (**d**., 64.; MP), | |
| (μ)₂ — 66. | **d**. 67. | (SBA ), | |
| 67. | (μ)₂ < μ (𝔟₁((μ )₂) = 0 C((μ)₂) | (**d**., 66.; MP), | |
| 68.-72. | Ant₂^𝔟₁(μ) (𝔟₁((μ )₁) = 0 C((μ)₁)) | (49., 65.; ch.in.), | |
| 73. | 72. (22. 75.) | (Imp 2), | |
| 74. | 22. 75. | (72., 73; MP), | |
| 75. | Ant₂^𝔟₁(μ) C((μ)₁) | (22., 74; MP), | |
| 76.-80. | Ant₂^𝔟₁(μ) (𝔟₁((μ )₂) = 0 C((μ)₂) | (56., 67.; ch.in.), | |
| 81. | 80. (15a. 83.) | (Imp 2), | |
| 82. | 15a. 83. | (80., 81.; MP), | |
| 83. | Ant₂^𝔟₁(μ) C((μ)₂) | (15a., 82; MP), | |
| 84. | **b**. 85. | (WBA ), | |
| 85. | Ant₂^𝔟₁(μ) (C((μ)₁) ( C((μ)₂) C(μ))) | (**b**., 84.; MP), | |
| 86. | 85. (75. 88.) | (Imp 2), | |
| 87. | 75. 88. | (85., 86.; MP), | |
| 88. | Ant₂^𝔟₁(μ) (C((μ)₂) C(μ)) | (75., 87.; MP), | |

| | | | | | |
|---|---|---|---|---|---|
| 89. | | 88. | (83. → 91.) | | (Imp 2), |
| 90. | | | 83. → 91. | | (88., 89; MP), |
| 91. | | Ant$_2$$\ell_i$($\mu$) → C($\mu$) | | | (83., 90; MP), |
| 92. | | 14. → (91. → 94.) | | | (Con 3), |
| 93. | | 91. → 94. | | | (14., 92.; MP), |
| 94. | | 14. & 91. | | | (91., 93.; MP), |
| 95. | | 94. → (12. → 97.) | | | (Dis 3), |
| 96. | | 12. → 97. | | | (94., 95; MP), |
| 97. | C($\mu$) | | | | (12., 96.; MP), |

and the discharging of **c**. gives

    c1.   $\ell_1(\mu) = 0$ → C($\mu$)                                       (1.-97.; **c**. discharged);

the discharging of **d**. followed by **Gen** give

    g1.   ( < $\mu$ → ($\ell_1$( ) = 0 → C( ))) → ($\ell_1(\mu) = 0$ → C($\mu$))      (c1.; **d**. discharged),

    g2.   $\mu$ ( ( < $\mu$ → ($\ell_1$( ) = 0 → C( ))) → ($\ell_1(\mu) = 0$ → C($\mu$)))      (g1.; **Gen**)

That is, the antecedent, g2., of the implication 63. is, deducible, in Ari, from the closed assumptions **a**., **b**., Therefore, from these assumptions also the consequent

    f1.   $\ell_1(\mu) = 0$ → C($\mu$)

of 63. can be deduced in the *open* Ari, and the discharging of the assumptions **b**., **a**. (in this order) gives the open Ari-proof of the **mp**-ind. formula (End of proof).



# APPENDIX **D1**.



In this Appendix we shall indicate how to Ari-prove the schema

(t-ax): $\quad \text{ax}_{\text{Ari}}(v_E) = 0 \quad$ E.

In the proof, 1.-61., below, there occur explicitly only three non-logical axioms, *viz.* formulae 1., 17., and 26., while implicitly there occurs only those non-logical axioms which occur along the branches $br_{46.}$ and $br_{51.}$ growing from the formulae 46. and 51., respectively. In the case of $br_{46.}$ it consists of the elementary calculation of ex with the arguments 2 and $\pi_E$ where the notation $\pi_E$ refers to the notion of En-N-generator explained in FOOTNOTE A. This case is covered as a special case of Lemma 1 of Appendix D2 and this lemma is cited in the proof comment as the justification for formula 46.. In the case of $br_{51.}$ it consists of a proof of the numeroidal reflection principle (nrp) which uses the p.r. character of $b$ and the formula

$b\text{calc}_b$ : $\quad ( b(\text{calc}_b(\ ,\ ),\ \overset{*}{b}(\ ,\mu) = (\ ,\ , b(\ ,\ ))) = 0$

provides the Ari-provability of the implication

$b(v[2^{\pi_E}], v_E) = 1 \quad b(\text{calc}_b(2^{\pi_E}], \pi_E),\ \overset{*}{b}(\ ,\mu) = 1(2^{\pi_E}], v_E)) = 0$

which is required for the proof of this nrp. This case is covered as a special case of Lemma 2 of Appendix D2. Inspection of this special case shows that it is free of any Kleene notations and, hence, can be given in Ari. In order to draw attention to this distinction this special case of Lemma 2 shall, for the time being, be labeled as Lemma 2* in the justification of formula 51., below.



| | | | |
|---|---|---|---|
| 1. | $(2^{\nu_E})_0 = \nu_E$ | | $(Ga_0)$, |
| 2. | 1. 3. | | $(Sym_=)$, |
| 3. | $\nu_E = (2^{\nu_E})_0$ | | (1., 2; MP), |
| 4. | 3. 5. | | (LEA-rp), |
| 5. | $2^{\nu_E} = 2^{(2^{\nu_E})_0}$ | | (3., 4; MP), |
| 6. | 3. 7. | | (LEA-rp), |
| 7. | $ax_{Ari}(\nu_E) = ax_{Ari}((2^{\nu_E})_0)$ | | (3., 6.; MP), |
| 8. | 7. 9. | | $(LEA_1^=)$, |
| 9. | $ax_{Ari}(\nu_E) = 0 \quad ax_{Ari}((2^{\nu_E})_0) = 0$ | | (7., 8.; MP), |
| 10. | 5. 11. | | (Con3), |
| 11. | $ax_{Ari}((2^{\nu_E})_0) = 0 \quad 2^{\nu_E} = 2^{(2^{\nu_E})_0}$ & $ax_{Ari}((2^{\nu_E})_0) = 0$ | | (5., 10; MP), |
| 12. | 11. 13. | | (Imp1), |
| 13. | $ax_{Ari}(\nu_E) = 0 \quad$ 11. | | (11., 12.; MP), |
| 14. | 13. (9. 16.) | | (Imp2), |
| 15. | 9. 16. | | (13., 14.; MP), |
| 16. | $ax_{Ari}(\nu_E) = 0 \quad 2^{\nu_E} = 2^{(2^{\nu_E})_0}$ & $ax_{Ari}((2^{\nu_E})_0) = 0$ | | (9., 15; MP), |
| 16a. | $ax_{Ari}(\nu_E) = 0 \quad Ant_1^{\ell_1}(2^{\nu_E})$ | | [$\mathbf{Def}_{Ant_1^{\ell_1}}$; 75.], |
| 17. | $(Ant_1^{\ell_1}(\ ) \quad \ell_1(\ ) = 0)$ | | $(df_1^{\ell_1})$, |
| | $2^{\nu_E} \dashv$ 18. \quad 17. \quad 19. | | (SBA), |



| | | | |
|---|---|---|---|
| 19. | Ant$_1$ $\mathcal{b}_1(2^{\nu E})$   $\mathcal{b}_1(2^{\nu E}) = 0$ | | (17., 18; MP), |
| 20. | $\mathcal{b}_1(2^{\nu E}) = 0$   (1.   $\mathcal{b}_1(2^{\nu E}) = 0$ & $2^{\nu E} = \nu_E$ ) | | (Con3), |
| 21. | 11.   22. | | (Imp1), |
| 22. | $\mathcal{b}_1(2^{\nu E}) = 0$   11. | | (11., 21.; MP), |
| 23. | 20.   (22.   25.) | | (Imp2), |
| 24. | 22.   25. | | (20., 23.; MP), |
| 25. | $\mathcal{b}_1(2^{\nu E}) = 0$   $\mathcal{b}_1(2^{\nu E}) = 0$ & $2^{\nu E} = \nu_E$ ) | | (22., 24.; MP), |
| 26. | ( $\mathcal{b}_1(\ ) = 0$ & $(\ )_0 =$    $\mathcal{b}(\ ,\ ) = 0$) | | (def$_2^{\mathcal{b}}$:), |
| | $2^{\nu E} \dashv 27.$    26.    28. | | (SBA ), |
| 28. | ( $\mathcal{b}_1(2^{\nu E}) = 0$ & $(2^{\nu E})_0 =$    $\mathcal{b}(2^{\nu E},\ ) = 0$) | | (26., 27; MP), |
| | $\nu_E \dashv 29.$    28.    30. | | (SBA ), |
| 30. | $\mathcal{b}_1(2^{\nu E}) = 0$ & $(2^{\nu E})_0 = \nu_E$    $\mathcal{b}(2^{\nu E}, \nu_E) = 0$ | | (28., 29; MP), |
| 31. | 30.   32. | | (Imp1), |
| 32. | $\mathcal{b}_1(2^{\nu E}) = 0$   30. | | (30., 31.; MP), |
| 33. | 32.   (25.   35.) | | (Imp2), |
| 34. | 25.   35. | | (32., 33.; MP), |
| 35. | $\mathcal{b}_1(2^{\nu E}) = 0$   $\mathcal{b}(2^{\nu E}, \nu_E) = 0$ | | (25., 34; MP), |
| 36. | 35.   37. | | (Imp1), |
| 37. | Ant$_1$ $\mathcal{b}_1(2^{\nu E})$    35. | | (35., 36.; MP), |
| 38. | 37.   (19.   40.) | | (Imp2), |
| 39. | 19.   40. | | (37., 38; MP), |



| | | | | |
|---|---|---|---|---|
| 40. | Ant$_1^{\mathscr{b}}$(2$^{\nu}$E) | $\mathscr{b}$(2$^{\nu}$E, $\nu_E$) = 0 | | (19., 39; MP), |
| 41. | 40. | 42. | | (Imp1), |
| 42. | ax$_{Ari}$($\nu_E$) = 0 | 40. | | (40., 41.; MP), |
| 43. | 42. | (16a. | 45.) | (Imp2), |
| 44. | 16a. | 45. | | (42., 43: MP), |
| 45. | ax$_{Ari}$($\nu_E$) = 0 | $\mathscr{b}$(2$^{\nu}$E, $\nu_E$) = 0 | | (16a., 44.; MP), |
| 46. | ($\nu$2)$^{\nu\pi}$E = $\nu$[2$^{\pi}$E] | | | (Lemma 1), |
| 46a. | $\nu$2$^{\nu\pi}$E = $\nu$[2$^{\pi}$E] | | | ($\nu$2 is 2: 46), |
| 46b. | $\nu$2$^{\nu}$E = $\nu$[2$^{\pi}$E] | | | ($\nu_E$ is $\nu\pi_E$: 46a.), |
| 47. | 46b.. | 48. | | (LEA-rp), |
| 48. | $\mathscr{b}$(2$^{\nu}$E, $\nu_E$) = $\mathscr{b}$($\nu$[2$^{\nu}$E, $\nu_E$]) | | | (46b., 47.: MP), |
| 49. | 48. | 50. | | (LEA$_{\bar 1}$), |
| 50. | $\mathscr{b}$(2$^{\nu}$E, $\nu_E$) = 0 | $\mathscr{b}$($\nu$[2$^{\nu}$E, $\nu_E$]) = 0 | | (48., 49.: MP), |
| 51. | $\mathscr{b}$($\nu$[2$^{\pi}$E], $\nu_E$) = 0 | E | | (Lemma 2*), |
| 52. | 51. | 53. | | (Imp1), |
| 53. | $\mathscr{b}$(2$^{\nu}$E, $\nu_E$) = 0 | 51. | | (51., 52: MP), |
| 54. | 53. | (50. | 56.) | (Imp2), |
| 55. | 50. | 56. | | (53., 54.: MP), |
| 56. | $\mathscr{b}$(2$^{\nu}$E, $\nu_E$) = 0 | E | | (50., 55: MP), |
| 57. | 56. | 58. | | (Imp1), |



| | | | | | |
|---|---|---|---|---|---|
| 58. | $ax_{Ari}(\nu_E) = 0$ | 56. | | | (56., 57.: MP), |
| 59. | | 58. (45. | 61.) | | (Imp2), |
| 60 | | 45. | 61. | | (58., 59.: MP), |
| 61. | $ax_{Ari}(\nu_E) = 0$ | E | | | (56., 57.: MP). |



APPENDIX **D2**.



In order to eliminate T-Ax as an axiom we give a sketch of its consequence, formula 7. which shall depend on two lemmas. Lemma 1 is a form of Gödel's Theorem V. in [GÖDEL, 1931] while Lemma 2 shall be a form of the proof-truth implication restricted to Ari 's axioms and quite unproblematic.

Although neither of these two lemmas pose any proof theoretical problem as such, their full proofs cannot be located in the published litterature and, thus cannot simply be refered to. A reason for their absence in the litterature is probably a function of the lengths of their (full) proofs each of which exceeds the length of more than 700 hundred lines. Therefore, it has been decided that these proofs will be more conveniently given as Part II of the present paper.

However, as only a special case for each of the two lemmas is needed for the proof of this Appendix, we shall refrain from formulating them in the full generality.

LEMMA 1. For any arguments of the function $\text{ex}(\_, \_)$ there exists a proof, $\text{Calc}_{\text{ex}}$, in Ari, the root of which is the equality $\text{ex}(\nu\_, \nu\_) = \nu\text{ex}(\_, \_)$. In particular, for the arguments 2, $\pi_{F_{(\mu)_0}}$ of $\text{ex}(\_, \_)$ there exists a proof $\text{Calc}_{\text{ex}}(2, \pi_{F_{(\mu)_0}})$.

PROOF. See Part II of this paper.

**Remark**. This lemma is needed for formula 110., below. (See Footnote A for the notations $\text{ex}$ and $\pi_{F_{(\mu)_0}}$).

LEMMA 2. .The formula $\mathcal{b}(\nu[2^{\pi_{F_{(\mu)_0}}}], \mathbf{v}_{F_{(\mu)_0}}) = 0 \quad F_{(\mu)_0}$ is Ari - provable as an instance of a proof-schema applicable to Ari - proofs of the proof-truth implication $\mathcal{b}(\mu, \mathbf{v}_{F_{(\mu)_0}}) = 0 \quad F_{(\mu)_0}$.

PROOF. See Part II of this paper.

**Remark**. Observe that the validity of the formulation of this lemma is restricted to the condition

$$\text{ax}_{\text{Ari}}(\pi_{F_{(\mu)_0}}) = 0. \qquad \pi_{F_{(\mu)_0}} = \nu\pi_{F_{(\mu)_0}}$$

The lemma is needed for formula 115., below. (See Footnote A for the notation $\pi_{F_{(\mu)_0}}$)

1.      $\mu = 2^{(\mu)_0}$ & $ax_{Ari}((\mu)_0) = 0$    $ax_{Ari}((\mu)_0) = 0$             (Con2),

     1a.   $Ant_1^{\ell_1}(\mu)$    $ax_{Ari}((\mu)_0) = 0$             $[\mathbf{Def}_{Ant_1^{\ell_1}}; 1]$,

2.      ( $ax_{Ari}(\ ) = 0$    $cfor(\ ) = 0$)             [...],

     $(\mu)_0 \dashv$ 3.    2.    4.             (SBA ),

4.      $ax_{Ari}((\mu)_0) = 0$    $cfor((\mu)_0) = 0$             (2., 3.; MP),

5.      4.    6.             (Imp1),

6.      $Ant_1^{\ell_1}(\mu)$    4.             (4., 5.; MP),

7.      6.    (1a.    9.)             (Imp2),

8.      1a.    9.             (6., 7.; MP),

9.      $Ant_1^{\ell_1}(\mu)$    $cfor((\mu)_0) = 0$             (1a., 8; MP),

10.      ( $cfor(\ ) = 0$    $fl(\ ) =$ )             $(\mathbf{def}_1^{fl})$,

     $(\mu)_0 \dashv$ 11.    10.    12.             (SBA ),

12.      $cfor((\mu)_0) = 0$    $fl((\mu)_0) = (\mu)_0$             (10., 11.; MP),

13.      12.    14.             (Imp1),

14.      $Ant_1^{\ell_1}(\mu)$    12.             (12., 13.; MP),

15.      14.    (9.    17.)             (Imp2),

16.      9.    17.             (14., 15; MP),

17.      $Ant_1^{\ell_1}(\mu)$    $fl((\mu)_0) = (\mu)_0$             (9., 16.; MP),

18.      $fl((\mu)_0) = (\mu)_0$    $v_{F_{fl((\mu)_0)}} = v_{F_{(\mu)_0}}$             $(LEA_2^v)$,

19.      18.    20.             (Imp1),

20.      $Ant_1^{\ell_1}(\mu)$    18.             (18., 19.; MP),

| | | | |
|---|---|---|---|
| 21. | 20. (17. 23.) | | (Imp2), |
| 22.. | 17. 23. | | (20., 21; MP), |
| 23.. | $Ant_1^{\ell_1}(\mu)$ $\nu_{F_{fl((\mu)_0)}} = \nu_{F_{(\mu)_0}}$ | | (17., 22.; MP), |
| 24. | $fl((\mu)_0) = fl((\mu)_0)$ | | (Ref), |
| 24a. | $\nu_{F_{fl((\mu)_0)}} = fl((\mu)_0)$ | | $[\mathbf{Def}\nu_{F_{fl((\mu)_0)}}]$, |
| 25. | $\nu_{F_{fl((\mu)_0)}} = \nu_{F_{(\mu)_0}}$ (24a. $\nu_{F_{(\mu)_0}} = fl((\mu)_0)$) | | $(LEA_{\overline{1}})$, |
| 26. | 25. 27. | | (Imp1), |
| 27. | $Ant_1^{\ell_1}(\mu)$ 25. | | (25., 26; MP), |
| 28. . | 27. (23. 30.) | | (Imp2), |
| 29. | 23. 30. | | (27., 28; MP), |
| 30. | $Ant_1^{\ell_1}(\mu)$ (24a. $\nu_{F_{(\mu)_0}} = fl((\mu)_0)$) | | (23., 29; MP), |
| 31. | 24a.. 32 | | (Imp1), |
| 32.. | $Ant_1^{\ell_1}(\mu)$ 24a. | | (24a., 31; MP), |
| 33.. | 30. (32. 35.) | | (Imp2), |
| 34. | 32. 35. | | (30, 33; MP), |
| 35.. | $Ant_1^{\ell_1}(\mu)$ $\nu_{F_{(\mu)_0}} = fl((\mu)_0))$ | | (23., 29; MP), |
| 36. | $fl((\mu)_0) = (\mu)_0$ $(\nu_{F_{(\mu)_0}} = fl((\mu)_0)$ $\nu_{F_{(\mu)_0}} = (\mu)_0)$ | | $(LEA_{\overline{2}})$, |
| 37.. | 36. 38. | | (Imp1), |
| 38.. | $Ant_1^{\ell_1}(\mu)$ 36. | | (36., 37.; MP), |
| 39.. | 38. (17. 41.) | | (Imp2), |
| 40. | 17. 41. | | (38, 39; MP), |

| | | | | |
|---|---|---|---|---|
| 41. | $\text{Ant}_1^{b_1}(\mu)$ | $(\nu_{F_{(\mu)_0}} = \text{fl}((\mu)_0)$ | $\nu_{F_{(\mu)_0}} = (\mu)_0)$ | (17., 40; MP), |
| 42.. | | 41. (35. 44.) | | (Imp2), |
| 43. | | 35. 44. | | (41, 42; MP), |
| 44.. | $\text{Ant}_1^{b_1}(\mu)$ | $\nu_{F_{(\mu)_0}} = (\mu)_0$ | | (35., 43; MP), |
| 45.. | $\nu_{F_{(\mu)_0}} = (\mu)_0$ | $\text{ax}_{\text{Ari}}(\nu_{F_{(\mu)_0}}) = \text{ax}_{\text{Ari}}((\mu)_0)$ | | (LEA-rp), |
| 46. | | 45. 47. | | (Imp1), |
| 47. | | $\text{Ant}_1^{b_1}(\mu)$ 45. | | (45., 46; MP), |
| 48. | | 47. (44. 50.) | | (Imp2), |
| 49. | | 44. 50. | | (47., 48; MP), |
| 50. | $\text{Ant}_1^{b_1}(\mu)$ | $\text{ax}_{\text{Ari}}(\nu_{F_{(\mu)_0}}) = \text{ax}_{\text{Ari}}((\mu)_0)$ | | (44., 49; MP), |
| 51. | $\text{ax}_{\text{Ari}}((\mu)_0) = 0$ | $(\text{ax}_{\text{Ari}}(\nu_{F_{(\mu)_0}}) = \text{ax}_{\text{Ari}}((\mu)_0)$ | $\text{ax}_{\text{Ari}}(\nu_{F_{(\mu)_0}}) = 0)$ | (LEA$_{\overline{2}}$), |
| 52. | | 51. 53. | | (Imp1), |
| 53. | | $\text{Ant}_1^{b_1}(\mu)$ 51. | | (51., 52; MP), |
| 54. | | 53. (1a. 56.) | | (Imp2), |
| 55. | | 1a. 56. | | (53., 54; MP), |
| 56. | $\text{Ant}_1^{b_1}(\mu)$ | $(\text{ax}_{\text{Ari}}(\nu_{F_{(\mu)_0}}) = \text{ax}_{\text{Ari}}((\mu)_0)$ | $\text{ax}_{\text{Ari}}(\nu_{F_{(\mu)_0}}) = 0)$ | (1a., 55.; MP), |
| 57. | | 56. (50. 59.) | | (Imp2), |
| 58. | | 50. 59. | | (56., 57; MP), |
| 59. | $\text{Ant}_1^{b_1}(\mu)$ | $\text{ax}_{\text{Ari}}(\nu_{F_{(\mu)_0}}) = 0$ | | (50., 58; MP), |
| 60· | $(2^{\nu_{F_{(\mu)_0}}})_0 = \nu_{F_{(\mu)_0}}$ | | | (Ga$_0$), |

| | | |
|---|---|---|
| 61. | 60.   62. | (Sym$_=$), |
| 62. | $\nu_{F(\mu)_0} = (2^{\nu}F_{(\mu)_0})_0$ | (60., 61; MP), |
| 63. | 62.   64. | (LEA-rp), |
| 64. | $2^{\nu}F_{(\mu)_0} = 2^{(2^{\nu}F_{(\mu)_0})_0}$ | (62., 63; MP), |
| 65. | 62.   66. | (LEA-rp), |
| 66. | $ax_{Ari}(\nu_{F_{(\mu)_0}}) = ax_{Ari}((2^{\nu}F_{(\mu)_0})_0)$ | (62., 65.; MP), |
| 67. | 60.   62. | (LEA$_1^=$), |
| 68. | $ax_{Ari}(\nu_{F_{(\mu)_0}}) = 0 \quad ax_{Ari}((2^{\nu}F_{(\mu)_0})_0) = 0$ | (66., 67.; MP), |
| 69. | 68.   70. | (Imp1), |
| 70. | $Ant_1^{\ell_1}(\mu)$   68. | (68., 69; MP), |
| 71. | 70.   (59.   73.) | (Imp2), |
| 72. | 59.   73. | (70., 71; MP), |
| 73. | $Ant_1^{\ell_1}(\mu) \quad ax_{Ari}((2^{\nu}F(\mu)_0)_0) = 0$ | (59., 72.; MP), |
| 74. | 64.   75. | (Con3), |
| 75. | $ax_{Ari}((2^{\nu}F_{(\mu)_0})_0) = 0 \quad 2^{\nu}F_{(\mu)_0} = 2^{(2^{\nu}F_{(\mu)_0})_0} \ \& \ ax_{Ari}((2^{\nu}F_{(\mu)_0})_0) = 0$ | (64., 74; MP), |
| 75a. | $ax_{Ari}((2^{\nu}F_{(\mu)_0})_0) = 0 \quad Ant_1^{\ell_1}(2^{\nu}F_{(\mu)_0})$ | [**Def**$_{Ant_1^{\ell_1}}$; 75.], |
| 76. | 75a.   77. | (Imp1), |
| 77. | $Ant_1^{\ell_1}(\mu)$   75a. | (75a., 76.; MP), |
| 78. | 77.   (73.   80.) | (Imp2), |
| 79. | 73.   80. | (77., 78; MP), |

| | | | |
|---|---|---|---|
| 80. | $\text{Ant}_1^{\mathscr{b}_1}(\mu)$ $\text{Ant}_1^{\mathscr{b}_1}(2^{\nu_F(\mu)_0})$ | | (73., 79.; MP), |
| 81. | $(\text{Ant}_1^{\mathscr{b}_1}(\ )\ \mathscr{b}_1(\ ) = 0)$ | | $(\text{df}_1^{\mathscr{b}_1})$, |
| | $2^{\nu_F(\mu)_0}$ ⊣ 82. 81. 83. | | (SBA ), |
| 83. | $\text{Ant}_1^{\mathscr{b}_1}(2^{\nu_F(\mu)_0})$ $\mathscr{b}_1(2^{\nu_F(\mu)_0}) = 0$ | | (81., 82; MP), |
| 84. | 83. 85. | | (Imp1), |
| 85. | $\text{Ant}_1^{\mathscr{b}_1}(\mu)$ 83. | | (83., 84.; MP), |
| 86. | 85. (80. 88.) | | (Imp2), |
| 87. | 80. 88. | | (85., 86.; MP), |
| 88. | $\text{Ant}_1^{\mathscr{b}_1}(\mu)$ $\mathscr{b}_1(2^{\nu_F(\mu)_0}) = 0$ | | (80., 87; MP), |
| 89. | $(\mathscr{b}_1(\ ) = 0\ \&\ (\ )_0 =\ \mathscr{b}(\ ,\ ) = 0)$ | | $(\text{def}_2^{\mathscr{b}}:)$, |
| | $2^{\nu_F(\mu)_0}$ ⊣ 90. 89. 91. | | (SBA ), , |
| 91. | $(\mathscr{b}_1(2^{\nu_F(\mu)_0}) = 0\ \&\ (2^{\nu_F(\mu)_0})_0 =\ \mathscr{b}(2^{\nu_F(\mu)_0},\ ) = 0)$ | | (89., 90; MP), |
| | $\nu_{F(\mu)_0}$ ⊣ 92. 91. 93. | | (SBA ), |
| 93. | $\mathscr{b}_1(2^{\nu_F(\mu)_0}) = 0\ \&\ (2^{\nu_F(\mu)_0})_0 =\ \nu_{F(\mu)_0}\ \mathscr{b}(2^{\nu_F(\mu)_0}, \nu_{F(\mu)_0}) = 0$ | | (91., 92; MP), |
| 94. | $\mathscr{b}_1(2^{\nu_F(\mu)_0}) = 0$ $((2^{\nu_F(\mu)_0})_0 =\ \nu_{F(\mu)_0}$ | | |
| | $\mathscr{b}_1(2^{\nu_F(\mu)_0}) = 0\ \&\ (2^{\nu_F(\mu)_0})_0 =\ \nu_{F(\mu)_0})$ | | (Con3), |
| 95. | 94. 96. | | (Imp1), |
| 96. | $\text{Ant}_1^{\mathscr{b}_1}(\mu)$ 94. | | (94., 95.; MP), |
| 97. | | 96. (88. 99.) | (Imp2), |
| 98. | | 88. 99. | (96., 97.; MP), |
| 99. | $\text{Ant}_1^{\ell_1(\mu)}$ | $((2^{\nu}F_{(\mu)_0})_0 = \nu_{F_{(\mu)_0}}$ | |
| | | $\ell_1(2^{\nu}F_{(\mu)_0}) = 0 \,\&\, (2^{\nu}F_{(\mu)_0})_0 = \nu_{F_{(\mu)_0}})$ | (88., 98.; MP), |
| 100. | | 60. 101. | (Imp1), |
| 101. | $\text{Ant}_1^{\ell_1(\mu)}$ | $(2^{\nu}F_{(\mu)_0})_0 = \nu_{F_{(\mu)_0}}$ | (60., 100.; MP), |
| 102. | | 99. (101. 104..) | (Imp2), |
| 103. | | 101. 104. | (99., 102.; MP), |
| 104. | $\text{Ant}_1^{\ell_1(\mu)}$ | $\ell_1(2^{\nu}F_{(\mu)_0}) = 0 \,\&\, (2^{\nu}F_{(\mu)_0})_0 = \nu_{F_{(\mu)_0}}$ | (101., 103.; MP), |
| 105. | | 93. 106. | (Imp1), |
| 106. | | $\text{Ant}_1^{\ell_1(\mu)}$ 93. | (93., 105.; MP), |
| 107. | | 106. (104. 109..) | (Imp2), |
| 108. | | 104. 109. | (106., 107.; MP), |
| 109. | $\text{Ant}_1^{\ell_1(\mu)}$ | $\ell(2^{\nu}F_{(\mu)_0}, \nu_{F_{(\mu)_0}}) = 0$ | (104., 108.; MP), |
| 109a. | $\text{Ant}_1^{\ell_1(\mu)}$ | $\ell(2^{\nu\pi}F_{(\mu)_0}, \nu_{F_{(\mu)_0}}) = 0$ | ($\nu_{F_{(\mu)_0}}$ is $\nu\pi_{F_{(\mu)_0}}$; 109.), |
| 109b. | $\text{Ant}_1^{\ell_1(\mu)}$ | $\ell([\nu 2]^{\nu\pi}F_{(\mu)_0}, \nu_{F_{(\mu)_0}}) = 0$ | ($\nu 2$ is 2; 109a.), |
| 110. | $[\nu 2]^{\nu\pi}F_{(\mu)_0} = \nu[2^{\pi}F_{(\mu)_0}]$ | | (Lemma 1), |
| 111. | | 110. 112. | (LEA¨rp), |
| 112. | $\ell([\nu 2]^{\nu\pi}F_{(\mu)_0}, \nu_{F_{(\mu)_0}}) = \ell([\nu[2^{\pi}F_{(\mu)_0}], \nu_{F_{(\mu)_0}})$ | | (110., 111.; MP), |
| 113. | | 112. 113. | (LEA$\bar{\bar{1}}$), |
| 114. | $\mathscr{E}([\nu 2]^{\nu \pi} F_{(\mu)_0}, \nu_{F_{(\mu)_0}}) = 0 \quad \mathscr{E}(\nu[2^{\pi} F_{(\mu)_0}], \nu_{F_{(\mu)_0}}) = 0$ | | (112., 113.; MP), |
| 115. | $\mathscr{E}(\nu[2^{\pi} F_{(\mu)_0}], \nu_{F_{(\mu)_0}}) = 0 \quad F_{(\mu)_0}$ | | (Lemma 2), |
| 116. | 115. 117. | | (Imp1), |
| 117. | $\mathscr{E}([\nu 2]^{\nu \pi} F_{(\mu)_0}, \nu_{F_{(\mu)_0}}) = 0 \quad$ 115. | | (115., 116.; MP), |
| 118. | 117. (114. 120..) | | (Imp2), |
| 119. | 114. 120. | | (117., 118.; MP), |
| 120. | $\mathscr{E}([\nu 2]^{\nu \pi} F_{(\mu)_0}, \nu_{F_{(\mu)_0}}) = 0 \quad F_{(\mu)_0}$ . | | (114., 119.; MP), |
| 121. | 120. 122. | | (Imp1), |
| 122. | $\text{Ant}_1 \mathscr{E}_{1(\mu)}$   120. | | (120., 121.; MP), |
| 123. | 122. (109b. 125..) | | (Imp2), |
| 124. | 109b. 125. | | (122., 123.; MP), |
| 125. | $\text{Ant}_1 \mathscr{E}_{1(\mu)} \quad F_{(\mu)_0}$ | | (109b., 124.; MP). |

*connection with the Foundations of Mathematics)*, Constructive Mathematics, Proceedings, New Mexico, 1980, ed. F Richman, LNM 873, Springer Verlag, Berlin Heidelberg New York, 1981, pp. 274-313.

[YESSENIN-VOLPIN, 2000]: Yessenin-Volpin, Alexander S., *A Completeness Proof for Ari-like Systems* (in preparation).